% Dmitri Pavlov's Plain TeX macros
% Copyright 2017, 2018 Dmitri Pavlov

% Remove \outer
\edef\inewcount{\noexpand\csname newcount\endcsname}
\edef\inewdimen{\noexpand\csname newdimen\endcsname}
\edef\inewskip{\noexpand\csname newskip\endcsname}
\edef\inewmuskip{\noexpand\csname newmuskip\endcsname}
\edef\inewbox{\noexpand\csname newbox\endcsname}
\edef\inewhelp{\noexpand\csname newhelp\endcsname}
\edef\inewtoks{\noexpand\csname newtoks\endcsname}
\edef\inewread{\noexpand\csname newread\endcsname}
\edef\inewwrite{\noexpand\csname newwrite\endcsname}
\edef\inewfam{\noexpand\csname newfam\endcsname}
\edef\inewlanguage{\noexpand\csname newlanguage\endcsname}
\edef\inewinsert{\noexpand\csname newinsert\endcsname}
\edef\inewif{\noexpand\csname newif\endcsname}

% UTF-8 encoding for Plain TeX.
% Copyright 2008, 2015, 2018 Dmitri Pavlov.
% You may redistribute this file under the terms of GNU General Public License version 3.
% Report bugs and suggestions to me by email (host: math.berkeley.edu, user: pavlov).
% So far this file includes support for most of TeX symbols and Latin 1 letters.

% Setup catcodes for UTF-8 characters: legitimate initial octets are active, intermediate octets are letters, other octets are illegal.
\countdef\ch=253
\ch="80 \loop\ifnum\ch<"100 \lccode\ch=0 \uccode\ch=0 \advance\ch1 \repeat % for LaTeX
\ch="80 \loop\ifnum\ch<"C0 \catcode\ch=11 \advance\ch1 \repeat % make all intermediate UTF-8 octets (top 2 bits are 10) letters
\ch="C0 \loop\ifnum\ch<"100 \catcode\ch=15 \advance\ch1 \repeat % make all octets with the top 2 bits set illegal
\ch="C2 \loop\ifnum\ch<"F5 \catcode\ch=\active \advance\ch1 \repeat % make all legitimate initial UTF-8 octets active
\ch="C2 \loop\ifnum\ch<"E2 \uccode\ch=1 \advance\ch1 \repeat % make uccode=1 for all initial UTF-8 octets of letters, to allow macros to recognize Unicode letters
\let\ch\underfined

% \cdef\macro+ is like \def\macro{+}, but the catcode of + will be 11
\catcode0=12 \def\cdef#1#2{\begingroup\lccode0=`#2 \lowercase{\endgroup \def#1{^^@}}} \catcode0=9

% Preliminary macros for UTF-8: every legitimate UTF-8 sequence produces a control sequence whose name consists precisely of the octets in the sequence.
\catcode`@=\active
\def@#1{\cdef\chr#1 \edef#1##1{\noexpand\csname\chr##1\noexpand\endcsname}} % 2-byte
@^^c2@^^c3@^^c4@^^c5@^^c6@^^c7@^^c8@^^c9@^^ca@^^cb@^^cc@^^cd@^^ce@^^cf@^^d0@^^d1@^^d2@^^d3@^^d4@^^d5@^^d6@^^d7@^^d8@^^d9@^^da@^^db@^^dc@^^dd@^^de@^^df
\def@#1{\cdef\chr#1 \edef#1##1##2{\noexpand\csname\chr##1##2\noexpand\endcsname}} % 3-byte
@^^e0@^^e1@^^e2@^^e3@^^e4@^^e5@^^e6@^^e7@^^e8@^^e9@^^ea@^^eb@^^ec@^^ed@^^ee@^^ef
\def@#1{\cdef\chr#1 \edef#1##1##2##3{\noexpand\csname\chr##1##2##3\noexpand\endcsname}} % 4-byte
@^^f0@^^f1@^^f2@^^f3@^^f4
\let\chr\undefined
\let\cdef\undefined

% The following three lines establish a concise syntax for defining characters: each line contains a UTF-8 character followed by its definition.
\def\grabfuturelet{\futurelet\next\grabexamine}
\def\grabexamine{\ifx\next\csname\expandafter\grab\fi}
\obeylines \def\grab\csname#1\endcsname#2^^M{\expandafter\def\csname#1\endcsname{#2}\expandafter\grabfuturelet} \expandafter\grabfuturelet%
% Coverage: all TeX symbols except those not present in Unicode and the part of Latin 1 and Latin Extended-A that is present in the Computer Modern fonts.
 ~
¢{\hbox{\rm\rlap/c}}
£{\it\$}
%¤
%¥
%¦
%¨
%ª
«\leftguillemet
­\-
%®
%¯
%°
%²
%³
%´
%µ
%·\cdotp
%¸\c
%¹
%º
»\rightguillemet
%¼
%½
%¾
À{\`A}
Á{\'A}
Â{\^A}
Ã{\~A}
Ä{\"A}
Ç{\c C}
È{\`E}
É{\'E}
Ê{\^E}
Ë{\"E}
Ì{\`I}
Í{\'I}
Î{\^I}
Ï{\"I}
% must define \ETH, \THORN, \eth, \thorn separately
Ð\ETH
Ñ{\~N}
Ò{\`O}
Ó{\'O}
Ô{\^O}
Õ{\~O}
Ö{\"O}
Ù{\`U}
Ú{\'U}
Û{\^U}
Ü{\"U}
Ý{\'Y}
Þ\THORN
à{\`a}
á{\'a}
â{\^a}
ã{\~a}
ä{\"a}
ç{\c c}
è{\`e}
é{\'e}
ê{\^e}
ë{\"e}
ì{\`\i}
í{\'\i}
î{\^\i}
ï{\"\i}
ð\eth
ñ{\~n}
ò{\`o}
ó{\'o}
ô{\^o}
õ{\~o}
ö{\"o}
ù{\`u}
ú{\'u}
û{\^u}
ü{\"u}
ý{\'y}
þ\thorn
ÿ{\"y}
% Latin Extended-A
Ā{\=A}
ā{\=a}
Ă{\u A}
ă{\u a}
%Ą
%ą
Ć{\'C}
ć{\'c}
Ĉ{\^C}
ĉ{\^c}
Ċ{\.C}
ċ{\.c}
Č{\v C}
č{\v c}
Ď{\v D}
ď{\v d}
%Đ
%đ
Ē{\=E}
ē{\=e}
Ĕ{\u E}
ĕ{\u e}
Ė{\.E}
ė{\.e}
%Ę
%ę
Ě{\v E}
ě{\v e}
Ĝ{\^G}
ĝ{\^g}
Ğ{\u G}
ğ{\u g}
Ġ{\.G}
ġ{\.g}
Ģ{\c G}
ģ{\c g}
Ĥ{\^H}
ĥ{\^h}
%Ħ
%ħ
Ĩ{\~I}
ĩ{\~\i}
Ī{\=I}
ī{\=\i}
Ĭ{\u I}
ĭ{\u\i}
%Į
%į
İ{\.I}
%ı\i
ĲIJ
ĳij
Ĵ{\^J}
ĵ{\^\j}
Ķ{\c K}
ķ{\c k}
%ĸ
Ĺ{\'L}
ĺ{\'l}
Ļ{\c L}
ļ{\c l}
Ľ{\v L}
ľ{\v l}
%Ŀ
%ŀ
%Ł\L
%ł\l
Ń{\'N}
ń{\'n}
Ņ{\c N}
ņ{\c n}
Ň{\v N}
ň{\v n}
%ŉ
%Ŋ
%ŋ
Ō{\=O}
ō{\=o}
Ŏ{\u O}
ŏ{\u o}
Ő{\"O}
ő{\"o}
%Œ\OE
%œ\oe
Ŕ{\'R}
ŕ{\'r}
Ŗ{\c R}
ŗ{\c r}
Ř{\v R}
ř{\v r}
Ś{\'S}
ś{\'s}
Ŝ{\^S}
ŝ{\^s}
Ş{\c S}
ş{\c s}
Š{\v S}
š{\v s}
Ţ{\c T}
ţ{\c t}
Ť{\v T}
ť{\v t}
%Ŧ
%ŧ
Ũ{\~U}
ũ{\~u}
Ū{\=U}
ū{\=u}
Ŭ{\u U}
ŭ{\u u}
%Ů
%ů
Ű{\H U}
ű{\H u}
%Ų
%ų
Ŵ{\^W}
ŵ{\^w}
Ŷ{\^Y}
ŷ{\^y}
Ÿ{\"Y}
Ź{\'Z}
ź{\'z}
Ż{\.Z}
ż{\.z}
Ž{\v Z}
ž{\v z}
%ſ
’'
‘`
”{''}
“{``}
‐-
–{--}
—{---}
% these ligatures should not be used in TeX text
%ﬀ{ff}
%ﬁ{fi}
%ﬂ{fl}
%ﬃ{ffi}
%ﬄ{ffl}
¡{!`}
¿{?`}
−-
′'
% these ligatures should not be used in TeX text
%″{''}
%‴{'''}
%⁗{''''}
ß\ss
æ\ae
Æ\AE
œ\oe
Œ\OE
ø\o
Ø\O
å\aa
Å\AA
ł\l
Ł\L
% defined below
%ı\i
%ȷ\j
†\dag
‡\ddag
§\S
¶\P
% combining characters are left undefined because they come after, not before the accented character
%◌̣\d
%◌̱\b
%◌̧\c
©\copyright
…\dots
%◌̀\`
%◌́\'
%◌̌\v
%◌̆\u
%◌̄\=
%◌̂\^
%◌̇\.
%◌̋\H
%◌̃\~
%◌̈\"
%◌͡\t
%\brace[lr][du]
α\alpha
β\beta
γ\gamma
δ\delta
ϵ\epsilon
ζ\zeta
η\eta
θ\theta
ι\iota
κ\kappa
λ\lambda
μ\mu
ν\nu
ξ\xi
οo
π\pi
ρ\rho
σ\sigma
τ\tau
υ\upsilon
ϕ\phi
χ\chi
ψ\psi
ω\omega
ε\varepsilon
ϑ\vartheta
ϖ\varpi
ϱ\varrho
ς\varsigma
φ\varphi
Γ\Gamma
Δ\Delta
Θ\Theta
Λ\Lambda
Ξ\Xi
Π\Pi
Σ\Sigma
Υ\Upsilon
Φ\Phi
Ψ\Psi
Ω\Omega
ℵ\aleph
ℏ\hbar
ı\relax\ifmmode\imath\else\i\fi
ȷ\relax\ifmmode\jmath\else\j\fi
ℓ\ell
℘\wp
ℜ\Re
ℑ\Im
∂\partial
∞\infty
%'\prime
∅\emptyset
∇\nabla
√\surd
⊤\top
⊥\bot
∠\angle
△\triangle
∀\forall
∃\exists
¬\neg
♭\flat
♮\natural
♯\sharp
♣\clubsuit
♢\diamondsuit
♡\heartsuit
♠\spadesuit
∐\coprod
⋁\bigvee
⋀\bigwedge
⨄\biguplus
⋂\bigcap
⋃\bigcup
∫\int
%!\intop
∏\prod
∑\sum
⨂\bigotimes
⨁\bigoplus
⨀\bigodot
∮\oint
%!\ointop
⨆\bigsqcup
%!\smallint
◁\triangleleft
▷\triangleright
△\bigtriangleup
▽\bigtriangledown
∧\wedge
∨\vee
∩\cap
∪\cup
%‡\ddagger : \ddag
%†\dagger : \dag
⊓\sqcap
⊔\sqcup
⊎\uplus
⨿\amalg
⋄\diamond
∙\bullet
≀\wr
÷\div
⊙\odot
⊘\oslash
⊗\otimes
⊖\ominus
⊕\oplus
∓\mp
±\pm
∘\circ
%white circle:
○\Orb
%large circle:
◯\bigcirc
∖\setminus
⋅\cdot
∗\ast
% Latin 1
×\times
% Unicode
⨯\times
⋆\star
∝\propto
⊑\sqsubseteq
⊒\sqsupseteq
∥\parallel
‖\|
∣\divides
%|\mid
⊣\dashv
⊢\vdash
↗\nearrow
↘\searrow
↖\nwarrow
↙\swarrow
⇔\Leftrightarrow
⇐\Leftarrow
⇒\Rightarrow
≠\neq
≤\leq
≥\geq
≻\succ
≺\prec
≈\approx
≽\succeq
≼\preceq
⊃\supset
⊂\subset
⊇\supseteq
⊆\subseteq
∈\in
∋\ni
≫\gg
≪\ll
%◌̸\not
↔\leftrightarrow
←\leftarrow
→\rightarrow
↦\mapsto
%↦\mapstochar : \mapsto
∼\sim
≃\simeq
%⊥\perp : \bot
≡\equiv
≍\asymp
⌣\smile
⌢\frown
↼\leftharpoonup
↽\leftharpoondown
⇀\rightharpoonup
⇁\rightharpoondown
↪\hookrightarrow
↩\hookleftarrow
%\lhook
%\rhook
⋈\bowtie
⊨\models
⟹\Longrightarrow
⟶\longrightarrow
⟵\longleftarrow
⟸\Longleftarrow
⟼\longmapsto
⟷\longleftrightarrow
⟺\Longleftrightarrow
%⇔\iff : \Longleftrightarrow
%.\ldotp : .
%⋅\cdotp : \cdot
%:\colon : :
%…\ldots : \dots
⋯\cdots
⋮\vdots
⋱\ddots
%\acute, \grave, \ddot, \tilde, \bar, \breve, \check, \hat, \vec, \dot, \widetilde, \widehat
%\overrightarrow, \overleftarrow, \overbrace, \underbrace
%\lmoustache, \rmoustache, \lgroup, \rgroup, \arrowvert, \Arrowvert, \bracevert
∥\Vert
%|\vert
↑\uparrow
↓\downarrow
↕\updownarrow
⇑\Uparrow
⇓\Downarrow
⇕\Updownarrow
%\\backslash
⟩\rangle
⟨\langle
%{\lbrace
%}\rbrace
⌉\rceil
⌈\lceil
⌋\rfloor
⌊\lfloor
≅\cong
∉\notin
⇌\rightleftharpoons
≐\doteq
% Combined symbols from Unicode math blocks
∄\not\exists
∌\not\ni
∔\dot+
∕/
∣|
∤\not|
∦\not\|
∬\int\!\!\!\int
∭\int\!\!\!\int\!\!\!\int
∮\oint
∸\dot-
≁\not\sim
≄\not\simeq
≆\not\cong
≇\not\cong
≉\not\approx
≐\dot=
≔:=
≕=:
≢\not\equiv
≭\not\asump
≮\not<
≯\not<
≰\not\le
≱\not\ge
⊀\not\prec
⊁\not\succ
⊄\not\subset
⊅\not\supset
⊈\not\subseteq
⊉\not\supseteq
⊦\vdash
⊧\models
⊬\not\vdash
⊭\not\models
⊲\triangleleft
⊳\triangleright
⋠\not\preceq
⋡\not\succeq
⋤\not\sqsubseteq
⋥\not\sqsupseteq
⋪\not\triangleleft
⋫\not\triangleright
◻\square

\catcode`\^^M=5 %
\let\grabfuturelet\undefined \let\grabexamine\undefined \let\grab\undefined

% Generic macros for UTF-8: every legitimate UTF-8 sequence produces a control sequence whose name consists precisely of the octets in the sequence.
% An undefined control sequence produces an error message.
\let\xcsname=\csname
\let\xendcsname=\endcsname
\def@#1{\def#1##1{\expandafter\ifx\csname\string#1##1\endcsname\relax\errmessage{Undefined UTF-8 sequence \string#1##1}\else\xcsname\string#1##1\xendcsname\fi}}
@^^c2@^^c3@^^c4@^^c5@^^c6@^^c7@^^c8@^^c9@^^ca@^^cb@^^cc@^^cd@^^ce@^^cf@^^d0@^^d1@^^d2@^^d3@^^d4@^^d5@^^d6@^^d7@^^d8@^^d9@^^da@^^db@^^dc@^^dd@^^de@^^df
\def@#1{\def#1##1##2{\expandafter\ifx\csname\string#1##1##2\endcsname\relax\errmessage{Undefined UTF-8 sequence \string#1##1##2}\else\xcsname\string#1##1##2\xendcsname\fi}}
@^^e0@^^e1@^^e2@^^e3@^^e4@^^e5@^^e6@^^e7@^^e8@^^e9@^^ea@^^eb@^^ec@^^ed@^^ee@^^ef
\def@#1{\def#1##1##2##3{\expandafter\ifx\csname\string#1##1##2##3\endcsname\relax\errmessage{Undefined UTF-8 sequence \string#1##1##2##3}\else\xcsname\string#1##1##2##3\xendcsname\fi}}
@^^f0@^^f1@^^f2@^^f3@^^f4
\let@\undefined
\catcode`@=12

\newif\ifscroll % scroll vs codex
\newif\ifsuppressunusedbib % suppress unused bibliography items?

% Diagonostics
\tracinglostchars=2

\def\printerr#1{\immediate\write17{#1}}
\def\warningline#1#2{\printerr{! #2}\printerr{l.#1}\printerr{}}
\def\warning{\warningline{\the\inputlineno}}
%\let\warning\errmessage % turn warnings into errors

% Syntactic macros
\long\def\gobble#1{}
\ifx\gobbleinit\undefined{\long\gdef\gobbleinit#1\par{}}\fi
\def\expand#1{\edef\expandmacro{#1}\expandmacro\let\expandmacro\undefined}
\def\setetok#1#2{\expand{\noexpand#1{#2}}}
\def\expandtoks#1{\expandafter\edef\expandafter\expandmacro\expandafter{\the#1}#1\expandafter{\expandmacro}}
\def\appendexpand#1#2{\setetok#1{\the#1#2}}
\long\def\append#1#2{#1\expandafter{\the#1#2}}
\long\def\appendtoksexpand#1#2{#1\expandafter\expandafter\expandafter{\expandafter\the\expandafter#1\the#2}}
\long\def\appendonceexpand#1#2{#1\expandafter\expandafter\expandafter{\expandafter\the\expandafter#1#2}}

% Macros that use \special, known to work with xdvi, dvips, dvipdfm

%\let\printlink\print % print links and anchors
\def\link#1#2{\lhighlight{#2}}
\def\llink#1{\printlink{llink #1}\link{\ohash#1}}
\catcode`\#=11 \def\ohash{#}\catcode`\#=6
\catcode`\&=11 \def\ampersand{&}\catcode`\&=4
\def\anchor#1#2{\printlink{anchor #1 #2}#2}

\def\setpapersize#1#2{} % \number#1 sp triggers bugs in dvips and does not work with dvipdfm(x)
\def\dumpbox#1#2#3{\shipout\vbox{\setpapersize{#1}{#2}\unvbox#3}}
\def\mps#1{\epsfbox{#1}}
\def\metadata#1#2{}
\def\src{} % for proper match between bibliographic references and source code

% A simplified version of the epsfbox macro from epsf.tex, exclusively for METAPOST output, with HiResBoundingBox
% Normally should be invoked through \mps, which works also with PDF
\newread\epsffilein
\newif\ifepsfbbfound\inewif\ifepsffilecont
\newdimen\epsfxsize\inewdimen\epsfysize
\newdimen\pspoints\pspoints1bp
\let\runmp\errmessage % will be set to \warning if the METAPOST file cannot be compiled
\def\epsfbox#1{\openin\epsffilein=#1 \ifeof\epsffilein\runmp{Could not open file #1}\else
	{\def\do##1{\catcode`##1=12}\dospecials\catcode`\ =10\epsffileconttrue
		\epsfbbfoundfalse
		\loop\read\epsffilein to\epsffileline \ifeof\epsffilein\epsffilecontfalse\else\expandafter\epsfaux\epsffileline :. \\\fi\ifepsffilecont\repeat
		\ifepsfbbfound\else\errmessage{No HiResBoundingBox comment found in file #1}\fi}%
	\closein\epsffilein
	\epsfysize\epsfury\pspoints \advance\epsfysize-\epsflly\pspoints
	\epsfxsize\epsfurx\pspoints \advance\epsfxsize-\epsfllx\pspoints
	% We create a box with height and depth corresponding to the two vertical dimensions in the bounding box and width given by the total width
	\setbox0\hbox{\vbox to\epsfury\pspoints{\vfil\hbox to\epsfxsize{\dimen0=\epsfllx\pspoints \kern-\dimen0 \includegraphics{#1}\hfil}}}%
	\dp0=\epsflly\pspoints \dp0=-\dp0
	\box0 \fi}
%	\hbox{\vbox to\epsfysize{\vfil\hbox to\epsfxsize{\special{psfile=#1}\hfil}}}}
% llx=\epsfllx\space lly=\epsflly\space urx=\epsfurx\space ury=\epsfury
%	    rwi=\number\epsftmp\space rhi=\number\epsfrsize
\catcode`\%=12 \def\epsfbblit{%%HiResBoundingBox} \catcode`\%=14
\let\dummybrace=} % to compensate for { on the previous line when it is read by \read
\def\epsfaux#1:#2\\{\def\testit{#1}\ifx\testit\epsfbblit \epsfgrab #2 . . . \\\epsffilecontfalse\epsfbbfoundtrue\fi}
\def\empty{}
\def\epsfgrab #1 #2 #3 #4 #5\\{\gdef\epsfllx{#1}\ifx\epsfllx\empty\epsfgrab #2 #3 #4 #5 .\\\else\gdef\epsflly{#2}\gdef\epsfurx{#3}\gdef\epsfury{#4}\fi} % gdef because epxfbox wraps us in a group

% pdfTeX completely ignores established conventions for \special
\newif\ifpdf \pdffalse \ifx\pdfoutput\undefined\else\ifx\pdfoutput\relax\else\ifnum\pdfoutput<1 \else\pdftrue\fi\fi\fi
\ifpdf
\pdfcompresslevel=0 % HTTP offers better compression methods anyway
\pdfobjcompresslevel=0
\def\pdflink#1#2{\leavevmode \lhighlight{\pdfstartlink user { /Subtype /Link /Border [0 0 0] /A << /S #1 >> }#2\pdfendlink}}
\def\llink#1{\pdflink{/GoTo /D (#1)}}
\def\link#1{\pdflink{/URI /URI (#1)}}
\def\anchor#1#2{\pdfdest name {#1} xyz #2}

\def\setpapersize#1#2{\pdfpagewidth#1 \pdfpageheight#2 }
\def\dumpbox#1#2#3{\setpapersize{#1}{#2}\shipout\box#3}
\def\metadata#1#2{\pdfinfo{/Title (#1) /Author (#2)}}
\input supp-pdf % the only external dependence of DPMAC: a converter from METAPOST's output to PDF
\def\mps#1{\convertMPtoPDF{#1}{1}{1}}
%\chardef\makeMPintoPDFobject=1
\fi

% Embedded METAPOST
\newtoks\buffertoks
\def\addcode{\immediate\write\mpout}
\def\addunexpandedcode#1{{\toks0={#1}\addcode{\the\toks0}}}
\def\addcodebuffer#1{\edef\tmp{#1}\buffertoks\expandafter\expandafter\expandafter{\expandafter\the\expandafter\buffertoks\tmp}}
\def\addunexpandedcodebuffer#1{\buffertoks\expandafter{\the\buffertoks#1}}

\def\grabcode{\catcode`\#=12 \endlinechar=10 
	\afterassignment\dumpcode\outtoks} % make sure the newline characters get written
\def\dumpcode{\addcode{\the\outtoks}\endinlinemp\gobble} % \gobble the newline with character code 10
\def\begininlinemp{\inimp\begingroup\catcode`\^=7 \iftypesetting\mps{\filestem.\the\figno}\let\addcode\gobble\fi \addcode{beginfig(\the\figno);}}
\def\endinlinemp{\addcode{endfig;}\addcode{}\endgroup\global\advance\figno1\relax}
\ifx\endprolog\undefined\let\endprolog\relax\fi

% Repair LaTeX damage
\def\plainfmtname{plain}\ifx\fmtname\plainfmtname\else
\edef\plainoutput{\the\output}
\global\chardef\itfam=4
\def\_{\leavevmode \kern.06em \vbox{\hrule width.3em}} % for LaTeX

\outputpenalty=0
\tracingstats=0
\newlinechar=-1
\maxdeadcycles=25
\showboxbreadth=5
\showboxdepth=3
\errorcontextlines=5
\overfullrule=5pt
\maxdepth=4pt
\parindent=20pt
\abovedisplayskip=12pt plus 3pt minus 9pt
\belowdisplayskip=12pt plus 3pt minus 9pt
\belowdisplayshortskip=7pt plus 3pt minus 4pt

\font\teni=cmmi10 % math italic

\font\tensy=cmsy10

\catcode"18=12
\catcode`@=11
{
\global\let\end\@@end
\global\let\input\@@input
}
%\let\everymath\frozen@everymath
%\let\everydisplay\frozen@everydisplay
% copied from plain.tex:
\def\eqalign#1{\null\,\vcenter{\openup\jot\m@th
  \ialign{\strut\hfil$\displaystyle{##}$&$\displaystyle{{}##}$\hfil
      \crcr#1\crcr}}\,}
\catcode`@=12
%\def\line{\hbox to\hsize}
%\everymath{}
%\everydisplay{}
\fi

% AMS fonts: family 8, 9, 10
\font\tenmsa=msam10 \font\sevenmsa=msam7 \font\fivemsa=msam5 \newfam\msafam \textfont\msafam=\tenmsa \scriptfont\msafam=\sevenmsa \scriptscriptfont\msafam=\fivemsa
   %\newfam\msbfam \textfont\msbfam=\tenmsb \scriptfont\msbfam=\sevenmsb \scriptscriptfont\msbfam=\fivemsb
%\def\Bbb{\fam\msbfam}
\font\teneufm=eufm10 \font\seveneufm=eufm7 \font\fiveeufm=eufm5 \newfam\eufmfam \textfont\eufmfam=\teneufm \scriptfont\eufmfam=\seveneufm \scriptscriptfont\eufmfam=\fiveeufm

\font\teneufb=eufb10 \font\seveneufb=eufb7 \font\fiveeufb=eufb5 \newfam\eufbfam \textfont\eufbfam=\teneufb \scriptfont\eufbfam=\seveneufb \scriptscriptfont\eufbfam=\fiveeufb

\font\teneurm=eurm10 \font\seveneurm=eurm7 \font\fiveeurm=eurm5 \newfam\eurmfam \textfont\eurmfam=\teneurm \scriptfont\eurmfam=\seveneurm \scriptscriptfont\eurmfam=\fiveeurm
\def\eurm{\fam\eurmfam}
\font\teneurb=eurb10 \font\seveneurb=eurb7 \font\fiveeurb=eurb5 \newfam\eurbfam \textfont\eurbfam=\teneurb \scriptfont\eurbfam=\seveneurb \scriptscriptfont\eurbfam=\fiveeurb

\font\teneusm=eusm10 \font\seveneusm=eusm7 \font\fiveeusm=eusm5 \newfam\eusmfam \textfont\eusmfam=\teneusm \scriptfont\eusmfam=\seveneusm \scriptscriptfont\eusmfam=\fiveeusm

\font\teneusb=eusb10 \font\seveneusb=eusb7 \font\fiveeusb=eusb5 \newfam\eusbfam \textfont\eusbfam=\teneusb \scriptfont\eusbfam=\seveneusb \scriptscriptfont\eusbfam=\fiveeusb

%\font\teneuex=euex10 \font\seveneuex=euex7 \newfam\euexfam \textfont\euexfam=\teneuex \scriptfont\euexfam=\seveneuex

\font\tenss=cmss10 \font\sevenss=cmss7 \font\fivess=cmss5 \inewfam\ssfam \textfont\ssfam\tenss \scriptfont\ssfam\sevenss \scriptscriptfont\ssfam\fivess
\def\sf{\fam\ssfam}

\font\sevenit=cmti7 \scriptfont\itfam=\sevenit

% Font style
\let\articletitle\seventeenss
\let\chaptertitle\twelvebf
\let\sectiontitle\tenbf
\let\subsectiontitle\tenbfit
\let\subsubsectiontitle\tenit
\let\contchaptertitle\tenbf % chapter titles in the table of contents
\let\contsectiontitle\tenrm % ditto for sections
\let\contsubsectiontitle\sevenrm % ditto for subsections
\let\contsubsubsectiontitle\fiverm % ditto for subsubsections
\let\parnumfont\tenrm % paragraph numbers
\let\parbackreffont\fiverm % back references at the end of paragraphs
\let\proclaimfont\tenbf
\let\prooffont\tenit
\let\mainfont\tenrm

% Formatting of draft comments

\def\lhighlight{} % in principle, links can be highlighted, but this is usually done by the renderer
 % for screen
 % for printing drafts
\def\suppresscomments#1#2#3{} % for sharing drafts
\def\prohibitcomments#1#2#3{\errmessage{Draft comments are not allowed in the final version}} % for final versions

\ifx\format\undefined\else\format\fi % default paper size is letter
\ifx\comment\undefined\def\comment{\prohibitcomments}\fi

\newdimen\hmargin
\hmargin=1in
\newdimen\vmargin
\vmargin=1in
\newdimen\plaintextwidth
\plaintextwidth=6.5in
\newdimen\plaintextheight
\plaintextheight=8.9in

% Output routine for scrolls with variable height
\newdimen\totalht
\newdimen\totalwd
\def\shipbox#1{%
	\totalht\ht#1
	\advance\totalht2\vmargin
	\totalwd\wd#1
	\advance\totalwd2\hmargin
	\hoffset-1in
	\advance\hoffset\hmargin
	\voffset-1in
	\advance\voffset\vmargin
	\dumpbox\totalwd\totalht{#1}}

% Contents
\newtoks\cont
\def\contents{\begingroup\suppressbackreftrue\the\cont\endgroup}
\let\printcont\gobble
\def\contlinechapt#1#2{\printcont{#2}\smallskip\everypar{}\noindent{\contchaptertitle\llink{chapter.#1}{#2}}\par} % Link: contents to chapter
\def\contlinesect#1#2{\printcont{#1.  #2}\everypar{}\indent{\contsectiontitle\llap{#1.\enskip}\llink{section.#1}{#2}}\par} % Link: contents to section
\def\contlinesubsect#1#2{\printcont{#1.  #2}\everypar{}\indent{\contsubsectiontitle{#1.\enskip}\llink{paragraph.#1}{#2}}\par} % Link: contents to subsection
\def\contlinesubsubsect#1#2{\printcont{#1.  #2}\everypar{}\indent\indent\indent{\contsubsubsectiontitle\llap{#1.\enskip}\llink{paragraph.#1}{#2}}\par} % Link: contents to subsubsection
\def\addcont#1#2#3{\append\cont{#1}\appendexpand\cont{{#2}}\append\cont{{#3}}}

% Chapters
\newif\ifpresec \presecfalse % true if we are in a chapter, but before its first section
\def\chapter#1\par{%
	\def\chapname{#1}%
	\parn0
	\subparn0
	\presectrue
	\numbfalse
	\addcont\contlinechapt\chapname{#1}%
	\curverb{Chapter~}%
	\assignlabel{chapter}{\chapname}%
	\tchapter#1\par}
\def\tchapter#1\par{% typeset chapter title
	\everypar{\let\beforesect\beforesection}% do not number paragraphs before the first section
	\chapbreak\bigbreak
	\centerline{\plabel\chaptertitle\anchor{chapter.#1}{#1}}% Anchor: chapter
	\nobreak\medskip
	\let\beforesect\relax % do not break pages right after a chapter title
}

\let\chapbreak\relax % \let\chapbreak\tchapbreak to start chapters on a new page

% Sections
\newcount\secn \secn0
\def\section#1\par{%
	\ifx\sectionid\undefined\advance\secn1 \edef\sectionid{\the\secn}\fi
	\everypar{\numpar}% number paragraphs inside a section
	\parn0
	\subparn0
	\presecfalse
	\numbfalse
	\addcont\contlinesect\sectionid{#1}%
	\curverb{\S}%
	\assignlabel{section}{\sectionid}%
	\tsection#1\par
}
\def\tsection#1\par{
	\beforesect\let\beforesect\beforesection
	\typesetsection{#1}%
	\aftersection	
	\let\sectionid\undefined
}
\def\beforesection{\vskip0pt plus.3\vsize \penalty-250 \vskip0pt plus-.3\vsize \bigskip \vskip\parskip}
\def\aftersection{\nobreak\smallskip}
\def\typesetsection#1{\leftline{\sectiontitle\ifx\sectionid\undefined\indent\else\hbox to \parindent{\hss\plabel\anchor{section.\sectionid}{\sectionid}\enspace\hfill}\fi#1}} % Anchor: section
\let\beforesect\beforesection

% Subsections
\def\subsection#1\par{\bigbreak\numbtrue\curverb{\S}\subsectiontitle\noindent#1\/\mainfont\par\nobreak\medskip\addcont\contlinesubsect{\the\secn.\the\parn}{#1}}
\def\subsubsection#1\par{\bigbreak\numbtrue\curverb{\S}\subsubsectiontitle\noindent#1\/\mainfont\par\nobreak\medskip\addcont\contlinesubsubsect{\the\secn.\the\parn}{#1}}

% Theorems and proofs
\def\sskip#1{\ifdim\lastskip<\medskipamount \removelastskip\penalty#1\medskip\fi}
\def\slug{\hbox{\kern1.5pt\vrule width2.5pt height6pt depth1.5pt\kern1.5pt}}
\newif\ifqed \newif\ifneedqed
\def\qed{\unskip\nobreak\ \slug\ifhmode\spacefactor3000 \fi\global\qedtrue}
\def\proclaim{\medbreak\atendpar{\sskip{55}}\numbtrue\gproclaim\proclaimfont} % always number theorems etc.
\def\proof{\medbreak\atendpar{\sskip{-55}}\needqedtrue\gproclaim\prooffont}
\def\xproclaim#1.{\medbreak\atendpar{\sskip{55}}{\everypar{}\noindent}{\proclaimfont#1.\enspace}\ignorespaces} % like \proclaim, no numbering
\let\abstract\xproclaim % abstracts don't get numbered

% Paragraphs
 % pseudo paragraph, cannot be labeled, does not have backrefs
 % ersatz paragraph with numbering, cannot be labeled, does not have backrefs
\def\ppar{\endgraf{\everypar{}\indent}} % paragraph inside a proof, not seen by \proof, has backrefs
\newtoks\atendpar % tokens to be inserted after \endgraf
\newtoks\atendbr % paragraph backreferences
\newif\ifnumb \numbfalse % number this paragraph?
\def\finishpar{\ifhmode\ifneedqed\ifqed\else\qed\fi\qedfalse\needqedfalse\fi\iflist\endlist\fi\the\atendbr\atendbr{}\endgraf\the\atendpar\atendpar{}\numbfalse\fi}
\def\endlist{\iflist\listfalse\endgraf{\parskip\smallskipamount\everypar{}\noindent}\fi} % terminate a list without starting a new paragraph
\newcount\parn
\newcount\subparn
\newif\ifparbref
\def\nextpar{\ifnumb\advance\parn1 \printlabel{advancing paragraph number to \the\parn}\def\brt{}%
	\ifparbref\edef\cseq{\csname backreference-list.paragraph.\the\secn.\the\parn\endcsname}%
	\expandafter\ifx\cseq\relax\else\edef\brt{\cseq}\printbackref{back references for paragraph.\the\secn.\the\parn: \cseq}\fi\fi
	\assignlabel{paragraph}{\the\secn.\the\parn}\fi}
\def\numpar{\ifnumb\nextpar\expand{\noexpand\typesetparnum{\the\secn.\the\parn}\noexpand\typesetpbr{\brt}}\fi}
\newdimen\pindent \pindent\parindent

% Paragraph and theorem numbering
\ifx\draftnum\undefined % normal numbering
\def\gproclaim#1#2.{\curverb{#2~}% #1: font, #2: prefix
	\ifnumb\nextpar\fi
	{\everypar{}\noindent}%
	\plabel#1#2%
	\ifnumb\ \anchor{paragraph.\the\secn.\the\parn}{}\the\secn.\the\parn\fi.\enspace\mainfont % Anchor: proclaim
	\ifnumb\expand{\noexpand\typesetpbr{\brt}}\fi
	\ignorespaces}
\def\typesetparnum#1{\ifnumb{\plabel\parnumfont\anchor{paragraph.#1}{}#1.\enspace}\fi} % Anchor: explicitly numbered paragraph
\def\typesetpbr#1{\ifnumb\def\brtext{#1}\ifx\brtext\empty\else\setetok\atendbr{{\parbackreffont Used in \noexpand\stripcomma\brtext.}}\fi\fi}
\else % number everything for proofreading purposes
\let\numbfalse\relax \numbtrue
\def\gproclaim#1#2.{\curverb{#2~}\medbreak\noindent#1#2.\enspace\mainfont\ignorespaces}
\inewdimen\brwidth \brwidth.6in
\parindent0pt \parskip1ex plus 1ex minus 1ex
\def\typesetparnum#1{\ifnumb\llap{\plabel\anchor{paragraph.#1}{}\parnumfont#1\enspace}\fi} % Anchor: implicitly numbered paragraph
\def\typesetpbr#1{\ifnumb\def\brtext{#1}\ifx\brtext\empty\else
	\llap{\smash{\vtop{\everypar{}\raggedright\rightskip0pt plus 0pt \leftskip0pt plus 1fill \hsize\brwidth
	\parnumfont \strut \break % the first line contains paragraph number
	\parbackreffont\stripcomma#1}}\enspace}\fi\fi}
\fi

% Lists
\def\hang{\hangindent\pindent}
\newif\iflist
\def\textindent#1{{\everypar{}\parindent\pindent\indent}\llap{#1\enspace}\listtrue\ignorespaces}

\def\li{\item{$\bullet$}}

% \par replaced by \endgraf
\def\item{\endgraf\hang\textindent}
\def\filbreak{\endgraf\vfil\penalty-200\vfilneg}

\def\eject{\endgraf\break}
\def\supereject{\endgraf\penalty-20000}
\def\smallbreak{\endgraf\ifdim\lastskip<\smallskipamount
	\removelastskip\penalty-50\smallskip\fi}
\def\medbreak{\endgraf\ifdim\lastskip<\medskipamount
	\removelastskip\penalty-100\medskip\fi}
\def\bigbreak{\endgraf\ifdim\lastskip<\bigskipamount
	\removelastskip\penalty-200\bigskip\fi}

% Email obfuscation (source and output)
%\newdimen\plht \setbox0\hbox{\char`\_} \plht\ht0 \setbox1\hbox{.} \advance\plht-\ht1 % compute the amount by which a dot accent should be lowered to become a period
%\catcode`\.\active \def.{\lower\plht\hbox{\char`\_}} \catcode`\.=12 % emails are typeset with a lowered dot accent instead of a period
%\def\email{\bgroup\catcode`.\active\xemail}
%\def\xemail#1#2{\rlap{\hphantom{#2}@#1}#2\hphantom{@#1}\egroup}

% Back references for labels and bibliography
\let\printbackref\gobble
\def\predefbackref#1{%
	\printbackref{defining back reference list backreference-list.#1}%
	\expandafter\gdef\expandafter\cseq\expandafter{\csname backreference-list.#1\endcsname}
	\expandafter\ifx\cseq\relax\expandafter\gdef\cseq{}\else\printbackref{duplicate omitted}\fi
}
\newcount\backref \backref0
\newif\ifsuppressbackref \suppressbackreffalse
\def\firstletter#1#2\endletter{#1}
\newtoks\backreflist
\newif\ifaddbr \addbrfalse
\def\recordbackref#1{%
	\edef\params{{\ifpresec\expandafter\firstletter\chapname\endletter\else\the\secn\fi.\the\parn\ifnumb\else*\fi}{\the\backref}{\the\inputlineno}}%
	\printbackref{recording back reference \string#1 for future processing with params \params}%
	\edef\tmp{\the\backreflist\noexpand\processbackref\noexpand#1\params}%
	\global\backreflist\expandafter{\tmp}%
	\printbackref{new content of backreflist: \the\backreflist}} % \global is needed because this can be invoked inside {\it ...}, say
\def\processbackref#1#2#3#4{%
	\edef\key{\expandafter\gobble\string#1}%
	\edef\cseq{\csname id.\key\endcsname}%
	\expandafter\ifx\cseq\relax \warningline{#4}{Undefined reference \string#1}\else
		\edef\lseq{\csname backreference-list.\cseq\endcsname}%
		\printbackref{processing back reference number #3 \string#1, originating from #2, at line #4; adding to \lseq}%
		% Check for duplicates
		\edef\lastnumber{\csname lastnumber.\cseq\endcsname}%
		\edef\newnumber{#2}%
		\addbrtrue
		\printbackref{lastnumber: \lastnumber; newnumber: \newnumber;}%
		\expandafter\ifx\csname lastnumber.\cseq\endcsname\relax % we are the first back reference
		\else\ifx\lastnumber\newnumber % same as the last one
			\printbackref{Suppressing duplicate back reference #2.}%
			\addbrfalse
		\fi\fi
		\expandafter\edef\csname lastnumber.\cseq\endcsname{#2}% record the new number
		\ifaddbr
			\expandafter\expandafter\expandafter\gdef\expandafter\expandafter\csname backreference-list.\cseq\endcsname\expandafter{\lseq, \llink{backreference.#3}{#2}}% Link: from a back reference list to the point of origin
			% Step 1: \expandafter\gdef\expandafter\"backreference-list.\cseq"\expandafter{\lseq, \llink{backreference.#3}{#2}}\fi
			% Step 2: \gdef\"backreference-list.\cseq"{\"expanded lseq", \llink{backreference.#3}{#2}}\fi
		\fi
	\fi
}
\newtoks\labelinitlist
\def\xxstripcomma, {}
\def\xstripcomma{\if\ntok,\let\xcont\xxstripcomma\else\let\xcont\relax\fi\xcont}
\def\stripcomma{\futurelet\ntok\xstripcomma}

% Prevention of duplicate labels
\inewif\ifrecorddups \recorddupstrue
\def\checkduplicates#1#2{\edef\key{\expandafter\gobble\string#1}%
	\iftypesetting\else
	\expandafter\ifx\csname line:\key\endcsname\relax\printlabel{keydefline: relax}\ifrecorddups\expandafter\xdef\csname line:\key\endcsname{\the\inputlineno}\fi\else\edef\keydefline{\csname line:\key\endcsname}\errmessage{#2}\fi\fi}

% Labels
\let\printlabel\gobble
\newtoks\curverb % curverb stores the current block name, such as "Chapter", "Theorem", etc.
\ifx\draftlabel\undefined
\def\plabel{}
\else
\def\labeltext{} % label text for proofreading purposes
\def\plabel{\ifx\labeltext\empty\else\smash{\llap{\parbackreffont\labeltext\quad}}\gdef\labeltext{}\fi} % gdef because \plabel is used inside boxes
\fi
\def\assignlabel#1#2{% #1: name (e.g., reference, paragraph, section, chapter), #2: text (e.g., 2.1)
	\ifx\lastlabel\undefined\else
	\edef\key{\expandafter\expandafter\expandafter\gobble\expandafter\string\lastlabel}% label name without backslash
	\printlabel{label \key: id.\key\space = #1.#2, text.\key\space = #2}%
	\expandafter\xdef\csname id.\key\endcsname{#1.#2}% id for DVI hrefs and backreference lists
	\expandafter\xdef\csname text.\key\endcsname{#2}% text that is actually typeset
	\edef\tmp{\the\curverb}%
	\ifx\tmp\empty\else % if curverb is nonempty, define a "verbal" label with a prefix "v"
	\printlabel{label v\key: id.v\key\space = #1.#2, text.v\key\space = \the\curverb#2}%
	\expandafter\xdef\csname id.v\key\endcsname{#1.#2}% id for DVI hrefs and backreference lists
	\expandafter\xdef\csname text.v\key\endcsname{\the\curverb#2}% text that is actually typeset
	\fi
	\predefbackref{#1.#2}%
	\fi\let\lastlabel\undefined}
\newtoks\vlist % Verification list
\def\label#1{%
	\iftypesetting\else
	\checkduplicates#1{Label \string#1 was already defined at line \keydefline}%
	\addverunused#1\verifylabel
	\def\lastlabel{#1}%
	\fi
	\numbtrue % always number labeled paragraphs
}

% Bibliography
\newif\ifyearkey % use years as bibliographic keys?
\def\y{} % use in bibliography as \y{1967}
\newdimen\bibindent % the maximum width of a bibliographic key
\newtoks\bibt % token list for all bibliographic items
\def\tbib#1{% #1 = \Paper
	\checkduplicates#1{Bibliographic reference \string#1 already defined at line \keydefline}%
	\addverunused#1\verifybib
	\edef\key{\expandafter\gobble\string#1}% reference name without backslash
	\printlabel{reference \key: id.\key\space = reference.\key, text.\key\space = \key}%
	\expandafter\edef\csname id.\key\endcsname{reference.\key}% id for DVI hrefs and backreference lists
	\expandafter\edef\csname text.\key\endcsname{\key}% text that is actually typeset
	\predefbackref{reference.\key}%
	\ifyearkey\else	\setbox0=\hbox{[\key]}\ifdim\bibindent<\wd0 \bibindent=\wd0\fi \fi % \bibindent holds the maximum length of all reference [keys]
	\appendexpand\bibt{\noexpand\typesetbib\noexpand#1\src}%
	\ifyearkey\let\next\xxftbib\else\let\next\ftbibalpha\fi\next}
\def\ftbibalpha#1\par{\append\bibt{#1\par}}
\def\xxftbib{\futurelet\next\xftbib}
\def\xftbib{\if\next[\let\next\ftbibyear\else\let\next\ftbibnoyear\fi\next}
\def\ftbibyear[#1]{\edef\yearkey{#1}\expandafter\ftbibyearbis\ignorespaces}
\def\ftbibyearbis#1\par{\append\bibt{#1\par}\ftbibend}
\def\ftbibnoyear#1\par{\append\bibt{#1\par}\edef\yearkey{\extractyear#1\par}\ftbibend}
\def\ftbibend{\expandafter\ifx\csname year.\yearkey\endcsname\relax
		\edef\yearindex{0}%
	\else
		\edef\yearindex{\csname year.\yearkey\endcsname}%
	\fi
	\count0=\yearindex
	\advance\count0 by 1
	\expandafter\edef\csname year.\yearkey\endcsname{\the\count0 }%
	\edef\alphakey{\ifcase \count0 ?\or a\or b\or c\or d\or e\or f\or g\or h\or i\or j\or k\or l\or m\or n\or o\or p\or q\or r\or s\or t\or u\or v\or w\or x\or y\or z\else .\the\count0 \fi}%
	\printlabel{key \key, year key \yearkey, alpha key \alphakey}%
	\expandafter\edef\csname text.\key\endcsname{\noexpand\typesetyearalpha{\key}{\yearkey}{\alphakey}}%
	\printlabel{reference \key\space adjustment: id.\key\space = reference.\key, text.\key\space = \yearkey.\alphakey}%
	\setbox0=\hbox{[\yearkey.\alphakey]}\ifdim\bibindent<\wd0 \bibindent=\wd0 \fi % \bibindent holds the maximum length of all reference [keys]
}
\def\typesetyearalpha#1#2#3{%
	\edef\yearindex{\csname year.#2\endcsname}%
	\ifnum\yearindex=1 #2\else#2.#3\fi}
\def\extractyear#1\y#2#3\par{#2}
\newif\iftype
\def\typesetbib#1#2\par{\edef\key{\expandafter\gobble\string#1}%
	\edef\bibbr{\csname backreference-list.reference.\key\endcsname}%
	\typetrue\ifsuppressunusedbib\ifx\bibbr\empty\typefalse%\warning{Suppressing unused bibliography item \string#1}
	\fi\fi
	\iftype
	\noindent\hbox to \bibindent{[\anchor{reference.\key}{\csname text.\key\endcsname}]\hfil}#2% Anchor: reference
	\ifx\bibbr\empty\else\expandafter\stripcomma\bibbr.\fi
	\hangindent\bibindent\filbreak
	\fi}

% Diagnostic messages for unused references and labels
\let\printverify\gobble
\def\addverunused#1#2{\appendexpand\vlist{\noexpand#2\noexpand#1{\the\inputlineno}}}
\def\verifyref#1#2#3{\printverify{verifying for #3 \string#1 (line #2)}%
	\edef\key{\expandafter\gobble\string#1}%
	\edef\cseq{\csname id.\key\endcsname}%
	\edef\tmp{\csname backreference-list.\cseq\endcsname}%
	\ifx\tmp\empty\warningline{#2}{#3 \string#1}\fi}
\def\verifylabel#1#2{\verifyref#1{#2}{Unused label}}
\def\verifybib#1#2{\verifyref#1{#2}{Unused reference}}

% URLs
\let\printurl\gobble
\newtoks\urltext
\newtoks\urlt
\newif\ifpunct
\def\urldash{-}
\def\urltilde{{\tensy^^X}} % like \sim
\def\ndash{\def\urldash{--}}
\def\http://{\hfil\penalty900\hfilneg\urltext={http://}\urlt={http:/\negthinspace/}\punctfalse\urlgrab}
\def\https://{\hfil\penalty900\hfilneg\urltext={https://}\urlt={https:/\negthinspace/}\punctfalse\urlgrab}
\def\urlgrab{\catcode`\#=11 \catcode`\&=11 \futurelet\ntok\urldispatch}
\def\urldispatch{%
	\ifx\ntok~\let\proceed\urlcont\else
	\ifcat\noexpand\ntok\space\let\proceed\urlfinish\else
	\ifcat\noexpand\ntok\relax\let\proceed\urlfinish\else
	\let\proceed\urlcont
	\fi\fi\fi\proceed}
\def\urlcont#1{\ifpunct\appendexpand\urltext\punctc\appendexpand\urlt\punctc\punctfalse\fi
	\ifx\ntok~\appendexpand\urltext{\noexpand~}\appendexpand\urlt\urltilde
	\else\if\ntok\ampersand\appendexpand\urltext{&}\appendexpand\urlt{\&}%
	\else\if\ntok\ohash\appendexpand\urltext\ohash\appendexpand\urlt\#%
	\else\if\ntok_\appendexpand\urltext_\appendexpand\urlt\_%
	\else\if\ntok-\appendexpand\urltext-\appendexpand\urlt\urldash
	\else\if\ntok.\puncttrue\def\punctc{.}%
	\else\if\ntok,\puncttrue\def\punctc{,}%
	\else\if\ntok;\puncttrue\def\punctc{;}%
	\else\appendexpand\urltext{#1}\appendexpand\urlt{#1}%
	\fi\fi\fi\fi\fi\fi\fi\fi\urlgrab}
\def\urlfinish{\catcode`\#=6 \catcode`\&=4 \hbox{\printurl{\the\urltext}\link{\the\urltext}{\the\urlt}}\ifpunct\punctc\punctfalse\fi\def\urldash{-}}
\def\idgrab{\futurelet\ntok\iddispatch}
\def\iddispatch{\ifcat\noexpand\ntok\space\let\proceed\urlfinish\else\if\ntok,\let\proceed\urlfinish\else\let\proceed\idcont\fi\fi\proceed}
\def\idcont#1{\ifpunct\appendexpand\urltext.\appendexpand\urlt.\punctfalse\fi
	\if\ntok.\puncttrue\def\punctc{.}\else\appendexpand\urltext{#1}\appendexpand\urlt{#1}\fi\idgrab}

% Mathematical string grabbing
\let\printgrab\gobble
\newtoks\grabname
\newtoks\grabtoks % beginning
\newtoks\grabcseq % control sequence name
\newtoks\subsuptoks % ^ or _
\newtoks\dtoks % first letter of a subscript
\newcount\grabsize
\newif\ifgrabsubscript % include subscript?
\newif\ifgrabsupscript % include superscript?
\def\grabsequence{\bgroup % \bgroup allows us to say things like $C^@op$, with "op" being a superscript
	\grabsubscripttrue\grabsupscripttrue % do grab sub/superscripts
	\grabstring}
\def\grabalpha{\bgroup % same
	\grabsubscriptfalse\grabsupscriptfalse % do not grab sub/superscripts
	\grabstring}
\def\grabingroup{\ifinfont\errmessage{Already inside a math token}\fi\append\grabtoks{\bgroup\grablink}\infonttrue}
\def\graboutgroup{\ifinfont\append\grabtoks{\endgrablink\egroup}\infontfalse\fi}
\def\grabstring#1#2#3{% #1 = descriptive name like cat, fun, trans; #2 = font like \bf, \rm, \it; #3 = postcommand like \nolimits
	\let\specialhat^
	\catcode`\^=7
	\aftergroup#3 % #3 could be \nolimits or \limits
	\ifx\specialaddon\undefined\else\expandafter\aftergroup\specialaddon\let\specialaddon\undefined\fi
	\inewif\ifdefine \inewif\ifinfont
	\printgrab{}\printgrab{grab a string of type #1, typeset using font \string#2, with postcommand \string#3}%
	\grabname{#1}\def\grabfont{#2}\grabsize0 \grabtoks={}\grabingroup \append\grabtoks{#2}\grabcseq={}%
	\futurelet\ntok\grabdeflookahead}
\def\grabdeflookahead{\if=\noexpand\ntok % @=Set creates an anchor, whereas @Set refers to it
	\definetrue\printgrab{defining}\expandafter\grabgobblefuturelet
	\else\printgrab{referencing}\definefalse\expandafter\grablookahead\fi}
\def\grabgobblefuturelet#1{\futurelet\ntok\grabtestforsilent} % gobble = and look for another =
\newif\ifsilentgrab
\def\grabtestforsilent{\if=\noexpand\ntok \silentgrabtrue \let\ncom\grabsilenteq \else \silentgrabfalse \let\ncom\grablookahead \fi \ncom}
\def\grabsilenteq={\grabfuturelet}
\def\grabfuturelet{\futurelet\ntok\grablookahead}
\def\grablookahead{\printgrab{futurelet token meaning: \meaning\ntok}%
	\let\ncom\grabfinish
	\if\bgroup\noexpand\ntok \printgrab{left brace, terminating}%
	\else \if\egroup\noexpand\ntok \printgrab{right brace, terminating}%
	\else \if\space\noexpand\ntok \printgrab{blank space, terminating}%
	\else \let\ncom\grabexamine \fi\fi\fi \ncom}
\def\grabexamine#1{\printgrab{grabexamine argument: \string#1, meaning \meaning#1}%
	\def\ncom{\grabfinish#1}%
	\ifcat$\ifcat*\string#1\fi$% is #1 not a command sequence?
		\ifcat _\noexpand#1 \ifgrabsubscript\printgrab{subscript, continuing}%
			\graboutgroup \append\grabtoks{#1}\subsuptoks{#1}\def\ncom{\grabsubsupfuturelet}%
						\else\printgrab{subscript, terminating}\fi
		\else \ifcat ^\noexpand#1 \ifgrabsupscript\printgrab{superscript, continuing}%
			\graboutgroup \append\grabtoks{#1}\subsuptoks{#1}\def\ncom{\grabsubsupfuturelet}%
						\else\printgrab{superscript, terminating}\fi
		\else \ifx \specialhat#1 \ifgrabsupscript\printgrab{specialhat superscript, continuing}%
			\graboutgroup \append\grabtoks{#1}\subsuptoks{#1}\def\ncom{\grabsubsupfuturelet}%
						\else\printgrab{superscript, terminating}\fi
		\else \ifcat\noexpand~\noexpand#1 \printgrab{active character \string#1, examining further}%
			\ifnum1=\uccode`#1 \printgrab{UTF-8 letter, continuing}%
				\advance\grabsize1 \append\grabtoks{#1}\appendexpand\grabcseq{\string#1}\def\ncom{\grabfuturelet}%
			\else
				\ifnum\the\grabsize=0 \printgrab{Nothing grabbed so far, continuing}%
					\advance\grabsize1 \append\grabtoks{#1}\appendexpand\grabcseq{\string#1}\def\ncom{\grabfuturelet}%
				\else\printgrab{Not a UTF-8 letter and not the first character in a string, terminating}%
				\fi
			\fi
		\else \ifcat a\noexpand#1 \printgrab{letter #1, continuing}%
			\advance\grabsize1 \append\grabtoks{#1}\append\grabcseq{#1}\def\ncom{\grabfuturelet}%
		\else\printgrab{nonactive character \string#1}%
			\ifnum\the\grabsize=0 \printgrab{sole argument, adding and terminating}%
				\advance\grabsize1 \append\grabtoks{#1}\append\grabcseq{#1}\def\ncom{\grabfinish}%
			\else\printgrab{terminating}\fi
		\fi\fi\fi\fi\fi
	\else \printgrab{command sequence \string#1, terminating}\fi
	\ncom}
\def\grabsubsupfuturelet{\futurelet\ntok\grabsubsuplookahead}
\newcount\dig
\newif\ifdigit
\def\grabsubsuplookahead{\printgrab{subsup futurelet token meaning: \meaning\ntok}%
	\if\bgroup\noexpand\ntok \printgrab{left brace, continuing}\let\ncom\grabentiresubsup%
	\else \if\egroup\noexpand\ntok \printgrab{right brace, continuing}\let\ncom\grabentiresubsup%
	\else \if\space\noexpand\ntok \errmessage{Blank space after \the\subsuptoks}%
	\else \let\ncom\grabsubsupexamine \fi\fi\fi \ncom}
\def\grabentiresubsup#1{\printgrab{subsup entire group added}\grabingroup\append\grabtoks{#1}\graboutgroup\grabfuturelet}
\def\grabsubsupexamine#1{\printgrab{examining subsup argument \string#1, meaning \meaning#1}%
	% We pass through (1) single letters; (2) command sequences
	\ifcat$\ifcat*\string#1\fi$% is #1 not a command sequence?
		\ifcat\noexpand~\noexpand#1 \printgrab{active character \string#1, continuing}%
			%\warning{math active C: \string#1, grabsize=\the\grabsize, grabtoks=\the\grabtoks, grabcseq=\the\grabcseq}%
			%\grabingroup\appendexpand\grabtoks{\grabfont\noexpand#1}\let\ncom\grabsubsupremainderfuturelet
			\grabingroup\appendexpand\grabtoks{\grabfont\noexpand#1}\let\ncom\grabfuturelet
		\else\ifnum"8000=\the\mathcode`#1 \printgrab{math active character \string#1, continuing}%
			%\warning{math active A: \string#1, grabsize=\the\grabsize, grabtoks=\the\grabtoks, grabcseq=\the\grabcseq}%
			%\grabingroup\appendexpand\grabtoks{\grabfont\noexpand#1}\let\ncom\grabsubsupremainderfuturelet
			\let\specialaddon\egroup
			\def\ncom{#1}%
			\grabtypeset
		\else\ifcat a\noexpand#1 \printgrab{letter #1, checking whether single or not}%
			\dtoks{#1}\def\ncom{\futurelet\ntok\grabsubsupsecondletterlookahead}%
		\else\printgrab{something else, inserting a single-character sub/superscript, continuing}
			\appendexpand\grabtoks{\bgroup\grabfont\noexpand#1\egroup}%
			\advance\grabsize2 \appendexpand\grabcseq{\the\subsuptoks\string#1}% possibly ignore digits here
			\def\ncom{\grabfuturelet}\fi\fi\fi
	\else \printgrab{command sequence \string#1, continuing}%
		\append\grabtoks{#1}\let\ncom\grabfuturelet\fi
	\ncom}
\def\grabsubsupsecondletterlookahead{\def\ncom{\appendexpand\grabtoks{\the\dtoks}\grabfuturelet}%
	\ifcat a\noexpand\ntok \printgrab{not a single letter, grabbing the entire subsupscript}%
		\advance\grabsize1 \appendexpand\grabcseq{\expandafter\string\the\subsuptoks}% append _ or ^ to the label
		\grabingroup
		\appendexpand\grabtoks{\grabfont\the\dtoks}%
		\appendexpand\grabcseq{\the\dtoks}%
		\def\ncom{\grabfuturelet}%
	\else \printgrab{single letter, continuing}\fi\ncom}
\def\grabfinish{\printgrab{grabfinish}\graboutgroup\grabtypeset\egroup}
\def\grabtypeset{\printgrab{grabtypeset grabsize=\the\grabsize, grabtoks=\the\grabtoks, grabcseq=\the\grabcseq}%
	\def\grablink##1\endgrablink{##1}%
	\ifnum\the\grabsize=0 \errmessage{No string to grab}\fi
	\ifnum\the\grabsize>1 % single-letter names are not references
		\ifdefine % defining a mathematical identifier
			\expandafter\checkduplicates\csname\the\grabname.\the\grabcseq\endcsname{Mathematical identifier \key\space already defined at line \keydefline}%
			\iftypesetting % if we are actually typesetting, create an anchor
				\ifsilentgrab
					\expandafter\gdef\csname silent:\the\grabname.\the\grabcseq\endcsname{}% record that this id is silent
				\else
					\anchor{\the\grabname.\the\grabcseq}{}% Anchor: definition of a mathematical identifier
				\fi % only create if not silent
			\else
				\expandafter\xdef\csname id.\the\grabname.\the\grabcseq\endcsname{paragraph.\the\secn.\the\parn}% paragraph id for a back reference list
				\predefbackref{paragraph.\the\secn.\the\parn}%
  			\fi
		\else % referencing a mathematical identifier
			\iftypesetting % we are actually typesetting
				\global\advance\backref1
				\expandafter\ifx\csname line:\the\grabname.\the\grabcseq\endcsname\relax % identifier is undefined
					\warning{Undefined mathematical identifier \the\grabname.\the\grabcseq}%
					\expandafter\gdef\csname\the\grabname.\the\grabcseq\endcsname{\relax}% report undefined references only once; gdef because inside bgroup..egroup
				\else % identifier has ben defined
					\expandafter\ifx\csname silent:\the\grabname.\the\grabcseq\endcsname\empty
						\edef\grablink##1\endgrablink{{##1}}%
					\else % not silent, need a hyperlink
						\edef\grablink##1\endgrablink{\noexpand\llink{\the\grabname.\the\grabcseq}{##1}% Link: from a mathematical identifier to its definition 
							\noexpand\anchor{backreference.\the\backref}{}}% Anchor: back reference for a mathematical identifier reference
					\fi
				\fi
			\else % not yet typesetting, just collecting back references
				\ifsuppressbackref\else
					\global\advance\backref1
					\printbackref{math back reference to \the\grabname.\the\grabcseq: backref.\the\backref\space at line \the\inputlineno}%
					\expandafter\recordbackref\csname\the\grabname.\the\grabcseq\endcsname
				\fi
			\fi
		\fi
	\fi
	\ifsilentgrab\else\the\expandafter\grabtoks\fi}
%		\ifcat 0\noexpand\ntok \printgrab{subsupscript is followed by a catcode 12 character, examining further}%
%			\dig=0 \loop \if\the\dig\ntok \digittrue \fi \ifnum\dig<9 \advance\dig1 \repeat
%			\ifdigit \printgrab{subsupscript is followed by a digit, passed through}%
%				\advance\grabsize2
%				\appendexpand\grabtoks{\subsup#1}%
%				\appendexpand\grabcseq{\expandafter\string\subsup#1}%
%				\def\ncom{\grabfinish}%
%			\else \printgrab{subsupscript is followed not by a digit or a letter, terminating}\fi

% index macros modeled after manmac.tex
% ^={...}: setup an anchor and typeset in italic font
% ^=:{...}: setup an anchor and typeset in roman font
% ^^={...}: setup an anchor, do not typeset
% ^{...}: record a back reference and typeset
% ^^{...}: record a back reference, do not typeset
% ^!{...}: reference a proclaimed statement in the format Theorem~3.5
\inewif\ifsilent \inewif\ifanchor
\newif\ifsuppresscs
\catcode`\^=7   \def\specialhat{\ifmmode\def\next{^}\else\let\next\beginxref\fi\next} \catcode`\^=\active \let^=\specialhat
\def\silentxref#1{\futurelet\next\silentxrefswitch}
\def\silentxrefswitch{\silenttrue\xref}
\def\beginxref{\futurelet\next\beginxrefswitch}
\def\beginxrefswitch{\ifx\next\specialhat\let\next\silentxref \else\silentfalse\let\next\xref\fi \next}
\def\xref{\leavevmode\futurelet\next\xrefswitch}
\def\xrefswitch{\ifx\next!\let\next\verbalxref \else \ifx\next=\let\next\anchorxref \else \anchorfalse \let\next\normalxref \fi \fi \next}
\newtoks\vtoksl
{\count0="C2 \loop\ifnum\count0<"F5 \catcode\count0=11 \advance\count0 by 1 \repeat % make all legitimate initial UTF-8 octets letters
\gdef\plainaccents{\suppresscstrue%
	\def\`##1{##1\empty ̀}%
	\def\'##1{##1\empty ́}%
	\def\^##1{##1\empty ̂}%
	\def\"##1{##1\empty ̈}%
	\def\~##1{##1\empty ̃}%
	\def\=##1{##1\empty ̄}%
	\def\.##1{##1\empty ̇}%
	\def\u##1{##1\empty ̆}%
	\def\v##1{##1\empty ̌}%
	\def\H##1{##1\empty ̋}%
	\def\t##1{##1\empty ͡}%
}}
\def\plainaccents{\let\xcsname=\empty \let\xendcsname=\empty}
\def\verbalxref!{\begingroup\plainaccents\verbalxrefaux}
\def\verbalxrefaux#1{%
	\lowercase{\vtoksl{#1}}%
	\expandtoks\vtoksl
	\iftypesetting
		\expandafter\gdef\expandafter\cseq\expandafter{\csname verbal.\the\vtoksl\endcsname}%
		\printlabel{verbal xref \the\vtoksl}%
		\expandafter\ifx\cseq\relax\errmessage{Undefined reference to \the\vtoksl}\else\cseq\fi
	\else
		\ifsuppressbackref\else
			\global\advance\backref1 %
			\printbackref{label verbal.#1: backref.\the\backref\space at line \the\inputlineno}%
			\blah
			\expandafter\recordbackref\csname verbal.\the\vtoksl\endcsname
		\fi
	\fi
	\endgroup}
\def\initverballabelcommand#1{%
	\printlabel{initializing verbal label #1 (\the\curverb\the\secn.\the\parn)}%
	\expandafter\xdef\csname id.verbal.#1\endcsname{paragraph.\the\secn.\the\parn}%
	\expandafter\xdef\csname text.verbal.#1\endcsname{\the\curverb\the\secn.\the\parn}%
	\expandafter\initlabelcommand\csname verbal.#1\endcsname
}
\def\anchorxref={\anchortrue\futurelet\next\anchorxrefswitch}
\def\anchorxrefswitch{\ifx\next:\let\next\nonitalicanchor\else\italictrue\let\next\normalxref\fi \next}
\def\nonitalicanchor:{\italicfalse\normalxref}
\newtoks\firsttoks \newtoks\secondtoks \inewif\ifplural \inewif\ifitalic
\def\parseplural#1[#2|#3]{\let\next\parseplural\ifx\hfuzz#2\hfuzz\ifx\hfuzz#3\hfuzz\let\next\relax\else\pluraltrue\fi\else\pluraltrue\fi
	\append\firsttoks{#1#2}\append\secondtoks{#1#3}\next}
\newtoks\firsttoksl
\newtoks\secondtoksl
\newtoks\nexttoks
\def\normalxref{\begingroup\plainaccents\normalxrefaux}
\def\normalxrefaux#1{\firsttoks{}\secondtoks{}\pluralfalse\parseplural#1[|]%
	\lowercase\expandafter{\expandafter\firsttoksl\expandafter{\the\firsttoks}}%
	\expandtoks\firsttoksl
	\lowercase\expandafter{\expandafter\secondtoksl\expandafter{\the\secondtoks}}%
	\expandtoks\secondtoksl
	\ifanchor
		\iftypesetting\else
			\initverballabelcommand{\the\firsttoksl}%
			\ifplural\initverballabelcommand{\the\secondtoksl}\fi
			\predefbackref{paragraph.\the\secn.\the\parn}%
			\expandafter\checkduplicates\csname\the\firsttoksl\endcsname{Verbal label \the\firsttoksl\space was already defined at line \keydefline}
			\ifplural\expandafter\checkduplicates\csname\the\secondtoksl\endcsname{Verbal label \the\secondtoksl\space was already defined at line \keydefline}\fi
		\fi
		\anchor{verbal.\the\firsttoksl}{}% Anchor: verbal reference
		\ifplural\anchor{verbal.\the\secondtoksl}{}\fi % Anchor: plural verbal reference
	\else
		\ifsuppressbackref\else
			\global\advance\backref1
			\printbackref{label normal.#1: backref.\the\backref\space at line \the\inputlineno}%
			\iftypesetting
			\else
				\expandafter\recordbackref\csname verbal.\the\secondtoksl\endcsname
			\fi
			\anchor{backreference.\the\backref}{}% Anchor: back reference for a verbal reference
		\fi
	\fi
	% The actual text is typeset after \endgroup
	\nexttoks{}%
	\ifsilent
		\nexttoks{\ignorespaces}%
	\else
		\ifanchor
			\ifitalic\nexttoks\expandafter{\expandafter\bgroup\expandafter\it\the\firsttoks\italcorr}%
			\else\nexttoks\expandafter{\the\firsttoks}%
			\fi
		\else
			\edef\tmp{{verbal.\the\secondtoksl}}%
			\nexttoks\expandafter{\expandafter\llink\tmp}% Link: from a verbal reference to its definition
			\nexttoks\expandafter\expandafter\expandafter{\expandafter\the\expandafter\nexttoks\expandafter{\the\firsttoks}}%
		\fi
	\fi
	\expandafter\endgroup\the\nexttoks}
\def\italcorr{\futurelet\next\italcorrtest}
\def\italcorrtest{\if,\noexpand\next\else\if.\noexpand\next\else\/\fi\fi\egroup}

% Individual repositories

%\def\mathjournals{\ndash\http://www.mathjournals.org}

%\def\matwbn{\http://matwbn.icm.edu.pl/ksiazki}
\def\gen:{\http://libgen.rs/book/index.php?md5=}
\def\jstor:{\https://www.jstor.org/stable/}
\def\eudml:{\https://eudml.org/doc/}

\def\arXiv:{\urltext={https://arxiv.org/abs/}\urlt={arXiv:}\punctfalse\idgrab}

\def\Zbl:{\urltext={https://zbmath.org/?q=an:}\urlt={Zbl:}\punctfalse\idgrab}
\def\doi:{\ndash\urltext={https://doi.org/}\urlt={doi:}\punctfalse\urlgrab}

% Squares
\def\sqr#1#2{{\thinspace\vbox{\hrule height.#2pt \hbox{\vrule width.#2pt height#1pt \kern#1pt \vrule width.#2pt} \hrule height0pt depth.#2pt}\thinspace}}
\def\square{\mathchoice\sqr64\sqr64\sqr{4.2}3\sqr33}

% Long arrows
\def\ltoarr#1{\mathop{\count0=#1 \loop\ifnum\count0>0 \smash-\mkern-7mu \advance\count0 -1 \repeat \mathord\rightarrow}\limits} % parametrized \rightarrowfill
\def\lto#1#2{\mathrel{\ltoarr{#1}^{#2}}} % parametrized \rightarrowfill, with a label
\def\longto#1^#2_#3{\mathrel{\ltoarr{#1}^{#2}_{#3}}} % parametrized \rightarrowfill, with a label above and below
\def\lgetsarr#1{\mathop{\mathord\leftarrow \count0=#1 \loop\ifnum\count0>0 \mkern-7mu\smash-\advance\count0 -1 \repeat}\limits} % parametrized \leftarrowfill
\def\lgets#1#2{\mathrel{\lgetsarr{#1}\limits^{#2}}} % parametrized \leftarrowfill, with a label
\def\longgets#1^#2_#3{\mathrel{\lgetsarr{#1}\limits^{#2}_{#3}}} % parametrized \leftarrowfill, with a label

% Double arrows

\def\toto{\mathrel{\vcenter{\hbox{$\to$}\kern-1.5ex \hbox{$\to$}}}}
\def\prearrfill{\smash-\mkern-7mu}
\def\postarrfill{\mkern-7mu\smash-}
\def\midarrfill#1{\cleaders\hbox{$\mkern-2mu\smash-\mkern-2mu$}\hskip0pt plus #1fil}
\def\rightarrfill{\mkern-7mu\mathord\rightarrow}
\def\leftarrfill{\mathord\leftarrow\mkern-7mu\midarrfill1\postarrfill}
\def\ltoto#1#2#3{\ifinner
	\mathrel{\vcenter{\hbox to #1em{$\prearrfill\midarrfill1{\scriptstyle#2}\midarrfill3 \rightarrfill$}%
		\kern-1.5ex \hbox to #1em{$\prearrfill\midarrfill3{\scriptstyle#3}\midarrfill1 \rightarrfill$}}}%
	\else
	\mathrel{\mathop{\vcenter{\hbox to #1em{\rightarrowfill}%
		\kern-1.5ex \hbox to #1em{\rightarrowfill}}}\limits^{#2}_{#3}}%
	\fi}
\def\ltogets#1#2#3{\ifinner
	\mathrel{\vcenter{\hbox to #1em{$\prearrfill\midarrfill1{\scriptstyle#2}\midarrfill3 \rightarrfill$}%
		\kern-1.5ex \hbox to #1em{$\leftarrfill\midarrfill3{\scriptstyle#3}\midarrfill1 \postarrfill$}}}%
	\else
	\mathrel{\mathop{\vcenter{\hbox to #1em{\rightarrowfill}%
		\kern-1.5ex \hbox to #1em{\leftarrowfill}}}\limits^{#2}_{#3}}%
	\fi}

% New macros for arrows
\def\ltogetscore#1#2{\dimen0=\fontdimen6 #1 2 \divide\dimen0 by 2 \multiply\dimen0 by #2 \vcenter{\hbox to \dimen0{\rightarrowfill}\kern-1.8ex \hbox to \dimen0{\leftarrowfill}}}
\def\ltogets#1#2#3{\mathrel{\mathop{\mathchoice{\ltogetscore\textfont{#1}}{\ltogetscore\textfont{#1}}{\ltogetscore\scriptfont{#1}}{\ltogetscore\scriptscriptfont{#1}}}^{#2}_{#3}}}

% Diagonal arrows
\def\rx#1#2{\rlap{\kern #1pt \raise#1pt \hbox{#2}}}
\def\dottednearrow{\rx{-8}. \rx{-6}. \rx{-4}. \rx{-2}. \rx0. \rx2. \rx4. \kern6pt \raise7.7pt \hbox{$\nearrow$}}

% Rectangular commutative diagrams
% \matrix from Plain TeX: #1 = \hfil$##$\hfil&&\quad\hfil$##$\hfil
\def\gmatrix#1#2{\null\,\vcenter{\normalbaselines
	\ialign{#1\crcr
		\mathstrut\crcr\noalign{\kern-\baselineskip}
		#2\crcr\mathstrut\crcr\noalign{\kern-\baselineskip}}}\,}
\def\cdmatrix{\gmatrix{\hfil$##$\hfil&&\enspace\hfil$##$\hfil\enspace&\hfil$##$\hfil}}
\def\sqmatrix{\gmatrix{\hfil$##$&\enspace\hfil$##$\hfil\enspace&$##$\hfil}}
\def\cdbl{\def\normalbaselines{\baselineskip20pt \lineskip3pt \lineskiplimit3pt }}

\def\cd{\cdbl\cdmatrix}
\def\sqcd{\cdbl\let\vagap\;\sqmatrix}

% Horizontal and vertical diagram arrows
\newcount\arrowsize \arrowsize3
\def\mapright#1{\smash{\lto\arrowsize{#1}}}
\def\mapleft#1{\smash{\lgets\arrowsize{#1}}}
\def\rvagap{\vagap} \def\lvagap{\vagap} \def\rvaskip{\vaskip} \def\lvaskip{\vaskip} \def\vaskip{} \def\vagap{}
\def\mapdown#1{\rvagap\Big\downarrow\rlap{$\vcenter{\hbox{$\scriptstyle#1$}}$}\rvaskip}
\def\mapup#1{\rvagap\Big\uparrow\rlap{$\vcenter{\hbox{$\scriptstyle#1$}}$}\rvaskip}
\def\lmapdown#1{\lvaskip\llap{$\vcenter{\hbox{$\scriptstyle#1$}}$}\Big\downarrow\lvagap}
\def\lmapup#1{\lvaskip\llap{$\vcenter{\hbox{$\scriptstyle#1$}}$}\Big\uparrow\lvagap}

% Wrappable diagrams
\newcount\forno \forno0
\def\arrno#1#2{\global\advance\arr1 \edef\eeqnno{\the\arr}%
	\global\advance\forno1 \edef\eforno{\the\forno}%
	\xdef#2{\noexpand\llink{equation.\eforno}{\eeqnno}}% Link: from an equation number to its definition; xdef to ensure it can be seen outside of \wrapdiagram
	#1{(\anchor{equation.\eforno}{\eeqnno})}} % Anchor: equation number
\newbox\mdiag
\def\wrapdiagram{%
	\setbox\mdiag\vtop\bgroup
	\null % an empty hbox to ensure proper vertical positioning
	\vskip\baselineskip
	\inewcount\arr \arr0
	\baselineskip0pt
	\lineskip4pt
	\lineskiplimit4pt
	\let\par\cr
	\obeylines
	\halign\bgroup\hfil$\displaystyle##$\hfil\cr
	\ewrapdiagram}

\def\ewrapdiagram#1{#1
	\egroup
	\egroup
	\vskip0pt plus \dp\mdiag \penalty-250 \vskip0pt plus-\dp\mdiag % ensure the diagram fits on a single page
	\hangafter-\dp\mdiag
	\divide\hangafter\baselineskip
	\advance\hangafter-2
	\hangindent-\wd\mdiag
	\advance\hangindent-2em
	\hbox to\hsize{\hfil\dp\mdiag0pt \box\mdiag}%
	\ignorespaces}

% METAPOST diagrams                                                                                                     
\newlinechar=10
\inewtoks\preamble
{\endlinechar=10 \catcode`#=12 \global\preamble{
prologues := 3;

verbatimtex
\let\endprolog

\expandafter\gobbleinit\input }\global\appendonceexpand\preamble\jobname\global\append\preamble{
\catcode`\^=7
etex

input cmarrows
setup_cmarrows(arrow_name = "texarrow"; parameter_file = "cmr10.mf"; macro_name = "drawarrow");
setup_cmarrows(arrow_name = "doublearrow"; parameter_file = "cmr10.mf"; macro_name = "drawdarrow");
def drawmarrow expr p = _apth:=p; _finmarr enddef;
rule_thickness#:=.4pt#;    % cmr10.mf: thickness of lines in math symbols
def _finmarr text t_ =
  drawarrow subpath(0, 0.5 * length(_apth)) of _apth t_;
  draw subpath(0.5 * length(_apth), length(_apth)) of _apth withpen pencircle scaled rule_thickness# t_;
enddef;

def object(suffix O)(expr x,y)(expr l) =
  save O;       
  pair O;
  O := (x,y) * u;
  picture O.tx;    
  O.tx := thelabel(l,O);
  draw O.tx;                         
enddef;                                    

def smorphism(suffix A,B) =
  save ss, tt;
  ss := xpart ((A..B) intersectiontimes bbox A.tx);
  tt := xpart ((A..B) intersectiontimes bbox B.tx);
  drawarrow subpath(ss,tt) of (A..B);
enddef;

def morphism(suffix A,B)(expr l)(expr f)(suffix $) =
  smorphism(A, B);
  label.$(l, point f[ss, tt] of (A..B));
enddef;
}}%

\inewcount\vertex                                      
\ifx\mathspecials\undefined\def\mathspecials{}\fi

\def\grabdiagramaux#1;{\vertex0 \dcom#1,*}                          
\def\dcom{\futurelet\next\dcomswitch}                                          
\def\dcomswitch{\ifcat a\noexpand\next \let\next\dcomalpha                                       
        \else\if*\noexpand\next \addcode{\the\buffertoks}\endinlinemp\let\next\gobble                                                 
        \else\if[\noexpand\next \let\next\grabscale                                                                
        \else\errmessage{Unrecognized diagram command \next}\let\next\relax\fi\fi\fi\next}
\def\grabscale[#1]{\dimen0=#1 \addcode{save u; u = \the\dimen0;}\dcom}
\def\dcomalpha#1{\def\objectname{#1}\futurelet\next\dcomalphaswitch}
\def\dcomalphaswitch{\ifcat a\noexpand\next \let\next\dcommorphism\else                                                          
        \expandafter\edef\csname vertex:\objectname\endcsname{\the\vertex}%              
        \if:\noexpand\next\let\next\grabcoords                                              
        \else\if=\noexpand\next\let\next\grabobjectequ                                     
        \else\errmessage{Expected : or = while processing a diagram object, got \meaning\next}\let\next\relax\fi\fi\fi\next}
\def\grabcoords:#1,#2={\toks\vertex{#1,#2}\grabobjectlabel}
\def\grabobjectequ={\edef\tmp{\the\toks\vertex}\ifx\tmp\empty\errmessage{No coordinates specified for vertex \the\vertex: \objectname}\fi\grabobjectlabel}
\def\grabobjectlabel#1,{\addcode{object(\objectname, \the\toks\vertex, }\addunexpandedcode{btex #1 etex);}\advance\vertex1 \dcom}
\def\dcommorphism#1{\def\tobjectname{#1}%
        \def\labelpos{.5}%                                      
        \futurelet\next\dcommorphismswitch}                     
\def\dcommorphismswitch{\if.\noexpand\next \expandafter\grabmorphismdir \else \setupmorphismlabel \expandafter\grabmorphismpos\fi}
\def\setupmorphismlabel{
        \edef\vlabeldira{direction.\objectname.\tobjectname}%
        \edef\vlabeldirb{direction.\tobjectname.\objectname}%
        \expandafter\ifx\csname\vlabeldira\endcsname\relax
        	\expandafter\ifx\csname\vlabeldirb\endcsname\relax
		        \edef\tmpa{\csname vertex:\objectname\endcsname}\expandafter\ifx\tmpa\relax\errmessage{No such vertex: \objectname}\fi
		        \edef\tmpb{\csname vertex:\tobjectname\endcsname}\expandafter\ifx\tmpb\relax\errmessage{No such vertex: \tobjectname}\fi
		        \ifnum \tmpa<\tmpb \edef\vlabeldir{direction.\tmpa.\tmpb}\else \edef\vlabeldir{direction.\tmpb.\tmpa}\fi
		        \expandafter\ifx\csname\vlabeldir\endcsname\relax\errmessage{No label direction specified for morphism \objectname->\tobjectname}\fi
		        \edef\labeldir{\csname\vlabeldir\endcsname}%                        
		\else
        		\edef\labeldir{\csname\vlabeldirb\endcsname}%                        
		\fi
	\else
        	\edef\labeldir{\csname\vlabeldira\endcsname}%                        
	\fi
}                                                                                        
\def\grabmorphismdir.{\def\labeldir{}\futurelet\next\grabmorphismdirec}
{\catcode`\@=13                                                               
\gdef\grabmorphismdirec{\ifx@\next \let\next\grabmorphismposition
        \else\if=\noexpand\next \let\next\grabmorphismequ               
        \else \let\next\grabmorphismdirect\fi\fi\next}                       
\gdef\grabmorphismdirect#1{\edef\labeldir{\labeldir#1}\futurelet\next\grabmorphismdirec}
\gdef\grabmorphismpos{\ifx@\next \expandafter\grabmorphismposition \else \expandafter\grabmorphismequ\fi} 
\gdef\grabmorphismposition@#1={\def\labelpos{#1}\grabmorphismlabel}
}
\def\grabmorphismequ={\grabmorphismlabel}
\def\grabmorphismlabel#1,{%
        \addcodebuffer{morphism(\objectname, \tobjectname, }%
        \addunexpandedcodebuffer{btex \everymath{\scriptstyle}#1 etex, }%
        \addcodebuffer{\labelpos, \labeldir);}%
        \dcom}

% Label extraction and definition
\let\printextract\gobble
\def\hinitlabelcommand#1{\printextract{initializing label command \string#1}%
	\gdef#1{\printlabel{invoked label \string#1}% #1: label command, #2: label id for DVI, #3: typesetted text
		\ifsuppressbackref
			\edef\key{\expandafter\gobble\string#1}%
			\csname text.\key\endcsname
		\else
			\global\advance\backref1 %
			\printbackref{label \expandafter\gobble\string#1: backref.\the\backref\space at line \the\inputlineno}%
			\iftypesetting
			\else
				\blah
				\recordbackref#1
			\fi
			\anchor{backreference.\the\backref}{}% Anchor: back reference for a label (theorem or bibliography)
			\edef\key{\expandafter\gobble\string#1}%
			\llink{\csname id.\key\endcsname}{\csname text.\key\endcsname}% Link: from a label (theorem or bibliography) to its definition
		\fi
		}}
\let\initlabelcommand\hinitlabelcommand
\def\pinitlabelcommand#1{\printextract{initializing label command \string#1}%
	\gdef#1{\printlabel{invoked plain label \string#1}% #1: label command, #2: label id for DVI, #3: typesetted text
		\iftypesetting
			\edef\key{\expandafter\gobble\string#1}%
			\csname text.\key\endcsname
		\else
			\blah
		\fi}}
\newread\labelin
\newif\iflabelcont
\let\terminate=\relax % allow \terminate to be read by \read
\long\def\labelauxaux#1\terminate
{}
\def\preprocesslabel#1{%
	\printextract{Processing label \string#1}%
	\edef\key{\expandafter\gobble\string#1}%
	\initlabelcommand#1
	\expandafter\initlabelcommand\csname v\key\endcsname
	\labelauxaux
}
\def\preprocessbib#1{%
	\printextract{Processing bib \string#1}%
	\initlabelcommand#1
	\labelauxaux
}
\long\def\labelaux#1{\ifx#1\label\let\next\preprocesslabel\else\ifx#1\bib\let\next\preprocessbib\else\let\next\labelauxaux\fi\fi\next}
\def\processoneline{\expandafter\labelaux\labelline\relax\relax\relax\terminate}
\def\preprocess#1{% preprocessing stage to define all labels, back references, bibliographic items, and the table of contents
	% first pass: predefine command sequences for labels and references
	\openin\labelin=#1 \labelconttrue
	\loop
		\read\labelin to\labelline
		\ifeof\labelin\let\next\labelcontfalse\else\let\next\processoneline\fi
		\next
	\iflabelcont\repeat
	% second pass: process labels, bibliography, back references, and the table of contents
	\expandafter\gobbleinit\input#1
	\ifhmode\par\fi\vfill\eject % for LaTeX
	\ifhmode\par\fi\vfill\supereject
}

\def\importlabels#1{ % import labels from an external document
	\recorddupsfalse
	\let\initlabelcommand\pinitlabelcommand
	\preprocess{#1}%
	\let\initlabelcommand\hinitlabelcommand
	\recorddupstrue
	\cont={}% discard the external table of contents
	\vlist={}% discard the verification list
	\let\chapname\undefined \secn0 \parn0 \backref0 % reset the numbering
}

% Primary typesetting routine

\newif\iftypesetting % are we in the final typesetting stage?

% Label and bibliography collection stage
\typesettingfalse % no typesetting at this stage
\def\blah{blah} % placeholder for future references
\output{\setbox0\box255 \setbox0\box\footins \deadcycles0 } % fake typesetting
\let\bib\tbib % collect bibliographic items at this stage
\let\refs\relax % do not typeset bibliography
\everypar{\numpar}
\parn0
\def\par{\finishpar} % no paragraph back references
% Collect embedded METAPOST files and write them to \jobname.gen.mp
\newif\ifmetapost \metapostfalse
\edef\filestem{\jobname.gen}
\newcount\figno
\newwrite\mpout
\newtoks\outtoks
\def\inimp{\ifmetapost\else\global\metaposttrue\immediate\openout\mpout=\filestem.mp \addcode{\the\preamble}\fi}
\figno0
% Temporarily set up large width and height to suppress over/underfull boxes
\newdimen\oldhsize
\oldhsize=\hsize
\newdimen\oldvsize
\oldvsize=\vsize
\hsize\maxdimen
\vsize\maxdimen
\hbadness10000
% Preprocess the file
\preprocess\jobname % collect all labels and bibliography from the main documents and possibly also external labels
% Restore the previously saved dimensions
\hsize\oldhsize
\vsize\oldvsize
\ifx\fmtname\plainfmtname
\hsize\plaintextwidth
\vsize\plaintextheight
\fi
\hbadness1000
% Process back references
\printerr{(\jobname.tex}
\printbackref{list of all back references: \the\backreflist}%
\the\backreflist % process all back references
\ifsuppressunusedbib\def\verifybib#1#2{}\fi
\the\vlist % verify if there are any unused references or bibliography items
\printerr{)}
% Finalize and compile the METAPOST file
\ifmetapost
\addcode{end}
\immediate\closeout\mpout
\immediate\write18{mpost -interaction nonstopmode \filestem.mp}% can also compile the METAPOST file separately
\ifeof18 \warning{Compile the METAPOST file \filestem.mp manually using mpost \filestem.mp}\let\runmp\warning\fi
\fi
\let\addcode\gobble % do not write to the closed file again

% The typesetting stage
\typesettingtrue
\def\importlabels#1{} % stop collecting labels and bibliographic items from other documents
\ifscroll \vsize\maxdimen\inewbox\abox\output{\setbox\abox\vbox{\unvbox255\unskip}\shipbox\abox}%
\else % real typesetting
	\ifx\fmtname\plainfmtname
	\totalht\plaintextheight
	\advance\totalht.1in
	\advance\totalht2\vmargin
	\totalwd\plaintextwidth
	\advance\totalwd2\hmargin
	\setpapersize\totalwd\totalht
	\hoffset-1in
	\advance\hoffset\hmargin
	\voffset-1in
	\advance\voffset\vmargin
	\fi
	\output{\plainoutput}
\fi
\let\blah\undefined % disallow placeholder references
\def\bib#1\par{} % do not collect bibliographic items at this stage
\def\refs{\raggedright\rightskip0em plus \maxdimen \advance\bibindent1em \everypar{}\the\bibt\ignorespaces} % typeset bibliography at this stage
\def\addverunused#1#2{} % do not collect verification data at this stage
\def\addcont#1#2#3{} % do not touch the table of contents at this stage
\def\prepass#1{} % do not collect labels and bibliography from other files at this stage
\parn0
\parbreftrue % typeset paragraph back references at this stage
\everypar{\numpar}
\let\chapname\undefined \secn0 \backref0 % reset the numbering
\figno0
\expandafter\gobbleinit\input\jobname\relax
\ifhmode\par\fi\vfill\eject % for LaTeX
\ifhmode\par\fi\vfill\supereject
\end

\catcode`^=7
\mathcode`@="8000 \catcode`@=13 \def@{\grabalpha{cat}\sf\empty} \catcode`@=12 % categories
\mathcode`"="8000 \catcode`"=13 \def"{\grabalpha{fun}\eurm\empty} \catcode`"=12 % functors
\mathcode`~="8000 \def~{\ifmmode\expandafter\grabtilde\else\penalty10000 \ \fi} % mathematical entities typeset in roman font
\def\grabtilde{\grabalpha{ent}\rm\empty}
\mathcode`?="8000 \catcode`?=13 \def?{\mathop\grabalpha{nat}\rm\nolimits} \catcode`?=12 % natural mathematical operators
\mathcode`\`="8000 \catcode`\`=13 \def`{\grabalpha{simp}\bf\empty} \catcode"60=12 % simplices are typeset in bold
\def\mathspecials{\catcode`@=13 \catcode`"=13 \catcode`!=13 \catcode`?=13 \catcode`\`=13\relax}

\mathchardef\togets="181D

\mathcode`\:="603A % colon is a punctuation mark, not a binary operation
\mathchardef\colon="303A % \colon is a binary operation, not a punctuation mark
\expandafter\def\csname\string≔\endcsname{\colon=} % must redefine ≔ because we redefined : above
\expandafter\def\csname\string∞\endcsname{\ifmmode\infty\else{\tensy1}\fi} % ∞ also works outside of formulas
\expandafter\def\csname\stringσ\endcsname{\ifmmode\sigma\else{\teni\char"1B}\fi} % σ also works outside of formulas
\expandafter\def\csname\stringτ\endcsname{\ifmmode\tau\else{\teni\char"1C}\fi} % τ also works outside of formulas
\catcode`^=\active

\endprolog

\hyphenation{co-do-main}

{\catcode`\'=\active \gdef'{^\bgroup\primes}}
\def\primes{\prime\futurelet\next\primesbis}
\def\primesbis{\ifx'\next\let\nxt\primesbisbis \else\ifx^\next\let\nxt\primet \else\ifx\specialhat\next\let\nxt\primet
  \else\let\nxt\egroup\fi\fi\fi \nxt}
\def\primesbisbis#1{\primes} \def\primet#1#2{#2\egroup}

\begingroup\setbox0=\hbox{$"==colim "==lim ?==id @==Banach @==Top @==Set$}\endgroup % predefine some names

\inewtoks\title
\inewtoks\author
\title{Gelfand-type duality for commutative von Neumann algebras}
\author{Dmitri Pavlov}
\metadata{\the\title}{\the\author}
\font\articletitle=cmss17

%$$
{\tabskip0pt plus 1fil \let\par\cr \obeylines \halign to\hsize{\hfil#\hfil
\articletitle \the\title\vadjust{\bigskip}
\chaptertitle \the\author\vadjust{\medskip}
%$$
Department of Mathematics and Statistics, Texas Tech University
\https://dmitripavlov.org/\vadjust{\medskip}
}}

\abstract Abstract.
We show that the following five categories are equivalent:
(1) the opposite category of commutative von Neumann algebras;
(2) compact strictly localizable enhanced measurable spaces;
(3) measurable locales;
(4) hyperstonean locales;
(5) hyperstonean spaces.
This result can be seen as a measure-theoretic counterpart of the Gelfand duality
between commutative unital C*-algebras and compact Hausdorff topological spaces.
This paper is also available as \arXiv:2005.05284v3.

\tsection Contents

\the\cont

\section Introduction

In 1939, Israel Gelfand [\NR] established a duality between
compact Hausdorff topological spaces and commutative unital ^{C*-algebras}.
This duality can be formulated (Negrepontis [\DT])
as a contravariant equivalence of categories,
where one takes continuous maps of topological spaces and unital ^{C*-homomorphisms} of ^{C*-algebras} respectively
as morphisms.

Based on this, one can conjecture an analogous duality
between an appropriate variant of measurable spaces and ^{commutative von Neumann algebras}.
It is fairly easy to guess some ingredients
for functors going in both directions for such a duality.
Given a measurable space~$X$,
one can construct a commutative complex *-algebra
of bounded complex measurable functions on~$X$.
If $X$ is equipped with a ^{σ-ideal} of ^{negligible sets}
(for example, induced from some ^{measure} by taking its sets of measure~0),
we can define the relation of ^{equality almost everywhere}
and take the quotient of the above algebra by this equivalence relation.
Under additional assumptions, such as ^{σ-finiteness} or, more generally, ^{[localizability|localizable enhanced measurable space]},
this *-algebra is a ^{von Neumann algebra}.
A ^{σ-ideal} of ^{negligible sets} is necessary for the above construction,
so we can expect it to be present in some form in the statement of the duality.
Thus, we consider triples $(X,M,N)$, where $X$ is a set, $M$ is a ^{σ-algebra} of ^{measurable subsets} of~$X$,
and $N⊂M$ is a ^{σ-ideal} of ^{negligible subsets} of~$X$.
We refer to such a triple as an ^{enhanced measurable space} (^!{enhanced measurable space}).

Conversely, given a ^{commutative von Neumann algebra},
one can take the ^{Gelfand spectrum} of its underlying ^{C*-algebra}
or, equivalently, the ^{Stone spectrum} [\BRGT] of its ^{Boolean algebra of projections}.
The resulting topological spaces are known as ^{hyperstonean spaces} (^!{hyperstonean spaces}).
They were introduced and studied by Dixmier [\HS].
Dixmier [\HS, Théorème~2] also proved that ^{hyperstonean spaces} are precisely the ^{Gelfand spectra} of ^{commutative von Neumann algebras}.
Any ^{hyperstonean space} can be equipped with a ^{σ-ideal} of negligible sets comprising precisely the ^{meager subsets}
(as already proposed by Loomis [\RepL] and Sikorski [\RepS])
and a ^{σ-algebra} of measurable subsets comprising precisely the symmetric differences of open and ^{meager subsets}.
This produces an ^{enhanced measurable space} (\vTMdefcorrect).
Only a ^{σ-ideal} of ^{negligible sets} can be defined canonically in this construction, not a specific measure.
Concerning morphisms of ^{commutative von Neumann algebras},
von Neumann [\MA],
C.~Ionescu Tulcea [\SCE],
Vesterstrøm–Wils [\PReal],
Edgar [\MWS],
Graf [\EPM],
as well as von Neumann–Maharam [\VN]
established criteria for lifting ^{homomorphisms of Boolean algebras} or ^{commutative von Neumann algebras}
to point-set measurable maps of measurable spaces.

However, a formulation that promotes these constructions to an actual equivalence of categories does not appear in the literature.
In particular, whereas on the ^{von Neumann algebra} side it is clear that one should take the category
of ^{commutative von Neumann algebras} and ^{normal *-homomorphisms},
on the measure theory side the situation is far from clear.
Some obvious choices for objects and morphisms,
such as ^{localizable enhanced measurable spaces}
with morphisms being equivalence classes of measurable maps modulo ^{equality almost everywhere},
fail to produce a category that is contravariantly equivalent to the ^{category of commutative von Neumann algebras}.

This paper resolves these issues by establishing the following result.

\label\mainthm
\proclaim Theorem.
The following categories are equivalent.
\li The category $@CSLEMS$ of ^{compact strictly localizable enhanced measurable spaces} (^!{compact strictly localizable enhanced measurable space}),
whose objects are triples $(X,M,N)$, where $X$ is a set, $M$ is a ^{σ-algebra} of measurable subsets of~$X$,
$N⊂M$ is a ^{σ-ideal} of negligible subsets of~$X$
such that the additional conditions of compactness (^!{compact enhanced measurable space}) and strict localizability (^!{strictly localizable enhanced measurable space}) are satisfied.
Morphisms $(X,M,N)→(X',M',N')$ are equivalence classes of maps of sets $f:X→X'$ such that $(f^*)_!M'⊂M$ and $(f^*)_!N'⊂N$
(superscript~$*$ denotes the preimage map, subscript~$!$ denotes the direct image map)
modulo the equivalence relation of ^{weak equality almost everywhere} (^^!{weak equality almost everywhere}): $f≈g$ if for all $m∈M'$ the symmetric difference $f^*m⊕g^*m$
belongs to~$N$.
\li The category $@HStonean$ of ^{hyperstonean spaces} and open maps (^^!{hyperstonean space}).
\li The category $@HStoneanLoc$ of ^{hyperstonean locales} and open maps (^^!{hyperstonean locale}).
\li The category $@MLoc$ of ^{measurable locales}, defined as the full subcategory of the ^{category of locales}
consisting of ^{complete Boolean algebras} that admit sufficiently many ^{continuous valuations} (^^!{measurable locale}).
\li The opposite category $@CVNA^@op$ of ^{commutative von Neumann algebras},
whose morphisms are ^{normal *-homomorphisms} of algebras in the opposite direction (^^!{commutative von Neumann algebras}).
\endlist Furthermore, as explained in the paper, the equivalences are implemented by the following four adjoint equivalences, providing a highly structured way to move between these five different settings:
$$@HStonean\ltogets{11}{"Ω}{"Sp}@HStoneanLoc\ltogets{11}{"COpen}{"Ideal}\vtop{
\lineskiplimit=1in \lineskip=.85ex
\setbox0=\hbox{$@MLoc$}
\dimen0=\wd0
\box0
\let\Big=\Bigg
\hbox to \dimen0{\hss\lmapup{"ML}\mapdown{"Spec}\hss}
\hbox to \dimen0{\hss$@CSLEMS$.\hss}
}\ltogets{11}{"L^∞}{"ProjLoc}@CVNA^@op$$

\proof Proof.
Combine
\vHStonean,
^!{hyperstonean duality},
\vCVNAMLoc,
and
\vCSLEMSMLoc.

The last four categories appear to be rather natural,
but the first category $@CSLEMS$ may appear unfamiliar,
both in terms of its objects and morphisms.
We offer some clarifying remarks that may help
to convince the reader that $@CSLEMS$ is a natural and essentially the only possible choice here.

First, the category $@CSLEMS$ is sensitive to the choices
like ^{[localizability|localizable enhanced measurable space]} versus ^{strict localizability},
the presence of ^{[compactness property|compact enhanced measurable space]},
or ^{equality almost everywhere} versus ^{weak equality almost everywhere}.
Only the choices made in \vmainthm\ produce a category contravariantly equivalent to ^{commutative von Neumann algebras}.
In particular, the equivalence relation of ^{weak equality almost everywhere}
differs from the equivalence relation of ^{equality almost everywhere}
if the ^{σ-algebra} is not ^{countably separated} (^!{countably separated}).
See \vcseplemma\ and \vweakeqnecessary.

The requirement that $N$ is a ^{σ-ideal} of~$2^X$ (as opposed to a ^{σ-ideal} of~$M$)
amounts to incorporating the ^{completeness} assumption in our definition of an ^{enhanced measurable space}.
The ^{completeness} condition does not change the resulting category,
but considerably simplifies the presentation, since in this case we can define
morphisms of ^{enhanced measurable spaces} as equivalence classes of point-set maps
whose associated preimage map preserves measurable and negligible subsets.
See ^!{completeness}, which explains how the noncomplete version works and why it is equivalent to the complete version.

The ^{strict localizability} property is a natural generalization of ^{σ-finiteness}
that is necessary since we allow non-^{σ-finite} ^{commutative von Neumann algebras}.
The ^{compactness} property is the abstract measure-theoretic formulation of Radon measures.
Indeed, a Radon measure on a topological space yields a ^{compact enhanced measurable space} by ^!{Radon enhanced measurable space}.
As explained there, the resulting space is also ^{strictly localizable}, which means that the overwhelming majority
of measurable spaces in analysis are in fact ^{compact strictly localizable enhanced measurable spaces}.
The properties of ^{compactness} and ^{strict localizability} are also essential for eliminating pathologies
that inevitably arise with larger classes of ^{enhanced measurable spaces}.
For instance, they are crucial for establishing ^!{measurable image},
which says that a ^{morphism of enhanced measurable spaces} has a well-defined image, which is a measurable set.
Another crucial point at which these conditions are relevant is \vlifting, which lies at the core
of the equivalence between $@CSLEMS$ and $@MLoc$.

The functor $"Spec∘"ProjLoc:@CVNA^@op→@EMS$
(\vdefSpec, ^!{projection locale functor})
produces the ^{Gelfand spectrum} of the underlying commutative ^{C*-algebra}
equipped with a ^{σ-algebra} of measurable subsets and ^{σ-ideal} of negligible subsets,
described above as the Loomis–Sikorski construction.
By \vspeccompactstrictloc, ^{enhanced measurable spaces} in the image of this functor are ^{compact} and ^{strictly localizable}.
Furthermore, by \vcslemsdetiso\ such spaces are closed under isomorphisms in the category $@LDEMS$ of ^{locally determined enhanced measurable spaces}
(^!{locally determined enhanced measurable spaces}).
This is not obvious, since the property of ^{compactness}
is formulated in terms of properties of the ^{σ-algebra} and ^{σ-ideal} that are not manifestly invariant under isomorphisms in $@LDEMS$.
In particular, \vstrictlocnotiso\ and \vcslemsnotiso\ show that ^{strictly localizable enhanced measurable spaces}
and ^{compact strictly localizable enhanced measurable spaces} are not closed under isomorphisms in $@EMS$.
The algebra of equivalence classes of bounded complex-valued functions on a ^{localizable enhanced measurable space} is a ^{commutative von Neumann algebra},
whose ^{Gelfand spectrum} is isomorphic to the original ^{enhanced measurable space} if the latter is ^{compact} and ^{strictly localizable}.
Thus, the category $@CSLEMS$ of ^{compact strictly localizable enhanced measurable spaces}
is the essential image of the functor $@CVNA^@op→@EMS$.
In other words, the category $@CSLEMS$ is forced on us by the nature of this equivalence and we essentially have no other choice.

Larger categories, like the category $@LEMS$ of ^{localizable enhanced measurable spaces},
can also be made equivalent to $@CVNA^@op$ with additional work,
but at the cost of losing the convenient point-set description of morphisms,
since the category of ^{localizable enhanced measurable spaces} that is equivalent to $@CVNA^@op$
is not $@LEMS$ itself, but rather $@LEMSc'$, which is obtained from $@LEMS$ by discarding certain morphisms
and performing a Gabriel–Zisman localization of the resulting nonfull subcategory $@LEMSc$
with respect to a certain class of morphisms, see \vstrictloc.
In the resulting category $@LEMSc'$, morphisms are not equivalence classes of maps of sets,
but rather abstract compositions of such morphisms and their inverses, which makes $@LEMSc'$ less convenient in practice.

Naturally, after \vmainthm\ one is led to wonder whether the various notions of a measure
defined on spaces in these equivalent categories are also equivalent to each other.
This is indeed true, although for simplicity and brevity, we formulate the following result for finite measures,
even though it continues to hold for infinite measures (semifinite or not).

\label\measureeq
\proclaim Theorem.
In the context of the chain of equivalences in \vmainthm,
we have canonical bijections between the following notions of a complex finite measure:
\li Finite ^{essential measures} on an ^{enhanced measurable space} (^!{essential measure});
\li Finite ^{normal measures} on a ^{hyperstonean space} (^!{normal measure});
\li ^{Normal valuations} on a ^{hyperstonean locale} (^!{normal valuation});
\li ^{Continuous valuations} on a ^{measurable locale} (^!{continuous valuation});
\li Elements of the ^{predual} of a ^{commutative von Neumann algebra} (^!{predual}).

\proof Proof.
Combine \vnormalmeasurevaluationhyperstonean,
\vhnormal,
\vpredualvaluation,
and \vmeasurevaluation.

The appearance of ^{essential measures} is motivated as follows.
A finite ^{measure} on an ^{enhanced measurable space} $(X,M,N)$
is a countably additive map $μ:M→`C$ that vanishes on~$N$ (^!{measure}).
Vanishing on~$N$ is analogous to absolute continuity.
A ^{measure}~$μ$ is ^{essential} (^!{essential}) if for any $m∈M$ such that $μ|_m≠0$
we can find a ^{σ-finite} measurable subset $m'⊂m$ such that $μ|_{m'}≠0$.
In other words, if an ^{essential measure} vanishes on all ^{σ-finite} measurable subsets of some $m∈M$, then it vanishes on all measurable subsets of~$m$.
Here a measurable subset $m∈M$ is ^{σ-finite} (^!{σ-finite}) if there is a ^{faithful finite measure}~$ν$ on the ^{induced enhanced measurable space} $(m,M_m,N_m)$
(^!{induced enhanced measurable space}),
and a ^{measure} is ^{faithful} if $ν|_m=0$ implies $m∈N$.

The condition of ^{essentiality} of a ^{measure} is nontrivial only in the non-^{σ-finite} case.
Without this condition, even the most fundamental theorems of measure theory,
such as the Radon–Nikodym theorem, fail in the non-^{σ-finite} case,
even assuming ^{compactness} and ^{strict localizability},
see \vessentialradonnikodym\ and \vessentialcounterexample.
If $(X,M,N)$ is ^{localizable} (^!{localizable}), i.e., the ^{Boolean algebra} $M/N$ is ^{complete},
being ^{essential} is equivalent to being a ^{completely additive measure} (^!{completely additive measure}).
Thus, given \vmeasureeq, we argue that ^{essential measures} are the correct notion of a measure in the non-^{σ-finite} setting.

We conclude with some remarks on the relation between ^{locales} and measure theory,
as featured in two out of five categories under consideration.
An overview of ^{locales} is given in \vreviewlocales, where one can find further references.
The category of ^{measurable locales} (^!{measurable locales}) was proposed by the author on December 15, 2010 in [\MLoc]
as a localic analog of ^{localizable enhanced measurable spaces}.
One passes from a ^{localizable enhanced measurable space} $(X,M,N)$ to the corresponding ^{measurable locale}
by taking the quotient ^{Boolean algebra}~$M/N$.
^{Boolean algebras} of such nature have been studied for a long time:
already in 1942, Halmos and von Neumann [\OMCM] proved that the real line as a measurable space
can be characterized in terms of its ^{Boolean algebra} of equivalence classes of measurable sets.
Volume~3 of Fremlin's {\it Measure theory\/} [\MTiii]
studies {\it localizable measure algebras},
which in our context are precisely ^{measurable locales} equipped with a possibly infinite positive ^{faithful} ^{continuous valuation}.
Jackson [\ShMeas, §3.3] observes that sheaves on a ^{σ-algebra} from a localic topos,
with the ^{equality almost everywhere} defining a Lawvere–Tierney topology,
whose sheaves again form a localic topos.
Assuming ^{[localizability|localizable enhanced measurable space]}, the resulting ^{locale} is isomorphic to the locale $M/N$ discussed above.
However, morphisms of measurable spaces are not discussed in this context.

The ^{category of measurable locales} $@MLoc$ is a full subcategory of the ^{category of locales} $@Loc$.
In particular, morphisms of measurable locales, which correspond to equivalences classes of measurable maps in point-set measure theory,
are in bijective correspondence with morphisms of locales, which correspond to continuous maps in point-set general topology.
This demonstrates the power of pointfree general topology: it is capable of treating both the traditional general topology
and measure theory using the same formalism of locales.

\subsection Structure of the paper

The structure of this paper follows the structure of the chain of functors in \vmainthm.

\vstoneduality\ explores in the detail the middle part of the chain of equivalences in \vmainthm,
concerning the categories $@HStonean$, $@HStoneanLoc$, and $@MLoc$
together with the two adjoint equivalences (\vHStonean, ^^!{hyperstonean duality}) relating them.
It obtains these two adjunctions by restricting the classical adjoint equivalences
between the categories $@Stone$ (^{Stone spaces}), $@StoneLoc$ (^{Stone locales}), and $@BAlg^@op$ (^{Boolean algebras}).
As an intermediate step, we pass through the known adjoint equivalences
between the categories $@Stonean$ (^{Stonean spaces}), $@StoneanLoc$ (^{Stonean locales}), and $@CBAlg^@op$ (^{complete Boolean algebras}).

\vvonneumannalgebras\ constructs the adjoint equivalence (\vCVNAMLoc) between $@MLoc$ and $@CVNA^@op$
using arguments resembling the Lebesgue integration theory.
Parts of this result appeared in the literature before
(e.g., the essential surjectivity of $"ProjLoc$ follows from the results by Segal [\DOA]),
but there does not appear to be a published source that presents a proof of the whole result.

\vpointset\ carefully defines the appropriate category $@CSLEMS$ of point-set measurable spaces
and establishes its elementary properties.

\veqmeasspaceloc\ constitutes the central core of this paper.
It constructs the functors $"Spec:@MLoc→@CSLEMS$ and $"ML:@CSLEMS→@MLoc$
and then constructs from them an adjoint equivalence between $@CSLEMS$ and $@MLoc$.

\subsection Prerequisites

We assume familiarity with the elementary theory
of categories, functors, natural transformations, adjunctions, equivalences,
limits, colimits, and quotient categories.
We also assume familiarity with the elementary theory of topological spaces and topological vector spaces.
We review the necessary facts about ^{locales} in \vreviewlocales\ and ^{von Neumann algebras} in \vreviewvna,
but our exposition is not self-contained and we give pointers to the literature where appropriate,
in particular we cite Johnstone [\SS] as much as we can for ^{Stone duality}
and Sakai [\CWstar] for several facts about ^{von Neumann algebras}.
Our presentation of point-set measure theory is self-contained
except that we refer to Corollary~341Q, Theorem~343B, and Lemma~451Q in Fremlin [\MTiii, \MTiv]
for the lifting theorems and a technical lemma about ^{compact} measure spaces,
and we also cite several counterexamples from Fremlin.

\subsection Conventions

All rings are by definition unital and all homomorphisms of rings by definition preserve units.
This applies, in particular, to C*-algebras and their *-homomorphisms,
which is different from the standard convention.

We use the sans serif font for categories like $@Top$, $@Set$ and the Euler font for functors like $"F$, $"G$, etc.

\proclaim Notation.
If $@C$ is a category, then $@C^@=op$ denotes its opposite category.
Likewise, if $"F:@C→@D$ is a functor, then $"F^"=op:@C^@op→@D^@op$ is its opposite functor.
Finally, if $t:"F→"G$ is a natural transformation of functors $"F,"G:@C→@D$,
then $t^~=op:"G^"op→"F^"op$ is the opposite natural transformation.
%We also use this notation for objects and morphisms:
%if $X∈@C$ is an object in a category~$@C$,
%then $X^@op$ denotes the corresponding object in $@C^@op$.
%Likewise, if $f:X→Y$ is a morphism in~$@C$,
%then $f^"op:Y^@op→X^@op$ denotes the corresponding morphism in~$@C^@op$.

If $f:X→Y$ is a map of sets, then $f^*B⊂X$ denotes the preimage of $B⊂Y$
and $f_!A⊂Y$ denotes the image of $A⊂X$.
The notation $f^{-1}$ is reserved for the inverse of~$f$, i.e., the map $f^{-1}:Y→X$ such that $f^{-1}∘f=?id_X$ and $f∘f^{-1}=?id_Y$.

We use the terms “meet” and “infimum”,
as well as “join” and “supremum” interchangeably for any poset.
Likewise, we make no distinction between
$⋀$ and $\inf$, or between $⋁$ and $\sup$.

Hyperlinks are used extensively in this paper,
not only for standard items, like numbered statements, bibliography, and back references,
including the “Used in” lists at the end of referenced theorems,
but also for terms inside the main text,
such as “^{enhanced measurable space}” or “^{hyperstonean locale}”,
as well as identifiers of mathematical objects like “$@EMS$” and “$"Spec$”,
which are all hyperlinked to their definition.

\subsection Future directions

After establishing \vmainthm,
one almost immediately runs into the question
whether various objects that can be defined
on measurable spaces, such as measures or measurable fields of Hilbert spaces,
also lead to a similar chain of equivalences.
This turns out to be true, but considerations of length prevent us from proving these results in this paper.
The next two subsections briefly sketch the relevant statements.
The third subsection explores the setting of elementary toposes.

\subsubsection Fiberwise measures and disintegrations

As explained in \vmeasureeq, each of the five categories introduced above allows for a notion of measure,
and all these notions turn out to be equivalent in a precise sense.
For applications, it is important to consider {\it relative measures\/} (alias {\it fiberwise measures\/}),
conditional expectations in probability theory being an important special case of relative measures.
Given a map $f:X\to Y$, a relative measure~$μ$ on~$f$ is roughly a choice of a measure~$μ_y$ on $f^*\{y\}$ for each point $y∈Y$
so that the family $\{μ_y\}_{y∈Y}$ is itself in some sense measurable,
for example, the map $y↦μ_y(m∩f^*\{y\})$ must be measurable for any measurable subset $m⊂X$.
If $Y=\{*\}$ is a point, then a relative measure on the unique map $f:X→Y$
is simply a measure on~$X$ in the usual sense.
On the other hand, a relative measure on the identity map $?id:X→X$
should be a measurable function $X→`C$, since a measure on a singleton space is just a number.

If $μ$ is a relative measure on $f:X→Y$ and $ν$ is a relative measure on $g:Y→Z$,
then we can construct a relative measure $ν∘μ$ on $g∘f:X→Z$
by assigning to $z∈Z$ the measure $(ν∘μ)_z$ on $f^*g^*\{z\}$
such that $(ν∘μ)_z(m)=∫_{g^*\{z\}}(y↦μ_y(m∩f^*\{y\})) ~d ν_z$ for any measurable subset $m⊂f^*g^*\{z\}$.
The associativity of composition is established using a generalized Fubini theorem.
Thus, we can take any of the above five categories,
equip its morphisms with an additional data of a relative measure,
and obtain a new category.
This allows us to give a rather general statement about equivalence of various notions of relative measures,
whose proof will appear elsewhere.

\proclaim Pretheorem.
In the context of the chain of equivalences in \vmainthm,
for each of the five categories~$@C$ in the statement we have a corresponding category $\hat@C$,
which has the same objects and whose morphisms are morphisms in the old sense
equipped with an additional fiberwise measure-like structure with individual fibers like in \vmeasureeq.
\ppar
The functors in the chain of equivalences in \vmainthm\ can be promoted to functors between such enhanced categories
(so that on the underlying old data they coincide with the previously defined functors), and they remain adjoint equivalences.
Restricting to real, positive, and/or finite measures also produces a chain of adjoint equivalences.
Finally, the forgetful functors from the above categories to the categories defined in \vmainthm\ form (contravariant) Grothendieck fibrations.

The last claim about Grothendieck fibrations
serves as a really neat formulation of the disintegration theorem of von Neumann, Rohlin, Pachl, and others.

\subsubsection Bundles of Hilbert spaces and topological vector spaces

For each of the five categories introduced above we can define Hilbert bundles,
and more generally, bundles of topological vector spaces,
and all these notions turn out to be equivalent in a precise sense.
We state the following theorem for the case of Hilbert spaces,
since this is the more familiar setting.

\proclaim Pretheorem.
In the context of the chain of equivalences in \vmainthm,
for each of the five categories $@C$ in the statement we have a corresponding category $\hat@C$,
whose objects are pairs $(X,V)$, where $X∈@C$ and $V$ is a bundle of Hilbert spaces over~$X$.
Morphisms $(X,V)→(X',V')$ are pairs $(f,g)$, where $f:X→X'$ is a morphism in~$@C$ and $g:V→f^*V'$
is an appropriately defined fiberwise linear map of bundles.
The functors in the chain of equivalences in \vmainthm\ can be promoted to functors between such enhanced categories
(so that on the $X$-component they coincide with the previously defined functors), and they remain adjoint equivalences.

The above result, for example, immediately implies the spectral theorem for families of commuting normal operators on a Hilbert space.

The more general version for bundles of topological vector spaces that we did not state above also has
many practical applications: it immediately implies the reduction theory for von Neumann algebras, for example.

\subsubsection The setting of arbitrary elementary toposes

Henry [\LMSLGD, Theorem~4.2.5]
shows that in any elementary topos with a natural numbers object
(which generalizes the usual Zermelo–Fraenkel category of sets)
one has a contravariant equivalence
between the category of commutative localic C*-algebras
and the category of compact regular locales.

This improves and strengthens the usual Gelfand duality for C*-algebras
by promoting it to arbitrary elementary toposes
that need not satisfy the axiom of choice or the law of excluded middle.

As an immediate practical application,
one obtains versions
of Gelfand duality for smooth or continuous bundles of C*-algebras,
as well as equivariant analogs.

Naturally, one is led to wonder to what extent the chain of adjoint equivalences in \vmainthm\
continues to hold for elementary toposes with a natural number object.
Of course, just like in the case of C*-algebras, there is no hope
for point-set notions to be equivalent to localic notions
due to the lack of the axiom of choice,
which is necessary in \vCoh\ to show that all ^{coherent locales} (in particular,
all ^{hyperstonean locales}) are ^{spatial},
essentially by ensuring the existence of sufficiently many maximal ideals using Zorn's lemma.
The axiom of choice is also essential for \vlifting, which is the main ingredient
of the equivalence of categories between $@CSLEMS$ and $@MLoc$.
The resulting theory makes it clear that the localic notions are the correct ones,
so the categories $@CSLEMS$ and $@HStonean$ must be discarded.

Of the remaining three categories $@HStoneanLoc$, $@MLoc$, and $@CVNA^@op$,
the first two are already defined in terms of ^{locales}.
Of course, instead of merely demanding the existence of sufficiently
many ^{normal valuations} respectively ^{continuous valuations}
one must instead construct the ^{locale} of ^{normal valuations}
respectively ^{continuous valuations} in a pointfree manner,
and then use it to state the localic analogs of the above properties.

One can also identify some plausible ingredients that are likely
to participate in the localic version of the category $@CVNA$ of ^{commutative von Neumann algebras}.
The first ingredient is a complex *-algebra object~$A$ in the category of ^{locales}
(with complex numbers themselves understood as a ^{locale}).
The underlying ^{locale} of~$A$ encodes the ultraweak topology on a von Neumann algebra.
The second ingredient is a Banach ^{locale}~$A_*$ (defined in Henry [\LMSLGD, Proposition~4.1.3])
corresponding to the ^{predual} equipped with the norm topology.
Finally, we should have an isomorphism $A→(A_*)^*$,
where the right side is interpreted as the localic weak-* dual of~$A_*$ (for the unit ball version see Henry [\LMSLGD, §4.2.3]).

Whether or not the above categories are actually equivalent in any elementary topos
with a natural numbers object is an open problem.

\subsection Acknowledgments

I thank the MathOverflow community,
where a part of \vmainthm\ was first stated in [\GelfandDuality],
with further developments in [\PushPull], [\Duality], and [\MLoc]
proving to be particularly influential in the writing of this paper.

I thank Andre Kornell, whose inquiries about the main theorem in 2012 and 2018 inspired \vstrictloc\ and provided an additional impetus to write this paper.
I thank Simon Henry for several conversations concerning localic Gelfand duality.
I thank Robert Furber for pointing out Lemma 451Q in Fremlin [\MTiv], \vMLnotcomplete, and a missing step in the proof of \vMLfactors.
I thank the anonymous referee of the Journal of Pure and Applied Algebra for suggesting numerous corrections and improvements.

\label\stoneduality
\section Stone-type dualities for locales

This section examines
the ^{hyperstonean duality}, which is used later to construct
^{enhanced measurable spaces} out of ^{measurable locales}.
The duality is expressed as a chain of adjoint equivalences of categories.

We start with a brief overview of locales and the chain of adjoint equivalences
(\vCoh, ^^!{Stone duality for distributive lattices})
$$@Coh\ltogets5{}{}@CohLoc\ltogets5{}{}@DLat^@op$$
between ^{coherent spaces}, ^{coherent locales}, and ^{distributive lattices}
that serves as a foundation for all Stone-type dualities considered later.
The equivalence between $@Coh$ and $@DLat^@op$ was established by Stone [\TRDL].

Restricting to the full subcategories of ^{Stone spaces}, ^{Stone locales}, and ^{Boolean algebras}
yields the classical ^{Stone duality for Boolean algebras} (Stone [\RBA, \BRGT]),
expressed via the following chain of adjoint equivalences (\vStone, ^^!{Stone duality for Boolean algebras}):
$$@Stone\ltogets5{}{}@StoneLoc\ltogets5{}{}@BAlg^@op.$$
Here ^{Stone locales} and ^{Stone spaces} are defined as compact zero-dimensional locales (respectively, compact zero-dimensional ^{sober topological spaces}).

Restricting to (nonfull) subcategories of ^{Stonean spaces}, ^{Stonean locales}, and ^{complete Boolean algebras}
produces another chain of adjoint equivalences (\vStonean, ^^!{Stonean duality}), known as the ^{Stonean duality} (Stone [\ACSBR]):
$$@Stonean\ltogets5{}{}@StoneanLoc\ltogets5{}{}@CBAlg^@op.$$
Here ^{Stonean locales} and ^{Stonean spaces} are defined as ^{extremally disconnected} ^{Stone locales} (respectively ^{Stone spaces}).
$@CBAlg$ is the ^{category of complete Boolean algebras} with supremum-preserving homomorphisms as morphisms.

Finally, restricting again to full subcategories of ^{hyperstonean spaces}, ^{hyperstonean locales}, and ^{localizable Boolean algebras}
produces the chain of adjoint equivalences from \vmainthm\ (\vHStonean, ^^!{hyperstonean duality}):
$$@HStonean\ltogets5{}{}@HStoneanLoc\ltogets5{}{}@LBAlg^@op=@MLoc.$$

\label\reviewlocales
\subsection Review of locales

For an accessible introduction to locales, see Chapter~1 of Borceux [\HCA].

\proclaim Definition.
The ^={category of frames} $@=Frm$ is defined as follows.
A ^={frame[|s]} is a poset~$L$ that admits finite infima (alias meets, $∧$), arbitrary suprema (alias joins, $∨$),
and for any $a∈A$ the map $x↦x∧a$ ($L→L$) preserves suprema.
A ^={map[|s] of frames} $L→L'$ is a map of sets $L→L'$ that preserves finite infima and arbitrary suprema.
^^={homomorphism[|s] of frames}
^^={morphism[|s] of frames}

\proclaim Definition.
The ^={category of locales} $@=Loc$ is defined as $@Frm^@op$.
Given a ^{locale}~$L$, we refer to the elements of the corresponding ^{frame}
as the ^={open[s|]} of~$L$.
^^={locale[|s]}
^^={morphism[|s] of locales}
^^={map[|s] of locales}

\proclaim Remark.
To avoid confusion between maps of frames and maps of locales,
we use a superscript~$*$ for maps of frames, e.g., $f^*$, $g^*$, etc.

\proclaim Remark.
The passage to the opposite category in the definition of $@Loc$ is motivated by the desire to make $@Loc$ like $@Top$.
However, most of the actual computations are easier to perform in the category $@Frm$.
We automatically transport all notions defined for ^{frames} or ^{locales} to the opposite category.
For instance, below we define ^{open maps of frames}, which automatically yields ^{open maps of locales}.

\proclaim Definition.
The functor $$"=Ω:@Top→@Loc$$ from the category of topological spaces and continuous maps to the ^{category of locales}
sends a topological space $(X,U)$ to its poset of open sets~$(U,⊂)$ equipped with the partial order given by inclusion of subsets.
It sends a continuous map of topological spaces $f:(X,U)→(X',U')$ to the induced preimage map $f^*:U'→U$.

The functor~$"Ω$ becomes fully faithful once restricted to the full subcategory of ^{sober topological spaces}.

\proclaim Definition.
A ^={sober topological space[|s]}
^^={sober space[|s]}
^^={sober}
is a topological space for which the operation of closure establishes a bijection
from the set of points to the set of irreducible closed subsets, i.e.,
nonempty closed subsets that cannot be represented as a union of two proper closed subsets.

All Hausdorff topological spaces are sober, and many non-Hausdorff spaces commonly occurring in mathematics,
such as the Zariski spectrum of a commutative ring, are also sober.

\proclaim Definition.
^{Locales} in the image of the functor~$"Ω$ are known as ^={spatial}.
^^={spatial locale[|s]}

Thus, to pass from $@Top$ to~$@Loc$,
one must discard the nonsober topological spaces and add the nonspatial locales.
Arguably, nonsober spaces carry little practical value for much of mathematics,
whereas nonspatial locales are quite important.
In fact, below we will define ^{measurable locales}, and the only ^{measurable locales} that are spatial
are precisely the discrete locales, i.e., those obtained by applying $"Ω$ to a discrete topological space.
Thus, we are forced to use the pointfree localic formalism in order to formulate our main theorem.

The functor~$"Ω$ is a left adjoint functor that fits into the adjunction
$$@Top\ltogets7{"Ω}{"Sp}@Loc$$
between topological spaces and ^{locales}.
This adjunction restricts to an adjoint equivalence
$$@=Sober\ltogets7{"Ω}{"Sp}@=SpatialLoc$$
between the categories of sober spaces and ^{spatial locales}.

\proclaim Definition.
The right adjoint functor $$"=Sp="=pt:@Loc→@Top$$ sends a ^{locale}~$L$
to the topological space $"Sp(L)$
whose set of points~$X$ is the set of morphisms of locales $p:1→L$
and the collection of open sets is constructed as the image of the map of sets $f^*:L→2^X$, where $f$ is the morphism of locales $$f:∐_{p:1→L}1→L.$$

The essential image of~$"Sp$ consists precisely of sober spaces.

Much of general topology can be developed in the setting of locales,
see, for instance, Picado–Pultr [\FLoc],
which offers a full coverage of basic topics in general topology,
including compactness, local compactness, uniformity, paracompactness, completion, metrics, connectedness, real numbers, and localic groups.
Furthermore, proofs in the pointfree setting are often more elegant and clear,
since we do not have to construct sets of points, which eliminates substantial chunks of arguments.
Additional benefits include the elimination of the axiom of choice from most proofs,
which allows one to use locales in equivariant or fibered settings codified by arbitrary toposes, where topological spaces simply do not work.
We refer the reader to the survey articles by Johnstone [\PPtop] and [\ArtP] for more information.

\proclaim Remark.
Any ^{map of frames} $f^*:L'→L$ preserves suprema, so it admits a right adjoint $f_*:L→L'$.
For ^{spatial locales}, $f_*U$ is the largest open subset of~$L'$ whose preimage is inside~$U$.

\proclaim Remark.
If $f:L→L'$ is a ^{map of locales},
then we denote the corresponding homomorphism of frames by $f^*:L'→L$
and its right adjoint by $f_*:L→L'$.
If $f^*$ happens to admit a left adjoint, it is denoted by $f_!:L→L'$.
If $L=2^X$ and $L'=2^{X'}$ are discrete locales
and the ^{map of locales} $L→L'$ is induced by a map of sets $f:X→X'$,
then the notation $f_!⊣f^*$ is consistent with the same notation for images and preimages of subsets with respect to the map of sets~$f$.

\proclaim Definition.
A ^={sublocale[|s]} of a ^{locale}~$L$ is a ^{map of locales} $f:S→L$
such that the ^{map of frames} $f^*:L→S$ is surjective on the underlying sets.

In particular, a subspace $S⊂X$ of a topological space $(X,U)$ equipped with the induced topology yields a ^{sublocale},
because open sets in the induced topology are precisely the intersections of $S$ and an element of~$U$.
Conversely, if $S$ is a ^{sublocale} of the underlying ^{locale} of a ^{sober topological space}~$X$
and $S$ is a ^{spatial locale}, then $S$ is induced by a subspace of~$X$.
(The underlying locale of a topological space typically also has a lot of sublocales that are not spatial.)

\proclaim Definition.
The ^={Heyting implication[|s]} of opens $x,y∈L$ of a ^{locale}~$L$
is the unique open $x→y$ in~$L$ such that $w≤(x→y)$ if and only if $(w∧x)≤y$.

We have $x→y=\sup_{w∧x≤y}w$.
For ^{spatial locales}, $x→y$ is the interior of the union of~$y$ and the complement of~$x$.

\proclaim Definition.
The ^={pseudocomplement[|s]} of an open $x∈L$ of a ^{locale}~$L$ is $¬x≔(x→0)$.

Thus, $¬x=\sup_{w∧x=0}w$.
For ^{spatial locales}, $¬x$ is the interior of the complement of~$x$.

\proclaim Definition.
An ^={open map[|s] of frames} $f^*:L'→L$ is a ^{map of frames}
that preserves infima and the Heyting implication (meaning $f^*(m→n)=(f^*m→f^*n)$), i.e., it is a morphism of complete Heyting algebras.
^^={open map[|s]}
^^={open map[|s] of locales}

\proclaim Remark.
If a ^{map of locales} $f:L→L'$ is open, the ^{map of frames} $f^*:L'→L$ admits a left adjoint $f_!:L→L'$.
One can think of $f_!$ as sending an open in~$L$ to its image in~$L'$, which is again open because $f$ is an open map.
By a theorem of Joyal and Tierney (which is not needed below), one can equivalently define ^{open maps of locales} $f:L→L'$
by requiring that $f^*$ admits a left adjoint $f_!:L→L'$
such that the Frobenius reciprocity condition is satisfied: $f_!(m∧f^*n)=f_!m∧n$,
or, equivalently, $f_!(a→f^*b)=f_!a→b$.

\label\reflectopen
\proclaim Remark.
Any open map $f:X→Y$ of topological spaces is sent to an open map $"Ω f:"Ω X→"Ω Y$ of locales by the functor $"Ω:@Top→@Loc$.
Conversely,
if every point $y∈Y$ has an open neighborhood $U⊂Y$ such that $U∖\{y\}⊂Y$ is open
(e.g., $Y$ is a $\rm T_1$-space),
then whenever $"Ω f$ is an open map of locales, the map~$f$ is an open map of topological spaces.

\subsection Valuations on locales

\proclaim Definition.
A positive ^={valuation[|s]} on a distributive lattice~$L$ (e.g., a ^{locale} or a ^{Boolean algebra})
is a map $ν:L→[0,∞)$ such that $ν(0)=0$, $ν(x)+ν(y)=ν(x∨y)+ν(x∧y)$, and $x≤y$ implies $ν(x)≤ν(y)$.
A positive ^{valuation} is ^={continuous} if $ν$ preserves existing suprema of directed subsets.
^^={continuous valuation[|s]}
A ^={complex valuation[|s]} is a map $ν:L→`C$ such that $ν(0)=0$, $ν(x)+ν(y)=ν(x∨y)+ν(x∧y)$,
and we say that it is continuous if $ν(\sup_{x∈I} x)=\lim_{x∈I} ν(x)$ for any directed subset $I⊂L$.
^^={continuous complex valuation[|s]}
A ^={real valuation[|s]} is a complex valuation that factors through $`R⊂`C$.
A ^{valuation} is ^={faithful} if $ν(x)=0$ implies $x=0$.
^^={faithful continuous valuation[|s]}

The following definition mimics the usual definition of a ^{measure} as a countably additive map from a ^{σ-algebra},
but with countable additivity promoted to complete additivity along the lines of \vessentialcompletelyadditive.

\proclaim Definition.
A ^={completely additive valuation[|s]} $ν:L→[0,∞)$ on a distributive lattice~$L$
is a ^{valuation} $ν:L→[0,∞)$ such that
for any family $x:I→L$ consisting of pairwise disjoint elements of~$L$
(i.e., $x_i∧x_j=0$ whenever $i≠j$)
such that $\sup x$ exists,
we have $$ν\left(\sup_{i∈I} x_i\right)=∑_{i∈I} ν(x_i),$$ where the right sum converges absolutely.

For arbitrary ^{locales} ^{continuous valuations} have better theoretical properties than ^{completely additive valuations}.
However, for ^{Boolean locales} there is no difference.

\label\complvaladd
\proclaim Lemma.
Assuming the axiom of choice,
a map $ν:A→[0,∞)$ on a ^{Boolean locale} (i.e., a ^{complete Boolean algebra})~$A$ is a ^{continuous valuation}
if and only if it is ^{completely additive valuation}.

\proof Proof.
Suppose $ν:A→[0,∞)$ is a ^{continuous valuation}.
Then for any disjoint family $x:I→A$ we have $$ν\left(\sup_{i∈I} x_i\right)=ν\left(\sup_{K⊂I}\sup_{k∈K}x_k\right)
=\sup_{K⊂I}ν(\sup_{k∈K}x_k)
=\sup_{K⊂I}∑_{k∈K}ν(x_k)
=∑_{i∈I} ν(x_i),$$
where $K$ runs over all finite subsets of~$I$.
This shows that $ν$ is a ^{completely additive valuation}.
\ppar
Conversely, if $ν:A→[0,∞)$ is a ^{completely additive valuation},
i.e.,
$$ν\left(\sup_{i∈I} x_i\right)=∑_{i∈I} ν(x_i),$$
then substituting $I=∅$ yields $ν(0)=0$.
The identity $z=(z∧w)∨(z∧¬w)$
implies $ν(x)=ν(x∧y)+ν(x∧¬y)$ and $ν(x∨y)=ν(y)+ν(x∧¬y)$,
which together imply $ν(x)+ν(y)=ν(x∨y)+ν(x∧y)$.
If $x≥y$, then $x∧y=y$, so $ν(x)=ν(x∧y)+ν(x∧¬y)$ implies $ν(x)≥ν(y)$.
Thus, $ν$ is a ^{valuation}.
\ppar
Finally, to show that $ν$ is continuous, suppose that $S⊂A$ is a directed subset of~$A$.
We may assume that if $s∈S$ and $s'≤s$, then also $s'∈S$.
Using Zorn's lemma, choose a maximal disjoint family~$P$ that refines~$S$.
By construction, $\sup P≤\sup S$.
If $ρ=(\sup S)∧¬(\sup P)≠0$, then there is $s∈S$ such that $s∧ρ≠0$
and the family $P∪\{s∧ρ\}$ is a disjoint refinement of~$S$, contradicting the maximality of~$P$.
Thus, $ρ=0$ and $\sup P=\sup S$.
\ppar
Denote by $Q$ the closure of~$P$ under finite joins.
Since $S$ is directed, $Q$ refines~$S$,
so $∑_{p∈P} ν(p)=\sup_{q∈Q} ν(q)≤\sup_{s∈S}ν(s)$.
On the other hand, for any $s∈S$ we have
$s=\sup_{p∈P}s∧p$,
so $ν(s)=∑_{p∈P}ν(s∧p)≤∑_{p∈P}ν(p)$.
Thus, $ν(\sup S)=ν(\sup P)=∑_{p∈P}ν(p)=\sup_{s∈S}ν(s)$.

\proclaim Remark.
In the context of topos theory, where the axiom of choice can be false,
only ^{continuous valuations} allow for a good theory,
with the proviso that $`R$ refers to the lower reals.

\subsection Stone duality for distributive lattices

In this section, we briefly review the chain of adjoint equivalences (\vCoh, ^^!{Stone duality for distributive lattices})
$$@Coh\ltogets5{}{}@CohLoc\ltogets5{}{}@DLat^@op$$
between ^{coherent spaces}, ^{coherent locales}, and ^{distributive lattices}
that serves as a foundation for all Stone-type dualities considered later.
Complete proofs of all cited facts are given by Johnstone [\SS, §II.3], so we only recall the relevant definitions.

\proclaim Definition.
A ^={lattice[|s]} is a poset that admits all finite meets and finite joins.
Lattices form a category whose morphisms are maps of posets that preserve finite meets and finite joins.
A ^={distributive lattice[|s]} is a lattice~$R$ that satisfies the finite distributive law
$$x∧(y∨z)=(x∧y)∨(x∧z)$$ for all $x,y,z∈R$.
Distributive lattices form a full subcategory~$@=DLat$ of the category of lattices.

\proclaim Definition.
A ^={compact open[|s]} of a ^{locale}~$L$ is an open $a∈L$ such that for any $S⊂L$ with $⋁S≥a$ we have
$⋁F≥a$ for some finite $F⊂S$.
A ^{locale}~$L$ is ^={compact} if the open $1∈L$ is a ^{compact open}.
A topological space~$T$ is ^{compact} if the locale~$"Ω(T)$ is compact,
i.e., any open cover of~$T$ has a finite subcover.
A ^={coherent locale[|s]} is a ^{locale}~$L$ such that compact opens of~$L$
^^={coherent}
are closed under finite meets and any open of~$L$ is a join of compact opens.
Thus, coherent locales are automatically compact.
Coherent locales form a category $@=CohLoc$
whose morphisms are ^={coherent map[s|] of locales}, defined as ^{maps of locales} $f:L→L'$ such that $f^*$ preserves compact opens.

\proclaim Definition.
A ^={coherent space[|s]} is a ^{sober topological space}~$S$ such that $"Ω(S)$ is a coherent locale.
Coherent spaces are automatically compact.
Coherent spaces form a category $@=Coh$ whose morphisms are ^={coherent map[s|] of spaces},
defined as continuous maps $f:S→S'$ such that $"Ω(f)$ is a ^{coherent map of locales}.

The axiom of choice implies that ^{coherent locales} are ^{spatial} (Johnstone [\SS, Theorem~II.3.4]),
which immediately yields the following result.

\label\Coh
\proclaim Proposition.
The adjoint equivalence
$$@Sober\ltogets5{"Ω}{"Sp}@SpatialLoc$$
restricts to the adjoint equivalence
$$@Coh\ltogets5{"Ω}{"Sp}@CohLoc.$$

\proclaim Definition.
The set of ^{compact opens} of a ^{locale}~$L$ is a ^{distributive lattice} denoted by $"COpen(L)$.
^{Coherent maps of locales} preserve ^{compact opens} and their finite joins and meets,
so we get a functor $$"=COpen:@CohLoc→@DLat^@op.$$

\proclaim Definition.
An ^={ideal[|s]} of a ^{distributive lattice}~$R$ is a subset $I⊂R$ that is closed under finite joins
and such that $a∈I$ and $b≤a$ imply $b∈I$.
The set $"Idl(R)$ of all ideals of a ^{distributive lattice}~$R$
equipped with the inclusion ordering is a ^{frame} (Johnstone [\SS, Corollary~II.2.11]).
A morphism $f:R→R'$ of ^{distributive lattices}
induces a morphism $"Idl(f):"Idl(R)→"Idl(R')$
that sends an ideal $I⊂R$ to the ideal of~$R'$ generated by the image of~$I$ under~$f$.
The resulting functor $$"=Idl:@DLat→@Frm$$
is left adjoint to the forgetful functor $@Frm→@DLat$ (as shown in the cited result).

\label\implideal
\proclaim Remark.
The formula $x→y=\sup_{w∧x≤y}w$ yields the following formula for the Heyting implication in~$"Idl(A)$:
$$I→J=\{k∈A\mid kI⊂J\}.$$
In particular,
$$¬I=\{k∈A\mid ∀i∈I: ki=0\}.$$

\proclaim Proposition.
(^=:{Stone duality for distributive lattices}.)
There is an adjoint equivalence of categories
$$@CohLoc\ltogets9{"COpen}{"Ideal}@DLat^@op.$$
The functor $"=Ideal$ is obtained from the functor $"Idl^"op$ by restricting its codomain to $@CohLoc$.
For a ^{coherent locale}~$L$, the unit $L→"Ideal("COpen(L))$ is a ^{map of locales} whose underlying ^{map of frames}
sends an ideal $I⊂"COpen(L)$ to the join of its elements.
For a ^{distributive lattice}~$R$,
the counit $"COpen("Ideal(R))→R$ is a morphism in $@DLat^@op$
such that its opposite morphism in $@DLat$
is the map $R→"COpen^"op("Idl(R))$ that sends $r∈R$ to the principal ideal of~$R$ generated by~$r$ (i.e., $\{b∈R\mid b≤r\}$),
which is automatically a ^{compact open} of the ^{frame} $"Idl(R)$.

\proclaim Remark.
The composition $$@DLat^@op\lto9{"Ideal}@CohLoc\lto9{"Sp}@Coh$$ is known as the ^={Stone spectr[um|a]} functor.

\subsection Stone duality for Boolean algebras

In this section we review how the ^{Stone duality for distributive lattices}
restricts to the following chain of adjoint equivalences (\vStone, ^^!{Stone duality for Boolean algebras}):
$$@Stone\ltogets5{}{}@StoneLoc\ltogets5{}{}@BAlg^@op,$$
essentially due to Stone [\RBA, \BRGT].

\proclaim Definition.
A ^={Boolean algebra[|s]} (alias ^={Boolean ring[|s]}) is
a ^{distributive lattice}~$R$ in which all elements are ^={complemented}: for any $x∈R$ there is $¬x∈R$ such that $x∧¬x=0$ and $x∨¬x=1$.
^^={map[|s] of Boolean algebras}
^^={homomorphism[|s] of Boolean algebras}
The ^={category of Boolean algebras} $@=BAlg$ is a full subcategory of~$@DLat$.

Equivalently, one could define ^{Boolean algebras} as rings such that $x^2=x$ for all~$x$
and the ^{category of Boolean algebras} as a full subcategory of the category of rings.
In this definition,
^{Boolean algebras} are automatically commutative $`Z/2$-algebras and the ordering can be recovered by defining $x≤y ≡ (xy=y)$.
Recall that by our convention, all rings are unital and all homomorphisms of rings preserve units.

\proclaim Definition.
A ^{locale}~$L$ is ^={regular} ^^={regularity} ^^={regular locale[|s]} if for any $y∈L$ we have $y=⋁_{x≺y}x$, where $x≺y$ means $¬x∨y=1$.
The category $@=StoneLoc$ of ^={Stone locale[s|]} is
^^={Stone frame[|s]}
^^={category of Stone locales}
the full subcategory of~$@CohLoc$
consisting of ^{regular} ^{coherent locales}.

\proclaim Definition.
A ^={Stone space[|s]} is a ^{coherent space}~$T$ such that $"Ω(T)$ is a ^{regular locale} (hence a ^{Stone locale}).
The ^={categor[y|ies] of Stone spaces} $@=Stone$ is the full subcategory of~$@Coh$ consisting of Stone spaces.

\proclaim Remark.
Stone originally defined what is now known as ^{Stone spaces} using totally disconnected spaces.
A topological space is ^={totally disconnected} if its only connected subspaces are singletons.
Johnstone [\SS, Theorem~II.4.2] shows that the class of ^{Stone spaces} as defined above
coincides with the class of compact totally disconnected Hausdorff topological spaces,
as used by Stone himself.
Alternatively, ^{Stone spaces} can be characterized as ^{compact} zero-dimensional ^{sober spaces}.

The following proposition follows immediately from our definition of ^{Stone spaces} and ^{Stone locales}.

\label\Stone
\proclaim Proposition.
The adjoint equivalence
$$@Coh\ltogets5{"Ω}{"Sp}@CohLoc$$
restricts to an adjoint equivalence
$$@Stone\ltogets5{"Ω}{"Sp}@StoneLoc.$$

\label\copenclopen
\proclaim Lemma.
For a ^{Stone locale}~$L$,
the set of ^{compact opens}
coincides with the set of ^={clopen[|s]} (closed and open) elements,
^^={clopen subset[|s]}
i.e., ^{complemented} elements of the underlying ^{frame} of~$L$.

\proof Proof.
In a ^{compact} ^{locale}~$L$ any ^{clopen} $a∈L$ is a ^{compact open}:
if for some $S⊂L$ we have $⋁S≥a$, then $S'=S∪\{¬a\}$ satisfies $⋁S=1$, so by compactness of~$L$,
there is a finite $F'⊂S'$ such that $⋁F'=1$, therefore taking $F=F'∖\{¬a\}⊂S$ we get $⋁F≥a$.
\ppar
Conversely, if $a∈L$ is a ^{compact open},
then by ^{regularity} of~$L$ we have
$a=⋁_{x≺a}x$, so using compactness of~$a$ we can extract a finite subset $F⊂\{x∈L\mid x≺a\}$ such that $⋁F=a$.
Since the set $\{z∈L\mid z≺a\}$ is closed under finite joins, we see that $⋁F=a≺a$, i.e., $¬a∨a=1$, so $a$ is ^{clopen}.

\proclaim Remark.
^{Maps of frames} preserve ^{complemented} elements.
Thus, ^{maps of locales} between ^{Stone locales} are automatically coherent,
so $@StoneLoc$ is a full subcategory of~$@Loc$, unlike $@CohLoc$.
Likewise, $@Stone$ is a full subcategory of~$@Top$, unlike $@Coh$.

\proclaim Proposition.
(^=:{Stone duality for Boolean algebras}.)
^^={Stone duality}
The adjoint equivalence of categories
$$@CohLoc\ltogets9{"COpen}{"Ideal}@DLat^@op$$
restricts to an adjoint equivalence
$$@StoneLoc\ltogets9{"COpen}{"Ideal}@BAlg^@op.$$

\proof Proof.
It suffices to show that $"COpen$ sends ^{Stone locales} to ^{Boolean algebras}
and $"Ideal$ sends ^{Boolean algebras} to ^{Stone locales}.
Suppose $L$ is a ^{Stone locale}.
By \vcopenclopen, ^{compact opens} in~$L$ are precisely the ^{clopen} elements.
In particular, every ^{compact open} in~$L$ is complemented; hence the ^{distributive lattice} $"COpen(L)$ is a ^{Boolean algebra}.
Conversely, suppose $A$ is a ^{Boolean algebra}.
Any principal ideal $(a)⊂A$ is a ^{clopen} (hence ^{compact open}) element of $"Ideal(A)$, with the complement $¬(a)=(1-a)$.
For any open $y∈L$ and a ^{clopen} $x∈L$ such that $x≤y$ we have $x≺y$ because $¬x∨y≥¬x∨x=1$.
Thus, for any $I∈"Ideal(A)$ we have $I=⋁_{a∈I}(a)$, where $(a)≺I$; hence the ^{coherent locale} $"Ideal(A)$ is ^{regular}.

\subsection Stonean duality

Next, we review how the ^{Stone duality for Boolean algebras} restricts to a chain of adjoint equivalences
(\vStonean, ^^!{Stonean duality})
$$@Stonean\ltogets5{}{}@StoneanLoc\ltogets5{}{}@CBAlg^@op$$
between the (nonfull) subcategories of ^{Stonean spaces}, ^{Stonean locales}, and ^{complete Boolean algebras},
which we refer to as the ^{Stonean duality}.
Stone's 1937 paper [\ACSBR] shows that ^{complete Boolean algebras}
correspond precisely to ^{extremally disconnected} ^{Stone spaces}
under Stone duality.

\proclaim Definition.
A ^{Boolean algebra} is ^={complete}
^^={complete Boolean algebra[|s]}
if any set of its elements admits a supremum.
A ^={morphism[|s] of complete Boolean algebras}
^^={complete homomorphism[|s]}
is a complete (i.e., continuous, or supremum-preserving) Boolean homomorphism.
The ^={category of complete Boolean algebras} $@=CBAlg$
is a (nonfull) subcategory of the ^{category of Boolean algebras}.

A crucial distinguishing feature of Stonean duality is that both categories involved in the adjunction
are subcategories of the ^{category of locales}.
We formulate this as a lemma.

\label\cbaloc
\proclaim Lemma.
The functor $@CBAlg^@op→@Loc$ defined as the opposite of the forgetful functor $@CBAlg→@Frm$ is a fully faithful functor.
Locales~$L$ in its essential image $@=BLoc$ are precisely ^={Boolean locale[s|]},
i.e., locales for which $a∨¬a=1$ for any $a∈L$,
or, equivalently, $¬¬a=a$ for any $a∈L$.
Any map of Boolean locales is automatically open.

\proof Proof.
If $¬¬a=a$ for all $a∈L$, then $1=¬¬a∨¬a=a∨¬a$ for all $a∈L$.
Conversely, if $a∨¬a=1$ for all $a∈L$, then $¬¬a=¬¬a∧(a∨¬a)=(¬¬a∧a)∨(¬¬a∧¬a)=(¬¬a∧a)∨0=¬¬a∧a=a$ for all $a∈L$.
This shows that the two conditions in the statement are equivalent.
\ppar
First, any $A∈@CBAlg$ is indeed a ^{frame}:
being a ^{complete Boolean algebra}, it admits arbitrary suprema and infima,
and the map $s↦s∧a$ preserves suprema for any $a∈A$
because it admits a right adjoint $t↦(¬a)∨t$.
\ppar
Secondly, any homomorphism $f:A→A'$ in $@CBAlg$ preserves arbitrary suprema and infima,
so is indeed a ^{map of frames}.
Conversely, if $f:A→A'$ preserves finite infima and arbitrary suprema,
then it is a ^{complete homomorphism}, so the functor is fully faithful.
Furthermore, $f$ preserves the Heyting implication, since in a ^{Boolean algebra}
$x→y=¬x∨(x∧y)=1+x+xy$, and all operations on the right side are preserved by~$f$.
Thus, $f$ preserves arbitrary infima and the Heyting implication,
so $f$ is an ^{open map of frames}.
\ppar
Finally, any ^{complete Boolean algebra} is a ^{Boolean locale} since $¬a=1-a$ in this case, so $a∨¬a=a∨(1-a)=1$.
Conversely, to show that any ^{Boolean locale} has a ^{complete Boolean algebra} as its underlying ^{frame},
we observe that the condition $¬¬a=a$ implies that the underlying Heyting algebra is a ^{Boolean algebra},
and hence also a ^{complete Boolean algebra}.

\proclaim Definition.
A ^{locale}~$L$ is ^={extremally disconnected} if $¬x∨¬¬x=1$ for any $x∈L$,
i.e., for any $x∈L$ the element $¬x∈L$ is ^{clopen}.
A ^={Stonean locale[|s]} is an ^{extremally disconnected} ^{Stone locale}.
^^={Stonean frame[|s]}
The ^={categor[y|ies] of Stonean locales} $@=StoneanLoc$ has ^{Stonean locales} as objects and ^{open maps of locales} as morphisms.
It is a (nonfull) subcategory of the ^{category of Stone locales}.

If we apply this definition to a locale that comes from a topological space,
the condition $¬x∨¬¬x=1$
boils down to saying that the closure~$\bar x$ is its own interior, i.e., the closure of any open~$x$ is open.

\proclaim Definition.
A topological space~$T$ is ^{extremally disconnected} if
^^={extremally disconnected space[|s]}
$"Ω(T)$ is an ^{extremally disconnected} locale.
A ^={Stonean space[|s]} is an extremally disconnected ^{Stone space},
equivalently, a ^{sober topological space}~$T$ such that $"Ω(T)$ is a ^{Stonean locale}.
The ^={categor[y|ies] of Stonean spaces} $@=Stonean$ has ^{Stonean spaces} as objects and open continuous maps as morphisms.
It is a (nonfull) subcategory of the ^{category of Stone spaces}.

A result of Edwards [\FLat, Theorem~5.1] shows that (in our notation)
morphisms in the image of the composition $"Sp∘"Ideal∘"ProjLoc$
are precisely the open maps, which motivates the choice of ^{open maps} as morphisms between ^{Stonean locales}.

Our definition of ^{Stonean locales} and ^{Stonean spaces} together with \vreflectopen\ immediately implies the following.

\label\Stonean
\proclaim Proposition.
The adjoint equivalence
$$@Stone\ltogets5{"Ω}{"Sp}@StoneLoc$$
restricts to an adjoint equivalence
$$@Stonean\ltogets5{"Ω}{"Sp}@StoneanLoc.$$

\proclaim Proposition.
(^=:{Stonean duality}.)
The adjoint equivalence
$$@StoneLoc\ltogets9{"COpen}{"Ideal}@BAlg^@op$$
restricts to an adjoint equivalence
$$@StoneanLoc\ltogets9{"COpen}{"Ideal}@CBAlg^@op=@BLoc.$$

\proof Proof.
It suffices to show that both functors preserve the given subcategories.
Johnstone [\SS, Lemma III.3.5] shows that both functors preserve objects of the given subcategories.
We include a short proof of this below together with a proof of preservation of morphisms.
\ppar
We work with the opposite adjunction
$$@CBAlg\ltogets9{"Idl}{"COpen^"op}@=StoneanFrm.$$
An element $x∈S$ of an ^{extremally disconnected} frame~$S$ is ^{clopen} if and only if $x=¬¬x$.
Furthermore, $¬¬x$ is the smallest ^{clopen} element greater than or equal to~$x$.
If $S∈@StoneanFrm$, then the ^{Boolean algebra} $"COpen^"op(S)$ is complete
because $\sup R=¬¬⋁R$ exists for any $R⊂"COpen^"op(S)$.
The supremum on the left is computed in the ^{Boolean algebra} $"COpen^"op(S)$,
whereas the join on the right is computed in the poset~$S$.
\ppar
Recall that a ^{map of frames} $f:S→S'$ is ^{open[| map]} if it preserves infima and Heyting implications, and hence also negations.
If $f:S→S'$ is an ^{open map} of ^{Stonean frames},
the induced map of Boolean algebras $$"COpen^"op(f):"COpen^"op(S)→"COpen^"op(S')$$
preserves suprema:
given $R⊂"COpen^"op(S)$, we have $\sup R=¬¬⋁R$,
so $$f(\sup R)=f(¬¬⋁R)=¬¬f(⋁R)=¬¬⋁f(R)=\sup f(R).$$
This shows that the functor $"COpen^"op$ restricts to the claimed subcategories.
\ppar
For any $A∈@CBAlg$, the ^{Stone frame} $"Idl(A)$ is ^{extremally disconnected} because $¬¬I∨¬I=A$ for any ideal $I⊂A$.
Indeed, according to \vimplideal, $¬I=\{k∈A\mid ∀i∈I: k∧i=0\}$.
Since $A$ is complete, we can rewrite this expression as $¬I=\{k∈A\mid k∧\sup I=0\}$,
i.e., $¬I$ is the principal ideal of~$A$ generated by the element $1-\sup I$.
Now $¬¬I$ is the principal ideal of~$A$ generated by the element $\sup I$.
Therefore, $¬I∨¬¬I=A$.
\ppar
Finally, the functor $"Idl$ sends a ^{complete homomorphism} $g:A→A'$ of ^{complete Boolean algebras}
to an ^{open map of frames} $"Idl(g):"Idl(A)→"Idl(A')$.
First, $"Idl(g)$ preserves infima;
the nontrivial part is $$\inf_k "Idl(g)(I_k)≤"Idl(g)(\inf_k I_k).$$
Take an element~$x'$ from the left side and recall that infima of ideals coincide with intersections,
i.e., for any $k$ we have $x'∈"Idl(g)(I_k)$, which means that there is $y_k∈I_k$ such that $x'≤g(y_k)$.
Now $x'≤\inf_k g(y_k)=g(\inf_k y_k)$ by completeness of~$g$, so $x'∈"Idl(g)(\inf_k I_k)$.
Secondly, $"Idl(g)$ preserves ^{Heyting implications};
the nontrivial part is $$("Idl(g)(I)→"Idl(g)(J))≤"Idl(g)(I→J).$$
Take an element~$x'$ from the left side;
by \vimplideal\ we have $x'∧"Idl(g)(I)≤"Idl(g)(J)$,
so for any $i∈I$ there is $j∈J$ such that
$$x'∧g(i)≤g(j)\iff x'≤(g(i)→g(j))=g(i→j) \iff g_!(x')≤i→j \iff g_!(x')∧i≤j.$$
Thus, $g_!(x')∧I≤J$, so $g_!(x')∈I→J$.
Now $x'≤g(g_!(x'))$ by adjunction $g_!⊣g$, so $x'∈"Idl(g)(I→J)$.

\proclaim Remark.
According to \vcbaloc,
both sides of the restricted adjunction $@StoneanLoc\togets@CBAlg^@op$ are subcategories of~$@Loc$,
with the subcategory $@CBAlg^@op=@BLoc$ being a full subcategory
and the subcategory $@StoneanLoc$ being a nonfull subcategory (its morphisms are ^{open maps of locales}).
As shown there, morphisms of $@CBAlg^@op$ are automatically open,
so both sides are full subcategories of the category of locales and open maps.
The two sides are drastically different in terms of their relation to topological spaces:
whereas in presence of the axiom of choice all ^{Stonean locales} are spatial, i.e., correspond
to topological spaces, the locales in $@CBAlg^@op=@BLoc$ are never spatial unless they are complete atomic Boolean algebras,
i.e., Boolean algebras of subsets of a given set~$X$.
In fact, any ^{complete Boolean algebra} splits as a product
of a complete atomic Boolean algebra and a complete atomless Boolean algebra.
The latter has no points at all when interpreted as a locale.
Another difference is that ^{Stonean locales} are always compact,
whereas the only compact locales in $@CBAlg^@op$ are ^{Boolean algebras} of subsets of a finite set.

\label\compactification
\proclaim Remark.
^{Compact opens} of any ^{Stonean frame}~$S$ are characterized by the property $s=¬¬s$,
so we have $"COpen^"op(S)=¬¬S$, where for any ^{locale}~$L$ we define $$¬¬L=\{s∈L\mid s=¬¬s\}=\{¬¬s\mid s∈L\}.$$
The double negation map $¬¬:L→¬¬L$ is a surjective homomorphism of ^{frames}
and we have a canonical embedding of ^{locales} $¬¬L→L$.
The sublocale $¬¬L$ of~$L$ is also known as the ^={double negation sublocale[|s]} of~$L$
and is the largest sublocale of~$L$ that is a ^{Boolean locale}.
It can also be characterized as the intersection of all dense sublocales of~$L$,
where a ^={dense sublocale} is a ^{sublocale} that is not contained in any closed ^{sublocale} other than~$L$.

\subsection Hyperstonean duality

Finally, we show how the ^{Stonean duality} restricts to a chain of adjoint equivalences from \vmainthm\ (\vHStonean, ^^!{hyperstonean duality})
$$@HStonean\ltogets5{}{}@HStoneanLoc\ltogets5{}{}@LBAlg^@op=@MLoc$$
between the full subcategories of ^{hyperstonean spaces}, ^{hyperstonean locales}, and ^{localizable Boolean algebras},
which we refer to as the ^{hyperstonean duality}.
^{Hyperstonean spaces} were introduced and studied by Dixmier [\HS],
who in particular shows [\HS, Théorème~2] that ^{hyperstonean spaces} are precisely the ^{Gelfand spectra} of ^{commutative von Neumann algebras},
equivalently, ^{Stone spectra} of ^{localizable Boolean algebras}.

\proclaim Definition.
A ^={localizable Boolean algebra[|s]} is a ^{Boolean algebra}~$A$
^^={localizable}
^^={localizability}
that is ^{complete} and in which the element $1∈A$ equals the supremum of all $a∈A$
such that the ^{complete Boolean algebra} $aA$ admits a (finite) ^{faithful continuous valuation} $aA→[0,∞)$.
The ^={category of localizable Boolean algebras} $@=LBAlg$
is a full subcategory of the ^{category of complete Boolean algebras} $@CBAlg$.
^^={map[|s] of localizable Boolean algebras}
^^={morphism[|s] of localizable Boolean algebras}

\label\lmcrit
\proclaim Lemma.
A ^{complete Boolean algebra}~$A$ is ^{localizable}
if and only if for any nonzero $b∈A$ there is a ^{continuous valuation} $μ:A→[0,∞)$
such that $μ(b)≠0$.

\proof Proof.
Suppose $A$ is ^{localizable}, i.e.,
the supremum of all $a∈A$ such that the ^{complete Boolean algebra} $aA$ admits a (finite) ^{faithful continuous valuation}
equals $1∈A$.
For any nonzero $b∈A$ we have $b=1⋅b=(\sup_a a)b=\sup_a(ab)$,
so one of $ab$ must be nonzero.
Fixing such an~$a$, there is a faithful ^{continuous valuation} $μ:aA→[0,∞)$.
It extends to a ^{continuous valuation} $ν:A→[0,∞)$ via the formula $ν(c)=μ(ac)$.
We have $ν(b)=μ(ab)≠0$ since $μ$ is faithful on $aA$.
\ppar
Conversely, if for any nonzero $b∈A$ there is a ^{continuous valuation} $μ:A→[0,∞)$
such that $μ(b)≠0$,
then denote by $c$ the supremum of all $a∈A$ such that the ^{complete Boolean algebra} $aA$ admits a ^{faithful continuous valuation}.
If $c≠1$, then there is a ^{continuous valuation} $μ:A→[0,∞)$ such that $μ(1-c)≠0$.
By continuity of~$μ$ there is a maximum $d∈A$ such that $μ(d)=0$, so $μ$ is faithful on $(1-d)A$ and vanishes on~$dA$.
By definition of~$c$, we must have $1-d≤c$,
but then $$μ(1-c)=μ((1-c)(1-d+d))=μ((1-c)(1-d))+μ((1-c)d)=μ((1-c)(1-d))≤μ((1-c)c)=0,$$ which contradicts $μ(1-c)≠0$.

\proclaim Proposition.
The category $@LBAlg$ admits all small products
and the forgetful functor $@LBAlg→@BAlg$ preserves small products.

\proof Proof.
^{Complete Boolean algebras} are closed under small products in the category $@BAlg$
because suprema in a product $∏_{i∈I}A_i$ of ^{complete Boolean algebras}~$A_i$ can be computed
indexwise with respect to~$I$.
^{Localizable Boolean algebras} are closed under small products in the category $@BAlg$
by \vlmcrit: if $b∈∏_{i∈I}A_i$ is nonzero, then for at least one $j∈I$ the projection $b_j∈A_j$ is nonzero.
Since $A_j$ is ^{localizable}, we can find a ^{continuous valuation} $μ_j:A_j→[0,∞)$
such that $μ_j(b_j)≠0$.
Composing $μ_j$ with the projection map $∏_{i∈I}A_i→A_j$,
we get a ^{continuous valuation} $μ:∏_{i∈I}A_i→[0,∞)$ such that $μ(b)≠0$.
Thus, $∏_{i∈I}A_i$ is ^{localizable}.

We now define the category of ^{measurable locales},
which is in the same relation to traditional measurable spaces
as ^{locales} are to traditional topological spaces.
Just like for ^{locales} we express everything in terms of the complete lattice of open sets,
for ^{measurable locales} we express everything in terms of the complete lattice of equivalence
classes of measurable sets modulo negligible sets.

\proclaim Definition.
The category $@=MLoc$ of ^={measurable locale[s|]}
is the full subcategory of the category $@Loc$ of ^{locales}
consisting of ^{localizable Boolean algebras}.
^^={category of measurable locales}
^^={morphism[|s] of measurable locales}
^^={map[|s] of measurable locales}

\proclaim Remark.
According to \vcbaloc, the inclusion $@CBAlg^@op→@Loc$ is a fully faithful functor,
and hence so is the composition $@LBAlg^@op→@CBAlg^@op→@Loc$.
In particular, $@MLoc$ is nothing else than $@LBAlg^@op$,
the opposite ^{category of localizable Boolean algebras}.
Also, $@MLoc$ is a full subcategory of the category~$@BLoc$ of ^{Boolean locales}
since $@LBAlg$ is a full subcategory of $@CBAlg$.

\proclaim Definition.
A ^={normal valuation[|s]} on a ^{locale}~$L$
is a ^{continuous valuation}~$ν$
such that $ν(¬¬a)=ν(a)$ for any $a∈L$.

\label\hnormal
\proclaim Lemma.
For any ^{Stonean locale}~$L$,
the restriction along the map of posets $"COpen(L)→L$
establishes a bijection from the set of ^{normal valuations} on~$L$
to the set of ^{continuous valuations} on~$"COpen(L)$.
The inverse map is given by precomposing with the ^{map of frames} $¬¬:L→"COpen(L)$.

\proof Proof.
The inclusion $"COpen(L)→L$ preserves finite meets and finite joins,
so any ^{valuation}~$μ$ on~$L$ restricts to a valuation~$ν$ on~$"COpen(L)$.
Given a directed subset $S⊂"COpen(L)$, we have $\sup S=¬¬⋁S$, where the supremum
on the left is taken in the ^{complete Boolean algebra} $"COpen(L)$,
whereas the join on the right is taken in the ^{locale}~$L$.
Given a ^{normal valuation}~$μ$ on~$L$, we have $μ(¬¬⋁S)=μ(⋁S)=\lim_{s∈S}ν(s)$,
which proves that the restriction~$ν$ of~$μ$ to $"COpen(L)$ is a ^{continuous valuation}.
\ppar
Conversely, given a ^{continuous valuation}~$ν$ on~$"COpen(L)$,
precomposing it with the ^{map of frames} $¬¬:L→"COpen(L)$
yields a ^{continuous valuation} on~$L$ that is normal by construction.
\ppar
The ^{map of frames} $¬¬:L→"COpen(L)$ is identity on $"COpen(L)$,
which immediately implies that the two constructed maps are inverse to each other,
so we indeed have a bijection between ^{normal valuations} on~$L$ and ^{continuous valuations} on $"COpen(L)$.

\label\hahnjordan
\proclaim Lemma.
(^=:{Hahn–Jordan decomposition of normal valuations}.)
Given a ^{normal valuation} $ν:L→`R$ on a ^{Stonean locale}~$L$,
define $ν_+$ respectively $ν_-$
as the supremum of all $a∈L$ such that $ν$ respectively $-ν$ restricted to $\{b∈L\mid b≤a\}$ is a positive ^{normal valuation}.
Then the above condition on $a∈L$ is equivalent to $a≤ν_+$ respectively $a≤ν_-$, in particular, both suprema are also maxima.
Furthermore, $ν_+∨ν_-=1$ (^={Hahn decomposition[|s]} of~$ν$)
and $ν_+∧ν_-$ is the maximal element $a∈L$ such that $ν$ vanishes on $\{b∈L\mid b≤a\}$.
Define $ν^+,ν^-:L→[0,∞)$ by setting $ν^+(x)=ν(ν_+∧x)$ and $ν^-(x)=-ν(ν_-∧x)$.
Then $ν^+$ and $ν^-$ are positive ^{normal valuations} on~$L$ and $ν=ν^+-ν^-$ (^={Jordan decomposition[|s]} of~$ν$).
We also set $|ν|=ν^++ν^-$ and $‖ν‖=|ν|(1)$, the latter known as the ^={variation norm} of~$ν$.
The above statements also hold if $ν:L→`R$ is a ^{continuous valuation} on a ^{Boolean locale}~$L$.
(^{Continuous valuations} on ^{Boolean locales} are automatically ^{normal valuations}.)
For complex valuations, we get a canonical decomposition into four summands.

\proof Proof.
\vhnormal\ establishes a bijection between ^{normal valuations} on a ^{Stonean locale}~$L$
and ^{continuous valuations} on the ^{complete Boolean algebra} $"COpen(L)$.
Since $"COpen(L)⊂L$ and $ν_+,ν_-∈"COpen(L)$ because $ν$ is a ^{normal valuation},
it suffices to treat the case of ^{continuous valuations} on ^{Boolean locales}.
\ppar
Proceeding exactly as in the proof of the classical Hahn–Jordan decomposition (see, for example, Fremlin [\MTii, Theorem~231E], whose proof we follow here),
we first show that $ν$ must be bounded.
Below, $a∖b$ means $a∧(1-b)$.
Set $N(x)=\sup\{|ν(y)|\mid y≤x\}$.
If $N(1)=∞$, construct a sequence $\{E_n\}_{n≥0}$ such that $N(E_n)=∞$ as follows.
Set $E_0=1∈L$.
Having constructed $E_n∈L$ with $N(E_n)=∞$,
choose $F_n≤E_n$ such that $|ν(F_n)|≥|ν(E_n)|+1$, then also $|ν(E_n∖F_n)|≥1$.
Since $N(E_n)=∞$, we must have $N(F_n)=∞$ or $N(E_n∖F_n)=∞$.
In the former case, set $E_{n+1}=F_n$ and in the latter, set $E_{n+1}=E_n∖F_n$.
We have $N(E_{n+1})=∞$ and $|ν(E_n∖E_{n+1})|≥1$.
The sequence $\{E_n∖E_{n+1}\}_{n≥0}$ is disjoint,
so $ν(⋁_{n≥0}E_n∖E_{n+1})=∑_{n≥0}ν(E_n∖E_{n+1})$ is finite,
which implies $ν(E_n∖E_{n+1})→0$ as $n→∞$, and yet $|ν(E_n∖E_{n+1})|≥1$ for all $n≥0$, a contradiction.
Thus, $N(1)<∞$ and $ν$ is bounded.
\ppar
Next, we show that $ν$ attains its supremum~$γ$.
Choose $E_n∈L$ such that $ν(E_n)≥γ-2^{-n}$, for all $n≥0$.
Set $F_{m,n}=⋀_{m≤i≤n}E_i$.
We have $ν(F_{m,n})≥γ-2^{1-m}+2^{-n}$ (shown by induction on~$n$).
Set $G_m=⋀_{n≥m}F_{m,n}$.
We have $ν(G_m)≥γ-2^{1-m}$ because $F_{m,n}≥F_{m,n+1}$ for all $m≥0$.
Set $H=⋁_{m≥0}G_m$.
Since $G_m≤G_{m+1}$, we have $ν(H)≥γ$, so $ν(H)=γ$.
For any $K≤H$ we have $ν(H)-ν(K)=ν(H∖K)≤γ=ν(H)$, so $ν(K)≥0$.
For any $K≤1-H$ we have $ν(H)+ν(K)=ν(H∨K)≤γ=ν(H)$, so $ν(K)≤0$.
\ppar
Consider the element $Q=1-(ν_+∨ν_-)$.
The above construction yields some $H≤Q$ on which $ν$ attains its supremum on $\{H\mid H≤Q\}$ and for any $K≤H$ we have $ν(K)≥0$.
Thus, $H=0$ since otherwise $ν_+$ could be replaced by $ν_++H$, yielding a contradiction with the definition of~$ν_+$.
Hence, $ν(K)≤0$ for all $K≤Q$.
Likewise, $ν(K)≥0$, and hence $ν(K)=0$ for all $K≤Q$.
Thus, $Q=0$ since otherwise $ν_+$ could be replaced by $ν_++Q$, yielding a contradiction with the definition of~$ν_+$.
Hence, $ν_+∨ν_-=1$.
\ppar
Both $ν^+$ and $ν^-$ are ^{continuous valuations} on~$L$ because $x↦ν_+∧x$ and $x↦ν_-∧x$ are morphisms of lattices that preserve suprema of directed subsets.
Also, $ν^+$ and $ν^-$ are positive because $ν^+(x)=ν(ν_+∧x)=\sup_a ν(a∧x)≥0$,
where $a$ runs over all elements of~$L$ such that $ν$ restricted to $\{b∈L\mid b≤a\}$ is a positive ^{continuous valuation}.
\ppar
We have $$ν(x)=ν((ν_+∨ν_-)∧x)=ν((ν_+∧x)∨(ν_-∧x))=ν(ν_+∧x)+ν(ν_-∧x)-ν(ν_+∧ν_-∧x)=ν^+(x)-ν^-(x),$$
which establishes the Jordan decomposition of~$ν$.

\proclaim Corollary.
Complex ^{normal valuations} on a ^{locale}~$L$
are precisely elements in the linear span of positive ^{normal valuations} on~$L$.

\proclaim Definition.
The category~$@=HStoneanLoc$ of ^={hyperstonean locale[s|]} is the full subcategory of $@StoneanLoc$
consisting of ^{Stonean locales}~$L$ that admit sufficiently many ^{normal valuations},
meaning that for any $a∈L∖\{0\}$ there is a ^{normal valuation} $ν:L→[0,∞)$ such that $ν(a)≠0$.

\proclaim Definition.
A finite ^={normal measure[|s]} $μ$ on a ^{Stonean space}~$S$
is a ^={Borel measure[|s]}~$μ$ on~$S$
(i.e., a ^{measure} on the σ-algebra generated by open subsets of~$S$)
whose restriction to open subsets of~$S$
is a ^{normal valuation} on~$"Ω(S)$.

Recall that a subset~$A$ of a topological space~$X$
is ^={rare} ^^={rare subset[|s]} ^^={rare set[|s]} (alias ^={nowhere dense}) if the interior of the closure of~$A$ is empty
and ^={meager} ^^={meager subset[|s]} ^^={meager set[|s]} (alias ^={first category}) if it is a countable union of rare sets.

\proclaim Remark.
The restriction of a finite Borel measure on~$S$ to $"Ω(S)$ is automatically a ^{valuation}.
The condition that this ^{valuation} is ^{continuous} is also known as τ-additivity, τ-regularity, or τ-smoothness in measure theory.
On a compact regular topological space, finite τ-smooth Borel measures coincide with finite Radon measures,
so we define ^{normal measures} (not necessarily finite) on a ^{Stonean space}~$S$ as Radon measures that vanish on all ^{rare subsets} of~$S$.
With such a reformulation it is clear that a ^{normal measure} on~$S$
yields an ordinary measure on the ^{enhanced measurable space} $"TM(S)$ constructed in \vTMdef.

\proclaim Remark.
For a ^{Boolean algebra}~$A$, finitely additive maps $A→[0,∞)$
are in a canonical bijective correspondence with (positive finite) Radon measures on the ^{Stone space} $"Sp("Ideal(A))$
by Fremlin [\MTiv, Proposition~416Q].
If $A$ is a ^{complete Boolean algebra},
then by \vcomplvaladd, ^{completely additive valuations} $A→[0,∞)$ are precisely ^{continuous valuations} $A→[0,∞)$,
by \vhnormal, ^{continuous valuations} on~$A$ are in a canonical bijective correspondence with ^{normal valuations} on $"Ideal(A)$,
and by \vnormalmeasurevaluation, ^{normal valuations} on $"Ideal(A)$ are in a canonical bijective correspondence with ^{normal measures} on $"Sp("Ideal(A))$.
Thus, the condition that a (positive finite) Radon measure on a ^{Stonean space}~$S$ vanishes on all ^{rare subsets},
i.e., is ^{normal[| measure]},
is precisely equivalent to the condition that the corresponding ^{valuation} on $"COpen("Ω(S))$ is ^{completely additive[| valuation]}, equivalently, ^{continuous}.

We now show that our definition of ^{normal measures} on ^{Stonean spaces}
is equivalent to Dixmier's definition [\HS, Definition~1].
This lemma is not used anywhere else.

\proclaim Lemma.
The set of ^{normal measures} on a ^{Stonean space}~$S$
coincides with the set
of positive Radon measures~$μ$ such that the integration map $C(S,`R)→`R$ with respect to~$μ$ preserves suprema of bounded directed subsets.

\proof Proof.
According to Dixmier [\HS, Proposition~1],
a measure~$μ$ is
a positive Radon measure for which the integration map $C(S,`R)→`R$ with respect to~$μ$ preserves suprema of bounded directed subsets
if and only if $μ$ is a positive Radon measure that vanishes on all ^{rare} (alias ^{nowhere dense}) subsets of~$S$.
Suppose $μ$ satisfies these equivalent conditions.
Denote by~$ν$ the restriction of~$μ$ to $"Ω(S)$, which is a ^{valuation} because $μ$ is finitely additive.
Radon measures are τ-smooth, which means precisely that $ν$ is a ^{continuous valuation}.
Given an open subset $U⊂S$, the set $\bar U∖U$ is ^{rare}, so $μ(\bar U)=μ(U)$,
which shows that $ν(¬¬a)=ν(a)$ for any $a∈L$ if we recall that $\bar U=¬¬U$ for any $U∈"Ω(S)$.
Thus, $ν$ is a ^{normal valuation} on~$"Ω(S)$.
\ppar
Conversely, consider a (finite positive) Borel measure~$μ$ on~$S$
whose restriction to the open subsets of~$S$
is a ^{normal valuation} on~$"Ω(S)$.
First, $μ$ vanishes on all ^{rare subsets} $R⊂S$
because $μ(S)=μ(S∖\bar R)$
since $¬¬(S∖\bar R)=S$, so $μ(\bar R)=0$ and $μ(R)=0$.
Secondly, the Borel measure~$μ$ is τ-smooth by definition of a ^{continuous valuation}.
Every τ-smooth Borel measure on a compact topological space is a Radon measure, which completes the proof.

\proclaim Definition.
A ^={hyperstonean space[|s]} is a ^{sober space}~$S$ such that $"Ω(S)$ is a ^{hyperstonean locale}.
The ^={category of hyperstonean spaces} $@=HStonean$ is a full subcategory of the ^{category of Stonean spaces}.
^^={morphism[|s] of hyperstonean spaces}

\label\meagerrare
\proclaim Lemma.
In a ^{hyperstonean space}~$S$,
the class of ^{rare subsets},
the class of ^{meager subsets},
and the class of measurable subsets (i.e., symmetric differences of open and ^{meager} subsets)
on which every ^{normal measure} vanishes
all coincide.

\proof Proof.
(See Dixmier [\HS, Proposition~5].)
All ^{rare subsets} are ^{meager} and ^{normal measures} vanish on all ^{meager subsets} by definition.
Suppose $A⊂S$ is a measurable subset on which all normal measures vanish.
We have to show that $A$ is ^{rare}.
The set $\bar A∖A$ is ^{rare}, so all ^{normal measures} vanish on $\bar A$,
and hence also on the ^{clopen} subset $\bar A^\circ$.
This means that $\bar A^\circ=∅$, i.e., $\bar A$ is ^{rare}, so $A$ is also ^{rare}.

\label\hyperstoneanmeasurable
\proclaim Lemma.
Given a ^{hyperstonean space}~$S$,
the following classes of sets (henceforth known as {\it measurable sets\/}) coincide:
\li symmetric differences of ^{clopen} sets and ^{rare sets};
\li symmetric differences of closed sets and ^{rare sets};
\li symmetric differences of open sets and ^{rare sets};
\li unions of open sets and ^{rare sets} (can be assumed to be disjoint);
\li differences of closed sets and ^{rare sets} (the former can be assumed to contain the latter).
\endlist
The class~$M$ of measurable sets is a ^{σ-algebra} of subsets of~$S$ (^!{σ-algebra}).
Every element~$X$ of~$M$ admits a unique presentation as the symmetric difference of a ^{clopen} set and a ^{rare set}.
Furthermore, for any measurable subset $X$ of~$S$,
the symmetric difference of any two of the five sets $X$, $\bar X$, $X^∘$, $\overline{X^∘}$, $(\bar X)^∘$ is a ^{rare set}.
The class~$N$ of ^{rare sets} is a ^{σ-ideal} (^!{σ-ideal}) of~$M$.

\proof Proof.
The interior of a closed subset~$F$ is a ^{clopen subset}~$F'$ and $F∖F'$ is a ^{rare set}.
The closure of an open subset~$U$ is a ^{clopen subset}~$U'$ and $U'∖U$ is a ^{rare set}.
Thus, $U=U'⊕(U'∖U)$ and $F=F'⊕(F∖F')$ are symmetric differences of a ^{clopen} set and a ^{rare} set, so the first three classes coincide.
\ppar
We have $U∪R=U⊕(R∖U)$ and $F∖R=F⊕(R∩F)$, so the last two classes consist of measurable sets.
Conversely, the relations $U⊕R=(U∖\bar R)∪((R∖U)∪((\bar R∖R)∩U))$ (the union is disjoint)
and $F⊕R=(F∪\bar R)∖((R∩F)∪(\bar R∖(R∪F)))$ (in the outer difference, the former set contains the latter)
show that any measurable set belongs to the fourth respectively fifth class.
Thus, all five classes coincide.
\ppar
By construction, ^{rare sets} are measurable
and by \vmeagerrare, the class of ^{rare sets} coincides with the class of ^{meager sets},
which is closed under passage to subsets and countable unions, and hence is a ^{σ-ideal}.
Since $S∖(U⊕R)=(S∖U)⊕R$ for any open set~$U$ and ^{rare set}~$R$, measurable sets are closed under complements.
Since $⋃_{i≥0}(U_i∪R_i)=⋃_{i≥0}U_i∪⋃_{i≥0}R_i$,
the set $⋃_{i≥0}U_i$ is open,
and the set $⋃_{i≥0}R_i$ is ^{meager}, and hence also ^{rare},
measurable sets are closed under countable unions.
Hence, the class of measurable sets is a ^{σ-algebra}.
\ppar
The presentation of a measurable set as the symmetric difference of a ^{clopen subset}~$U$ and a ^{rare} subset~$R$ is unique:
if $U⊕R=U'⊕R'$, then $U⊕U'=R⊕R'$ is ^{clopen} and ^{rare}.
Clopen rare sets are empty, so $U=U'$ and $R=R'$.
\ppar
The relations $\overline{F∖R}⊂F$ and $(U∪R)^∘⊃U$, combined with the fact that subsets of ^{rare} sets are again ^{rare},
show that the sets $X⊕\bar X$ and $X⊕X^∘$ are ^{rare},
which implies the claim about 5 sets obtained from~$X$ in the main statement.

\label\normalmeasurevaluationhyperstonean
\proclaim Lemma.
Given a ^{hyperstonean space}~$S$,
there is a canonical bijective correspondence between ^{normal measures} on~$S$
and ^{normal valuations} on the ^{hyperstonean locale} $"Ω(S)$,
given by restricting ^{normal measures} on~$S$ to open subsets of~$S$.

\proof Proof.
By definition, a ^{normal measure} on~$S$ restricts to a ^{normal valuation} on~$"Ω(S)$.
To show that a ^{normal valuation} on~$"Ω(S)$ extends to a unique ^{normal measure} on~$S$,
observe that by \vhyperstoneanmeasurable,
any measurable subset~$A$ of~$S$ has a unique presentation as the symmetric difference of a ^{clopen subset}~$U$ and a ^{rare subset}~$R$.
For any ^{normal measure}~$μ$ that extends~$ν$ we must have $μ(A)=μ(U⊕R)=μ(U)=ν(U)$, which shows uniqueness.
To show existence, we set $μ(A)=ν(U)$ and claim that $μ$ is a ^{normal measure}.
By construction, $μ$ extends~$ν$ and vanishes on all ^{meager subsets}.
To show additivity for a countable family of disjoint measurable subsets $A_i⊂S$ ($i≥0$),
present $A_i$ as the symmetric difference of a ^{clopen} subset $U_i⊂S$ and a ^{rare subset} $R_i⊂S$.
For any $i≠j$, the intersection $U_i∩U_j$ is ^{clopen} and ^{rare}, hence empty.
Thus, both $\{U_i\}_{i≥0}$ and $\{A_i\}_{i≥0}$ are disjoint families.
We have $$⋃_{i≥0}A_i⊂⋃_{i≥0}(U_i∪R_i)=⋃_{i≥0}U_i∪⋃_{i≥0}R_i,$$ where the subset $⋃_{i≥0}R_i$ is ^{meager}, hence also ^{rare}.
Hence, $$⋃_{i≥0}A_i⊕⋃_{i≥0}U_i$$ is a ^{rare subset} and $$μ\left(⋃_{i≥0}A_i\right)=μ\left(⋃_{i≥0}U_i\right).$$
Therefore, showing that $μ\left(⋃_{i≥0}A_i\right)=∑_{i≥0}μ(A_i)$
reduces to showing that $$μ\left(⋃_{i≥0}U_i\right)=∑_{i≥0}μ(U_i).$$
Since $⋃_i U_i=⋁_{k≥0} ⋁_{i<k} U_i$, invoking continuity and modularity (hence finite additivity) of~$ν$ completes the proof.

\label\normalmeasurevaluation
\proclaim Remark.
\vnormalmeasurevaluationhyperstonean\ is true for all ^{Stonean spaces}.
This result follows from Alvarez-Manilla [\ExtVal, Theorem~4.4],
which says that (finite in our case) continuous valuations on a regular topological space
extend uniquely to regular τ-smooth Borel measures.

We now show that our definition of ^{hyperstonean spaces}
is equivalent to Dixmier's definition [\HS, Definition~3].
This lemma is not used anywhere else and is only established to justify our terminology.

\proclaim Lemma.
The class of ^{hyperstonean spaces}
coincides with the class of ^{Stonean spaces} such that the union of supports of all positive ^{normal measures} is everywhere dense.

\proof Proof.
Suppose $S$ is a ^{Stonean space}
such that the union of supports of all positive ^{normal measures} is everywhere dense.
Then for any nonempty open subset $U⊂S$
there is a positive ^{normal measure}~$μ$ on~$S$
such that the support of~$μ$ is a ^{clopen} subset of~$S$
that intersects~$U$.
By definition of a ^{normal measure},
the measure~$μ$ restricts to a ^{normal valuation}~$ν:"Ω(S)→[0,∞)$ such that $ν(U)≠0$.
Since $U∈"Ω(S)∖\{∅\}$ is arbitrary, this proves that $"Ω(S)$ is a ^{hyperstonean locale} and therefore $S$ is a ^{hyperstonean space}.
\ppar
Suppose now that $S$ is a ^{hyperstonean space}.
By \vnormalmeasurevaluation, any ^{normal valuation} on~$"Ω(S)$ yields a ^{normal measure} on~$S$.
Since $"Ω(S)$ is hyperstonean, any open subset of~$S$ has a nonempty intersection with the support of some normal measure on~$S$.
This means that the union of all such supports is everywhere dense in~$S$,
and hence the union of supports of all positive ^{normal measures} on~$S$ is everywhere dense.

The following proposition follows immediately from our definition of ^{hyperstonean spaces}.

\label\HStonean
\proclaim Proposition.
The adjoint equivalence
$$@Stonean\ltogets5{"Ω}{"Sp}@StoneanLoc$$
restricts to an adjoint equivalence
$$@HStonean\ltogets5{"Ω}{"Sp}@HStoneanLoc.$$

\proclaim Proposition.
(^=:{Hyperstonean duality}.)
The adjoint equivalence
$$@StoneanLoc\ltogets9{"COpen}{"Ideal}@CBAlg^@op$$
restricts to an adjoint equivalence
$$@HStoneanLoc\ltogets9{"COpen}{"Ideal}@LBAlg^@op=@MLoc.$$

\proof Proof.
It suffices to show that the functor $"COpen$ sends
^{hyperstonean locales} to ^{localizable Boolean algebras}
and the functor $"Ideal$ sends ^{localizable Boolean algebras}
to ^{hyperstonean locales}.
Both claims follow immediately from \vhnormal.

\proclaim Remark.
Continuing \vcompactification, for a ^{hyperstonean locale}~$H$
we have a canonical dense inclusion of ^{locales}
$$"COpen(H)→H.$$
Equivalently, for a ^{measurable locale}~$L$ we have a canonical dense inclusion
$$L→"Ideal(L).$$
The inclusions exhibit their codomain, a ^{hyperstonean locale},
as the ^={hyperstonean compactification[|s]} of the domain, a ^{measurable locale}.
Furthermore, the domain is the smallest dense ^{sublocale} of the codomain.
The opens of the domain bijectively map to the ^{clopens} of the codomains.

\label\vonneumannalgebras
\section Duality between measurable locales and commutative von Neumann algebras

This section constructs in \vCVNAMLoc\ an adjoint equivalence of categories
$$@CVNA^@op\ltogets9{"ProjLoc}{"L^∞}@MLoc=@LBAlg^@op$$
of ^{measurable locales} and ^{commutative von Neumann algebras}
as the opposite adjoint equivalence of
$$@LBAlg\ltogets9{"Step}{"ProjAlg}@CVNA.$$
Here $"ProjLoc$ sends a commutative ^{von Neumann algebra} to its ^{localizable} ^{Boolean algebra of projections},
whereas $"L^∞$ sends a ^{localizable Boolean algebra} corresponding to a ^{measurable locale}~$L$
to the ^{commutative von Neumann algebra} of bounded ^{maps of locales} $L→`C$,
where $`C$ denotes the usual ^{locale} of complex numbers.
The essential surjectivity of $"ProjLoc$ is implied by the results of Segal [\DOA].

\label\reviewvna
\subsection Commutative von Neumann algebras

This section recalls the necessary definitions and facts from the theory of von Neumann algebras.
Although we try to give a reasonably detailed account,
our exposition is not self-contained,
and we refer the reader to Sakai's book [\CWstar] for more information.

All rings are unital and all homomorphisms of rings preserve units.
In particular, we require ^{C*-algebras} to be unital, since we do not need the nonunital ones.
We use $`C$ as our ring of coefficients throughout this paper.
(Real C*-algebras and real von Neumann algebras can be most easily defined and treated
through the complexification functor: $A$ is a real C*-algebra if $A⊗_`R `C$ is a complex C*-algebra,
and similarly for real von Neumann algebras.)

\proclaim Definition.
A ^={C*-algebra[|s]} is a complex unital algebra~$A$
equipped with a complex-antilinear map $*:A→A$ and a complete norm $‖{-}‖:A→[0,∞)$
such that $‖1‖≤1$, $‖xy‖≤‖x‖⋅‖y‖$,
$1^*=1$, $(xy)^*=y^*x^*$,
$‖x^*‖=‖x‖$, $‖x^*x‖=‖x^*‖⋅‖x‖$.
A ^={C*-homomorphism[|s]} $f:A→A'$ of ^{C*-algebras} is a unital homomorphism of complex algebras
such that $f(x^*)=f(x)^*$.
^^={*-homomorphism}
(Such homomorphisms are automatically contractive: $‖f(x)‖≤‖x‖$.)

\proclaim Definition.
A ^={von Neumann algebra[|s]} (alias ^={W*-algebra[|s]}) is a ^{C*-algebra}~$A$ that admits a ^={predual[|s]},
i.e., a complex Banach space~$A_*$ such that there is a (necessarily isometric) isomorphism $∫:A→(A_*)^*$ of complex Banach spaces.
A ^={morphism[|s] of von Neumann algebras} (alias ^={normal *-homomorphism[|s]} or ^={W*-homomorphism[|s]})
is a ^{C*-homomorphism} of ^{C*-algebras} $f:A→A'$ that admits a ^{predual},
which is given by preduals $(A_*,∫)$ and $(A'_*,∫')$ for $A$ and $A'$ respectively,
together with a (necessarily contractive) map of Banach spaces $f_*:A'_*→A_*$
such that the following square commutes:
\vadjust{\medskip}$$\arrowsize7 \cd{A&\mapright{f}&A'\cr
\lmapdown{∫}&&\mapdown{∫'}\cr
(A_*)^*&\mapright{(f_*)^*}&(A'_*)^*.\cr
}$$
The ^={category of von Neumann algebras} is denoted by $@=VNA$.
The ^={category of commutative von Neumann algebras} is denoted by $@=CVNA$.
^^={commutative von Neumann algebra[|s]}
^^={homomorphism[|s] of commutative von Neumann algebras}

\proclaim Remark.
If $(A,∫)$ and $(A',∫')$ are preduals of a ^{C*-algebra}~$A$,
then by Sakai's theorem [\CharW] there is a unique (necessarily isometric) isomorphism $g:A→A'$ such that $∫=g^*∘∫'$.
Thus all preduals of~$A$ are isomorphic via a unique isomorphism,
which allows us to talk about {\it the\/} predual of~$A$.
For a predual $(A,∫)$ of~$A$ the isometric isomorphism $∫:A→(A_*)^*$ endows $A$ with a topology
transferred from the weak topology on $(A_*)^*$.
This topology on~$A$ is independent of the choice of $(A,∫)$ and
is known as the ^={ultraweak topolog[y|ies]} ({\it σ-topology\/} in Sakai's book [\CWstar]).
Taking the dual of the isomorphism~$∫$ with these topologies,
we get an isomorphism $A_*→A^*$, i.e., $A_*$ is the ultraweak dual of~$A$.
Likewise, given a ^{morphism of von Neumann algebras} $f:A→A'$, we get a commutative square
\vadjust{\medskip}$$\arrowsize7 \cd{A^*&\mapleft{f^*}&A'^*\cr
\lmapup{∫^*}&&\mapup{∫'^*}\cr
A_*&\mapleft{f_*}&A'_*,\cr
}$$
where both vertical maps are isomorphism.
Thus, taking the dual in the ultraweak topology yields a functor,
the ^={predual functor[|s]}, which is denoted by $$(-)_*:@VNA^@op→@Banach.$$
Composing $(-)_*^"op$ with the dual Banach space functor
$$(-)^*:@Banach^@op→@Banach$$
yields a functor
$$((-)_*^"op)^*:@VNA→@Banach$$
and $∫$ is a natural isomorphism
$$∫:?id_{@VNA}→((-)_*^"op)^*.$$

\subsection Functors between commutative von Neumann algebras and measurable locales

\proclaim Definition.
For a von Neumann algebra~$A$ we define $$A^{≥0}=A^*A=\{a^*a\mid a∈A\}.$$
If $∫:A→(A_*)^*$ is the predual isomorphism for a ^{von Neumann algebra}~$A$,
then we define $$A_*^{≥0}=\{μ∈A_*\mid\Bigl(∫a\Bigr)(μ)≥0 \hbox{ for all $a∈A^{≥0}$}\}.$$

\label\ultraweaknormal
\proclaim Remark.
A positive linear functional $A→`C$ on a ^{von Neumann algebra}~$A$ is ultraweakly continuous
if and only if it is a ^={normal functional[|s]},
i.e., preserves suprema of bounded directed subsets.
Such functionals are in bijection with elements of $A_*^{≥0}$.
(Sakai [\CWstar, Theorem~1.13.2].)
Any element $μ∈A_*$ is a complex linear combination of at most four elements of $A_*^{≥0}$.
(Sakai [\CWstar, Theorem~1.14.3].)

\label\projiscomplete
\proclaim Proposition.
For any commutative ^{von Neumann algebra}~$A$,
its ^={Boolean algebra of projections} $$"=ProjAlg(A)=\{x∈A\mid x^*=x, x^2=x\}$$
equipped with the induced operations of addition $a⊕b=(a-b)^2$ and multiplication
is a ^{complete Boolean algebra}.

\proof Proof.
Idempotents in any commutative unital ring form a Boolean algebra with respect to the operations $a∧b=ab$ and $a∨b=a+b-ab$.
Self-adjoint idempotents are preserved under these operations, so they also form a Boolean algebra.
The ^{Boolean algebra} $"ProjAlg(A)$ is complete:
if $S⊂"ProjAlg(A)$ is an arbitrary directed subset of $"ProjAlg(A)$,
then for any $μ∈A_*^{≥0}$ the directed subset $μ_!(S)⊂[0,∞)$ is bounded by $μ(1)$, therefore $\lim_{s∈S} μ(s)$ exists
for all $μ∈A_*^{≥0}$.
Since any element $μ∈A_*$ is a complex linear combination of at most four elements of $A_*^{≥0}$ (\vultraweaknormal),
we deduce that $\lim_{s∈S}μ(s)$ exists for all $μ∈A_*$,
which by the definition of the ultraweak topology on~$A$ means that $\lim_{s∈S}s$ exists
and by continuity of the maps $a↦a^2-a$ and $a↦a-a^*$ must be a projection itself.

\label\predualvaluation
\proclaim Proposition.
For any ^{commutative von Neumann algebra}~$A$,
restriction along the inclusion of sets $"ProjAlg(A)→A$
establishes a bijection from $A_*$ to the vector space of complex ^{continuous valuations} on $"ProjAlg(A)$.

\proof Proof.
Given $μ∈A_*$, we show that its restriction~$ν$ to $"ProjAlg(A)$ is a complex ^{continuous valuation} on $"ProjAlg(A)$.
As established in \vprojiscomplete, the ^{Boolean algebra} $"ProjAlg(A)$ is complete.
The properties $ν(0)=0$ and $ν(x)+ν(y)=ν(x∨y)+ν(x∧y)$ follow from the linearity of~$μ$.
To show that $ν(\sup_{x∈I} x)=\lim_{x∈I} ν(x)$ for any directed subset $I⊂L$,
it suffices to observe that $\sup_{x∈I}x=\lim_{x∈I}x$, where the limit is taken in the ^{ultraweak topology}
(Sakai [\CWstar, Lemma 1.7.4]).
\ppar
Given a complex ^{continuous valuation} $ν$ on $"ProjAlg(A)$,
we show that it admits a unique extension to an ultraweakly continuous map $A→`C$.
By the ^{Hahn–Jordan decomposition of normal valuations} (^!{Hahn–Jordan decomposition of normal valuations}),
it suffices to treat the positive case.
First, since $ν$ is a ^{valuation}, it admits a unique extension to a linear functional on the linear span of $"ProjAlg(A)$ inside~$A$.
Secondly, since $ν$ is a ^{continuous valuation}, the associated linear functional~$f$ is a ^{normal functional} (^!{normal functional}),
hence also ultraweakly continuous by \vultraweaknormal.
Finally, the linear span of all projections is ultraweakly dense in~$A$,
so $f$ extends uniquely to an ultraweakly continuous functional $A→`C$.

\proclaim Definition.
The ^={projection locale functor} for ^{commutative von Neumann algebras} is a functor $$"=ProjLoc:@CVNA^@op→@MLoc=@LBAlg^@op$$ defined as the opposite
functor of $$"ProjAlg:@CVNA→@LBAlg$$ that
sends a commutative ^{von Neumann algebra}~$A$ to its ^{Boolean algebra of projections} (defined in ^!{Boolean algebra of projections}), which is ^{localizable}.
A ^{morphism[| of von Neumann algebras]} $f:A→A'$ is sent to its restriction $"ProjAlg(f):"ProjAlg(A)→"ProjAlg(A')$ to projections in its domain and codomain.

\proclaim Lemma.
This definition is correct.

\proof Proof.
The ^{Boolean algebra} $"ProjAlg(A)$ is ^{complete} by \vprojiscomplete.
It is also ^{localizable}:
given $m∈"ProjAlg(A)$ such that $m≠0$, its image under the isomorphism $∫:A→(A_*)^*$ is a nonzero element $κ∈(A_*)^*$,
meaning there is $ν∈A_*$ such that $κ(ν)=ν(m)≠0$.
By decomposing $ν$ as a complex linear combination of at most four elements of $A_*^{≥0}$,
we can assume $ν≥0$ and $ν(m)≠0$.
By \vpredualvaluation, $ν$ is a ^{continuous valuation} on $"ProjAlg(A)$, so $A$ is ^{localizable} by \vlmcrit.
\ppar
For any ^{normal *-homomorphism} $f:A→A'$
the map $"ProjAlg(f):"ProjAlg(A)→"ProjAlg(A')$
is a ^{homomorphism of Boolean algebras} because $f(x^*)=f(x)^*$ and $f(x^2)=f(x)^2$.
If $S⊂"ProjAlg(A)$ is a directed subset,
then for any $μ∈(A'_*)^{≥0}$
we have $μ(\sup f(S))=\sup μ(f(S))=μ(f(\sup S))$
because $μ$ and $μ∘f$ preserve suprema
by virtue of being elements of $(A'_*)^{≥0}$ respectively $A_*^{≥0}$.
Since $μ∈(A'_*)^{≥0}$ is arbitrary, this shows that $\sup f(S)=f(\sup S)$, so $f$ is a ^{complete homomorphism}.

\proclaim Remark.
The composition $$@CVNA^@op\lto9{"ProjLoc}@MLoc\lto9{"Ideal}@HStoneanLoc\lto9{"Sp}@HStonean$$ is the ^={Gelfand spectr[um|a]} functor
restricted to the ^{category of commutative von Neumann algebras},
as observed by de Groote in [\ObsSto, Theorem~3.2] and [\CQO, Theorem~5.2.1].

\proclaim Definition.
The functor $$"L^∞:@MLoc→@CVNA^@op$$ sends a ^{measurable locale}~$L$ to the ^{commutative von Neumann algebra}
of bounded (meaning factoring through some open ball) ^{maps of locales} $L→`C$,
where $`C$ denotes the usual ^{locale} of complex numbers.
The functor $"L^∞$ sends a ^{map of measurable locales} $f:L→L'$
to the ^{morphism of von Neumann algebras} $"L^∞(L')→"L^∞(L)$
that sends $h∈"L^∞(L')$ ($h:L'→`C$) to $h∘f∈"L^∞(L)$.
The opposite functor is denoted by $$"=Step=("L^∞)^"op:@LBAlg→@CVNA.$$

\label\predualbounded
\proclaim Lemma.
This definition is correct.

\proof Proof.
$"L^∞(L)$ is a commutative complex ^{C*-algebra} because $`C$ is.
Completeness follows from the existence of a predual shown below, so does not require a separate proof,
and the other properties are of purely algebraic nature.
\ppar
Below, we construct the predual in the real setting and then complexify by tensoring with~$`C$ over~$`R$.
We set the real predual $"L^∞(L)_*$ to the normed vector space $"L^1(L)$ of real ^{continuous valuations} on~$L$,
equipped with the ^{variation norm} defined in ^!{variation norm}.
\ppar
Consider the linear span of characteristic maps $χ_m:L→`C$ for all $m∈L$.
We refer to its elements as {\it step maps}.
Real step maps are order-dense in the real part of $"L^∞(L)$
because any real $f∈"L^∞(L)$ is the supremum of step functions given by finite subsums of $∑_{k∈`Z}k/n⋅χ_{f^*[k/n,(k+1)/n)}$ as $n→∞$.
Thus, $"L^∞(L)$ is the Dedekind completion of the Archimedean vector lattice of real step maps.
This justifies the notation $"Step$ for the functor $("L^∞)^"op$.
\ppar
Consider the linear span of characteristic functionals $ψ_m:"L^1(L)→`C$ ($μ↦μ(m)$) for all $m∈L$.
We refer to its elements as {\it step functionals}.
Real step functionals are order-dense in the real part of $("L^1(L))^*$ because
any real $F∈("L^1(L))^*$ is the supremum of step functionals given by finite subsums of $∑_{k∈`Z}k/n⋅ψ_{G_{(k+1)/n}-G_{k/n}}$,
where $G_α∈L$ for $α∈`R$ is defined as the supremum of all $p∈L$ such that $F(μ(p∧-))≤αμ(p)$ for all $μ∈"L^1(L)$.
%Since $F$ is a bounded functional on~$"L^1(L)$ (meaning $|F(μ)|≤‖F‖⋅μ(1)$ for any $μ≥0$ in $"L^1(L)$), we must have $G_α=0$ for all $α≤-‖F‖$ and $G_α=1$ for all $α>‖F‖$.
Indeed, if we fix some $μ∈"L^1(L)$, then
$$F(μ)-∑_{k∈`Z}k/n⋅ψ_{G_{(k+1)/n}-G_{k/n}}(μ)=∑_{k∈`Z}F(μ((G_{(k+1)/n}-G_{k/n})∧-))-∑_{k∈`Z}k/n⋅μ(G_{(k+1)/n}-G_{k/n})=A,$$
where $$k/n⋅μ(G_{(k+1)/n}-G_{k/n})≤F(μ((G_{(k+1)/n}-G_{k/n})∧-))≤(k+1)/n⋅μ(G_{(k+1)/n}-G_{k/n}),$$
so $A≥0$ and $$A≤∑_{k∈`Z}1/n⋅μ(G_{(k+1)/n}-G_{k/n})=μ(1)/n.$$
Thus, $A→0$ as $n→∞$, establishing the claim.
Hence, $("L^1(L))^*$ is the Dedekind completion of the Archimedean vector lattice of real step functionals.
\ppar
Consider the map that sends $χ_m$ to $ψ_m$ for any $m∈L$.
This map extends by linearity to an isomorphism of Archimedean vector lattices from step functions to step functionals.
Indeed, the resulting map~$ζ$ is surjective by construction.
For injectivity, observe first that $ψ_m≠0$ whenever $m≠0$.
Furthermore, an arbitrary linear combination of~$ψ_m$ can be analyzed
by restricting separately to each possible meet of elements of the form $m$ and $¬m$ participating in the linear combination,
after which both linear combinations take the form $αχ_m$ respectively $αψ_m$ for some $m∈L$.
Thus, the map~$ζ$ is an isomorphism and induces an isomorphism $$∫_L:"L^∞(L)→("L^1(L))^*$$ between the Dedekind completions of these Archimedean vector lattices,
which shows that $"L^1(L)$ is the predual of $"L^∞(L)$.
This is the localic analogue of the Lebesgue integral, traditionally denoted by $∫_L f ~d μ$ instead of $∫_L(f)(μ)$.
\ppar
Finally, given a ^{morphism of measurable locales} $f:L→L'$,
we construct a predual of the induced ^{C*-homomorphism} of ^{C*-algebras} $"L^∞(L')→"L^∞(L)$,
which has to be a map of the form $"L^1(f):"L^1(L)→"L^1(L')$.
This is precisely the pushforward map on ^{continuous valuations}:
given $μ∈"L^1(L)$, we set $"L^1(f)(μ)(m')=μ(f^*(m'))$.
This is a contractive map of Banach spaces,
and to verify the predual property, it suffices to verify the commutativity of the diagram
\vadjust{\medskip}$$\arrowsize7 \cd{"L^∞(L')&\mapright{"L^∞(f)}&"L^∞(L)\cr
\mapdown{∫_{L'}}&&\mapdown{∫_L}\cr
"L_1(L')^*&\mapright{"L_1(f)^*}&"L_1(L)^*\cr
}$$
on all characteristic maps~$χ_{m'}$, where $m'∈L'$.
Indeed, $$"L_1(f)^*\Bigl(∫_{L'}(χ_{m'})\Bigr)="L_1(f)^*(μ'↦μ'(m'))=(μ↦μ(f^*(m'))),$$
whereas $$∫_L"L^∞(f)(χ_{m'})=∫_Lχ_{f^*(m')}=(μ↦μ(f^*(m'))),$$
so both values coincide.

\subsection Equivalence of commutative von Neumann algebras and measurable locales

In this section, we construct the counit and unit of the adjunction
$$@CVNA^@op\ltogets9{"ProjLoc}{"L^∞}@MLoc=@LBAlg^@op$$
by constructing the unit~$χ$ and counit~$μ$ of the opposite adjunction
$$@LBAlg\ltogets9{"Step}{"ProjAlg}@CVNA$$
and prove that they are isomorphisms and satisfy the triangle identities.

\proclaim Definition.
The natural isomorphism $$χ:?id_@LBAlg→"ProjAlg∘"Step$$
sends a ^{localizable Boolean algebra}~$L$ to the ^{map of localizable Boolean algebras} $χ_L:L→"ProjAlg("Step(L))$
that sends $m∈L$ to the projection $χ_m$ in the ^{von Neumann algebra} $"Step(L)$ given by the characteristic map~$χ_m$ of~$m$.

\proclaim Lemma.
This definition is correct.

\proof Proof.
We have $χ_{mm'}=χ_m χ_{m'}$, $χ_1 = 1$, $χ_{m∨m'}=χ_m ∨ χ_{m'}$, and $χ_0 = 0$.
Thus, $χ_L$ is a homomorphism of lattices that preserves 0 and 1,
hence also a homomorphism of Boolean algebras.
The map $χ_L$ is injective because different $m∈L$ yield different characteristic maps.
It is surjective because any ^{map of locales} $f:L→`C$ such that $f^2=f$ and $f^*=f$
necessarily factors through the ^{map of locales} given by the inclusion of the discrete locale on~$\{0,1\}$ into~$`C$.
Since Boolean algebras form a variety of algebras, bijective homomorphisms of ^{Boolean algebras}
are isomorphisms of ^{Boolean algebras},
hence also complete isomorphisms of ^{complete Boolean algebras}.
\ppar
The naturality condition is satisfied because the following square commutes for any ^{morphism of localizable Boolean algebras} $h:L'→L$:
$$\arrowsize17 \sqcd{
"ProjAlg("Step(L))&\mapleft{"ProjAlg("Step(h))}&"ProjAlg("Step(L'))\cr
\mapup{χ_L}&&\mapup{χ_{L'}}\cr
L&\mapleft{h}&L'.\cr
}$$
Indeed, evaluating on an arbitrary $m'∈L'$ yields
$$χ_L(h(m'))=χ_{h(m')}$$
and
$$"ProjAlg("Step(h))(χ_{L'}(m'))="ProjAlg("Step(h))(χ_{m'})="Step(h)(χ_{m'})=χ_{h(m')}.$$
This establishes the naturality property and completes the proof.

\proclaim Definition.
The natural isomorphism $$μ:"Step∘"ProjAlg→?id_@CVNA$$
sends a commutative ^{von Neumann algebra}~$A$
to the ^{morphism of von Neumann algebras} $$"Step("ProjAlg(A))→A$$
whose value on the characteristic map~$χ_m$ of any projection $m∈A$ equals~$m$.
(The inverse isomorphism essentially encodes the spectral theorem for bounded operators:
it sends an element $a∈A$ to the ^{map of locales} $"ProjAlg(A)→`C$
whose inverse image map sends an open subset of~$`C$ to the corresponding spectral projection of~$a$.)

\proclaim Lemma.
This definition is correct.

\proof Proof.
We work with the real parts of $"Step("ProjAlg(A))$ and~$A$.
The canonical inclusion $"ProjAlg(A)→A$ induces a ^{*-homomorphism} to~$A$
from the linear span of all characteristic maps $χ_m:"ProjLoc(A)→`R$, where $m∈"ProjLoc(A)$ is a projection in~$A$.
This map preserves the natural partial order on both sides.
Recall from the proof of \vpredualbounded\ that this linear span is order-dense in the real part of $"Step("ProjAlg(A))$.
Therefore, we have a canonical order-preserving map $μ_A:"Step("ProjAlg(A))→A$, which is a ^{C*-homomorphism} by continuity.
\ppar
To show that $μ_A$ is a ^{morphism of von Neumann algebras}, we construct its ^{predual}
as a morphism of Banach spaces
$$A_*→"L^1("ProjAlg(A)),$$
where $"L^1$ denotes the vector space of real ^{continuous valuations} of $"ProjAlg(A)$ as in the proof of \vpredualbounded.
Indeed, \vpredualvaluation\ constructs such a map and proves that it is an isomorphism.
Therefore, its dual map, given by the bottom map in the following square diagram, is also an isomorphism.
The ^{predual} square
$$\cd{"Step("ProjAlg(A))&\mapright{μ_A}&A\cr
\lmapdown{∫_{"ProjAlg(A)}}&&\mapdown{∫}\cr
("L^1("ProjAlg(A)))^*&\mapright{}&(A_*)^*\cr
}$$
commutes because the two compositions coincide on step maps $χ_m$ ($m∈"ProjAlg(A)$)
by construction, therefore they coincide on their linear span and its order-closure.
Thus, $μ_A$ is an isomorphism of ^{von Neumann algebras} because the other three maps in the square are isomorphisms.
\ppar
The naturality condition is satisfied because the following square commutes for any homomorphism $h:A→A'$ of commutative ^{von Neumann algebras}:
$$\arrowsize17 \sqcd{
"Step("ProjAlg(A))&\mapright{"Step("ProjAlg(h))}&"Step("ProjAlg(A'))\cr
\mapdown{μ_A}&&\mapdown{μ_{A'}}\cr
A&\mapright{h}&A'.\cr
}$$
It suffices to verify commutativity when evaluated on projections in $"Step("ProjAlg(A))$.
Indeed, given such a projection $χ_m$ for some $m∈"ProjAlg(A)$, we have
$$h(μ_A(χ_m))=h(m)$$
and
$$μ_{A'}("Step("ProjAlg(h))(χ_m))=μ_{A'}(χ_m∘"ProjAlg(h))=μ_{A'}(χ_{h(m)})=h(m),$$
which completes the proof.

Assembling all the results of this section, we obtain the following result.
While several parts of this result were known before
(e.g., the essential surjectivity of $"ProjLoc$ falls out of results due to Segal [\DOA]),
there does not appear to be a published source that assembles the individual ingredients in a single proof.

\label\CVNAMLoc
\proclaim Theorem.
The functors
$$"Step:@LBAlg→@CVNA$$
and
$$"ProjAlg:@CVNA→@LBAlg$$
together with natural isomorphisms
$$χ:?id_@LBAlg→"ProjAlg∘"Step$$
and
$$μ:"Step∘"ProjAlg→?id_@CVNA$$
form an adjoint equivalence of categories.
Therefore, the opposite adjunction
$$@CVNA^@op\ltogets9{"ProjLoc}{"L^∞}@MLoc=@LBAlg^@op$$
is also an adjoint equivalence of categories.

\proof Proof.
It remains to show that the exhibited equivalence is an adjoint equivalence.
The two triangle identities imply one another, so it suffices to establish just one of them.
To show that the composition
$$"Step(L)\lto9{"Step(χ_L)}"Step("ProjAlg("Step(L)))\lto9{μ_{"Step(L)}}"Step(L)$$
equals identity,
it suffices to evaluate it on all characteristic maps~$χ_m$ for $m∈L$.
We have
$$"Step(χ_L)(χ_m)=χ_{χ_m},$$
i.e., the characteristic map $χ_{χ_m}:"ProjLoc("L^∞(L))→`R$ of the projection $χ_m∈"ProjLoc("L^∞(L))$.
Then
$$μ_{"Step(L)}(χ_{χ_m})=χ_m,$$
which completes the proof.

\label\pointset
\section Point-set measurable spaces

This section defines the category $@CSLEMS$ in \vmainthm,
namely, the category of ^{compact strictly localizable enhanced measurable spaces},
and explores its properties that will be necessary later.

\subsection Enhanced measurable spaces

\proclaim Definition.
A ^={σ-algebra[|s]} on a set~$X$ is a collection of subsets of~$X$ closed under complements and countable unions.

\proclaim Definition.
A ^={σ-ideal}~$N$ of a ^{σ-algebra}~$M$ on a set~$X$
is a subset of~$M$ closed under countable unions and passage to subsets in~$M$ (meaning $A⊂B$, $A∈M$, and $B∈N$ imply $A∈N$).

In particular, a ^{σ-ideal} of the ^{σ-algebra}~$2^X$ consisting of all subsets of~$X$
is a collection of subsets of~$X$ closed under passage to subsets and countable unions.

\proclaim Definition.
A (complete) ^={enhanced measurable space[|s]} is a triple $(X,M,N)$, where $X$ is a set, $M$ is a ^{σ-algebra} on~$X$,
and $N$ is a ^{σ-ideal} of~$2^X$ such that $N⊂M$.
A ^={measurable set[|s]} is an element of~$M$ and a ^={negligible set[|s]} is an element of~$N$.
A ^={conegligible set[|s]} is the complement of a negligible set with respect to~$X$.
^^={measurable subset[|s]}
^^={negligible subset[|s]}
^^={conegligible subset[|s]}
^^={conegligible}

\proclaim Definition.
If $(X,M,N)$ is an ^{enhanced measurable space}
and $m∈M$, then the ^={induced enhanced measurable space[|s]}
on~$m$ is the triple $(m,M_m,N_m)$,
where $M_m=\{m'∈M\mid m'⊂m\}$ and $N_m=\{n∈N\mid n⊂m\}$.

\proclaim Definition.
The category~$@=PreEMS$ is defined as follows.
Its objects are ^{enhanced measurable spaces}.
Morphisms $(X,M_X,N_X)→(Y,M_Y,N_Y)$ are ^={premap[s|] of enhanced measurable spaces},
^^={premap[|s]}
^^={premorphism[|s] of enhanced measurable spaces}
defined as maps of sets $f:X'→Y$ such that $X'⊂X$ is a ^{conegligible set} (denoted by $?=pdom f$, which stands for ^={point-set domain[|s]}),
for any $m∈M_Y$ we have $f^*m∈M_X$,
and for any $n∈N_Y$ we have $f^*n∈N_X$.

Wendt [\MIFHS, \CDis, \MHS, \CBase] refers to the last property (i.e., $(f^*)_!N_Y⊂N_X$) as “measure zero reflecting”.

\proclaim Proposition.
This definition is correct.

\proof Proof.
We define $?pdom(g∘f)=?pdom f∩f^*(?pdom g)$,
which is conegligible because both factors are,
the latter because $f^*$ preserves ^{conegligible subsets}.
Morphisms are composed like maps of sets, after we restrict them to the new ^{point-set domains}.
We have $$\eqalign{?pdom(h∘(g∘f))&=?pdom(g∘f)∩(g∘f)^*(?pdom h)\cr&=?pdom f∩f^*(?pdom g)∩f^*g^*(?pdom h)\cr}$$
and $$\eqalign{?pdom((h∘g)∘f)&=?pdom f∩f^*(?pdom(h∘g))\cr&=?pdom f∩f^*(?pdom g∩g^*(?pdom h))\cr&=?pdom f∩f^*(?pdom g)∩f^*g^*(?pdom h),\cr}$$
so composition is indeed associative.
The point-set composition is a ^{premap} because $(g∘f)^*=f^*∘g^*$.

\proclaim Remark.
One could drop the ^={completeness} condition and define a category $@=PreIEMS$ (I for incomplete)
of not necessarily complete ^{enhanced measurable spaces}.
Its objects are triples $(X,M,N)$,
where $X$ and $M$ are as above and $N$ is a ^{σ-ideal} of~$M$.
A morphism $(X,M_X,N_X)→(Y,M_Y,N_Y)$
is a map of sets $f:X'→Y$ such that $X'$ is ^{conegligible} (now meaning that $X∖X'$ is a subset of some $n∈N_X$),
for any $m_Y∈M_Y$ we have $f^*(m_Y)⊕m_X⊂n_X$ for some $m_X∈M_X$ and $n_X∈N_X$,
and for any $n_Y∈N_Y$ we have $f^*(n_Y)⊂n_X$ for some $n_X∈N_X$.
The inclusion functor $ι:@PreEMS→@PreIEMS$ exhibits $@PreEMS$ as a full subcategory
of $@PreIEMS$, with its image consisting precisely of {\it complete\/} objects,
i.e., objects $(X,M,N)$ for which
$N$ is a ^{σ-ideal} of~$2^X$, meaning $N$ is closed under passage to subsets in~$X$.
The completion functor $"C:@PreIEMS→@PreEMS$
sends a space $(X,M_X,N_X)$ to $(X,M'_X,N'_X)$,
where $n'∈N'_X$ if there is $n∈N_X$ such that $n'⊂n$
and $m'∈M'_X$ if there is $m∈M_X$ such that $m'⊕m∈N'_X$.
On morphisms, $"C$ retains the underlying map of sets.
The functor~$"C$ is fully faithful.
For any $\hat X=(X,M_X,N_X)$ the identity map of sets $X→X$
yields isomorphisms $\hat X→"C\hat X→\hat X$ in $@PreIEMS$.
Thus, the inclusion $ι:@PreEMS→@PreIEMS$ together with the completion functor $"C:@PreIEMS→@PreEMS$ is an adjoint equivalence of categories.
Hence, there is no benefit for us to consider noncomplete spaces.

\proclaim Remark.
We could require $?pdom f=X$ (i.e., $X'=X$), which would work in most of the paper.
However, a crucial point where it is absolutely necessary to allow $?pdom f≠X$
is \vlifting, where we construct a morphism from an ^{enhanced measurable space} $(X,M,N)$
to the ^{Stone spectrum} of the ^{localizable Boolean algebra} $M/N$.
If $X∈N$, then $M/N=0$ and its ^{Stone spectrum} is empty, so the only map of sets into it must have $?pdom f=∅$.
Otherwise we would have to manually exclude the case $X∈N$, $X≠∅$ from ^{enhanced measurable spaces},
which would create even more problems.

\subsection Equality almost everywhere

\proclaim Definition.
The category $@=StrictEMS$ of ^{enhanced measurable spaces} and ^={strict map[s|]}
^^={strict morphism[s|] of enhanced measurable spaces}
^^={strict map[|s] of enhanced measurable spaces}
is the quotient of the category $@PreEMS$
by the equivalence relation~$∼$ of ^={equality almost everywhere},
^^={equal almost everywhere}
for which
$$f∼f':(X,M_X,N_X)→(Y,M_Y,N_Y)$$
if $\{x∈X\mid f(x)≠f'(x)\}∈N_X$.
(If one of $f(x)$ or $f'(x)$ is undefined,
then the other one must be defined in order for $f(x)≠f'(x)$ to hold.
The easiest way to think about this is to assume $f(x)=*$ for some $*∉Y$ whenever $f(x)$ is undefined.)

\proclaim Proposition.
This definition is correct.

\proof Proof.
Throughout this proof we follow the convention that $f(x)=*∉Y$ whenever $f(x)$ is undefined.
Observe that $f∼f'$ is indeed an equivalence relation: the reflexivity and symmetry properties are obvious,
and if $f∼f'$ and $f'∼f''$, then $f∼f''$ because $$\{x∈X\mid f(x)≠f''(x)\}⊂\{x∈X\mid f(x)≠f'(x)\}∪\{x∈X\mid f'(x)≠f''(x)\}∈N_X.$$
The equivalence relation is compatible with composition:
if $$f∼f':(X,M_X,N_X)→(Y,M_Y,N_Y)$$ and $$g:(Y,M_Y,N_Y)→(Z,M_Z,N_Z),$$ then $g∘f∼g∘f'$ because $$\{x∈X\mid g(f(x))≠g(f'(x))\}⊂\{x∈X\mid f(x)≠f(x')\}∈N_X.$$
Likewise, if $f:(X,M_X,N_X)→(Y,M_Y,N_Y)$ and $g∼g':(Y,M_Y,N_Y)→(Z,M_Z,N_Z)$,
then $g∘f∼g'∘f$ because $$\{x∈X\mid g(f(x))≠g'(f(x))\}⊂f^*\{y∈Y\mid g(y)≠g'(y)\}∈N.$$
Thus, the quotient category exists.

\proclaim Definition.
The category $@=EMS$ of ^{enhanced measurable spaces}
^^={morphism[|s] of enhanced measurable spaces}
^^={map[|s] of enhanced measurable spaces}
is the quotient of the category $@PreEMS$
by the equivalence relation~$≈$ of ^={weak equality almost everywhere},
^^={weakly equal almost everywhere}
for which $$f≈f':(X,M_X,N_X)→(Y,M_Y,N_Y)$$ if for any $m∈M_Y$
we have $χ_m∘f∼χ_m∘f'$, where $$χ_m:(Y,M_Y,N_Y)→(\{0,1\},2^{\{0,1\}},\{∅\})$$
is the ^={characteristic map} of~$m$:
we have $χ_m(y)=1$ if and only if $y∈m$, for any $y∈Y$.
Equivalently, $f≈f'$ if for any $m∈M_Y$ the symmetric difference $f^*m⊕f'^*m$ belongs to~$N_X$.

\proclaim Proposition.
This definition is correct.

\proof Proof.
Observe that $f≈f'$ is an equivalence relation: the reflexivity and symmetry properties are obvious,
and if $f≈f'$ and $f'≈f''$, then $f≈f''$ because $χ_m∘f∼χ_m∘f'∼χ_m∘f''$ for all $m∈M_Y$,
and we already know that $∼$ is an equivalence relation.
The equivalence relation is compatible with composition:
if $f≈f':(X,M_X,N_X)→(Y,M_Y,N_Y)$ and $g:(Y,M_Y,N_Y)→(Z,M_Z,N_Z)$,
then $g∘f∼g∘f'$ because for any $m∈M_Z$ we have $$(g∘f)^*m⊕(g∘f')^*m=f^*(g^*m)⊕f'^*(g^*m)∈N_X$$
since $g^*m∈M_Y$.
Likewise, if $f:(X,M_X,N_X)→(Y,M_Y,N_Y)$ and $g≈g':(Y,M_Y,N_Y)→(Z,M_Z,N_Z)$,
then $g∘f≈g'∘f$ because for any $m∈M_Z$ we have $$(g∘f)^*m⊕(g'∘f)^*m=f^*(g^*m⊕g'^*m)∈N_X$$
since $g^*m⊕g'^*m∈N_Y$.
Thus, the quotient category exists.

\proclaim Remark.
If $Y≠∅$, any ^{premap} $f:(X,M_X,N_X)→(Y,M_Y,N_Y)$ with $?pdom f≠X$
can be extended to a ^{premap} $f':(X,M_X,N_X)→(Y,M_Y,N_Y)$
by setting $f(x)=y$ for all $x∈X∖?pdom f$,
for some fixed $y∈Y$.
We have $f∼f'$, so if $Y≠∅$, we could require ^{premaps} to be everywhere defined.
However, if $Y=∅$ and $X∈N_X$, there are no ^{premaps} $f:(X,M_X,N_X)→(Y,M_Y,N_Y)$ with $?pdom f=X$,
even though the two ^{enhanced measurable spaces} are isomorphic in the category $@StrictEMS$
via the empty map of sets $∅→∅$, since $∅$ is conegligible in both spaces.

\proclaim Definition.
We define two ^{enhanced measurable spaces}:
$`C_~=Borel=(`C,M_~Borel,\{∅\})$ and $`C_~=Lebesgue=(`C,M_~Lebesgue,N_~Lebesgue)$.
Here $M_~Borel$ is the set of all Borel subsets of~$`C$,
$M_~Lebesgue$ is the set of all Lebesgue-measurable subsets of~$`C$,
and $N_~Lebesgue$ is the set of all Lebesgue-negligible subsets of~$`C$.
In particular, $M_~Lebesgue$ is the set of unions of elements of $M_~Borel$ and $N_~Lebesgue$.
We define $`R_~Borel$ and $`R_~Lebesgue$ in a similar way.

\proclaim Remark.
$`C_~Borel$ and $`C_~Lebesgue$ are not isomorphic in $@EMS$.
As we will see later, the Boolean algebra $M/N$ is an isomorphism invariant of an object in $@EMS$.
For $`C_~Borel$ this Boolean algebra has atoms given by singleton subsets of~$`C$,
whereas for $`C_~Lebesgue$ this Boolean algebra is atomless.

\proclaim Example.
Recall that any Lebesgue-measurable function $`C→`C$ is ^{equal almost everywhere} to a Borel measurable function,
so after passing to the category $@EMS$ there is no difference between the two notions,
both give a morphism $`C_~Lebesgue→`C_~Borel$.
Morphisms $`C_~Lebesgue→`C_~Lebesgue$ form a proper subset of Lebesgue measurable functions.
(Preimages of negligible subsets of~$`C$ under Lebesgue measurable functions need not be negligible, as witnessed by constant functions.)
Morphisms $`C_~Borel→`C_~Borel$ are precisely Borel measurable functions, without any identification.
There are no morphisms $`C_~Borel→`C_~Lebesgue$ because preimages of ^{negligible subsets} must be negligible, hence empty,
but any singleton subset of the codomain is negligible.

\subsection Comparison of equivalence relations on morphisms

We now compare the two equivalence relations $∼$ and~$≈$.

\proclaim Lemma.
If $f∼f':(X,M_X,N_X)→(Y,M_Y,N_Y)$, then $f≈f'$.
In particular, $@EMS$ is a quotient of $@StrictEMS$.

\proof Proof.
We have $$\{x∈X\mid χ_m(f(x))≠χ_m(f'(x))\}⊂\{x∈X\mid f(x)≠f'(x)\}$$ for any $m∈M_Y$.

\proclaim Definition.
An ^{enhanced measurable space} $(X,M,N)$
is ^={countably separated}
^^={countable separability}
^^={countably separated enhanced measurable space[s|]}
if there is a countable subset $M_s⊂M$ such that
for any distinct $x,x'∈X$ there is $m∈M_s$ with $x∈m$ and $x'∉m$
and also for any $x∈X$ there is $m∈M_s$ with $x∈m$.
(The last condition is nontrivial only when $X$ is a singleton and is necessary below for ^{premaps} whose ^{point-set domain} is a proper subset.)

\label\cseplemma
\proclaim Lemma.
If $f≈f':(X,M_X,N_X)→(Y,M_Y,N_Y)$ and $(Y,M_Y,N_Y)$ is ^{countably separated}, then $f∼f'$.
In particular, the quotient map $@StrictEMS→@EMS$ becomes an equivalence of categories
if we restrict to countably separated spaces on both sides.

\proof Proof.
We have $$\{x∈X\mid f(x)≠g(x)\}⊂⋃_{m∈M_s}(f^*m⊕g^*m)∈N_X.\qed$$

\label\weakeqnecessary
\proclaim Remark.
\vcseplemma\ is false without the assumption of ^{countable separability},
see Fremlin [\MTiii, Example 343I], which constructs a morphism $f:X→X$ such that $f≈?id_X$ and $?pdom f=X$, but $f(x)≠x$ for all $x∈X$,
where $X=\{0,1\}^`R$ is the ^{σ-finite enhanced measurable space} given by the product of a continuum many copies of the two-point space $\{0,1\}$.
In particular, the quotient functor $@StrictEMS→@EMS$ is not an equivalence of categories.
The condition of countable separability considered above is very restrictive:
by Lemma~343E in Fremlin [\MTiii] it is equivalent to the existence of a ^{morphism of enhanced measurable spaces} $f:(X,M,N)→`R_~Borel$
whose underlying map of sets is injective.
Thus, for arbitrary ^{enhanced measurable spaces} (even if assumed to be ^{σ-finite})
the equivalence relation of ^{equality almost everywhere} is too strict:
it fails to identify distinct ^{premaps} of ^{enhanced measurable spaces}
that are sent to identical maps of the corresponding ^{measurable locales}
or ^{commutative von Neumann algebras} by the relevant functors,
since these functors factor through $@EMS$.

\subsection Coproducts of enhanced measurable spaces

\proclaim Lemma.
The category $@PreEMS$ admits small coproducts
and the forgetful functor $@PreEMS→@Set$ creates (i.e., preserves and reflects) small coproducts.

\proof Proof.
We construct the coproduct $∐_{i∈I}(X_i,M_i,N_i)$
by setting the underlying set to $Y=∐_{i∈I}X_i$
and declaring that $m∈M_Y$ (respectively $n∈N_Y$)
if for every $i∈I$ we have $m∩X_i∈M_i$ (respectively $n∩X_i∈N_i$).
The injection maps $ι_i:X_i→∐_{i∈I}X_i$ yield morphisms in $@PreEMS$.
To show the universal property of coproducts, it suffices to observe
that given a collection of morphisms $(f_i:(X_i,M_i,N_i)→(Z,M_Z,N_Z))_{i∈I}$ (with $?pdom f_i=X'_i⊂X_i$) in $@PreEMS$,
the induced map of sets $[f_i]_{i∈I}:∐_{i∈I}X'_i→Z$ is a morphism in $@PreEMS$, which follows from the definition of $M_Y$ and $N_Y$.
Since the forgetful functor reflects isomorphisms, reflection of coproducts is implied by preservation of coproducts.

\proclaim Lemma.
The category $@EMS$ admits small coproducts
and the quotient functor $@PreEMS→@EMS$ preserves small coproducts.

\proof Proof.
Since $@PreEMS$ admits small coproducts,
it suffices to show that given a coproduct cocone in $@PreEMS$,
its image in $@EMS$ is also a coproduct cocone.
We prove the latter claim by establishing the universal property of coproducts.
Existence follows from the fact that the quotient functor is full and surjective on objects.
Uniqueness amounts to showing that $f∘ι_i≈g∘ι_i$ for all~$i$ implies $f≈g$,
which follows from the definition of~$≈$ and the construction of negligible sets in the coproduct.

\subsection Measures on enhanced measurable spaces

\proclaim Definition.
A (complex infinite) ^={measure[|s]} on an ^{enhanced measurable space} $(X,M,N)$
is a map $μ:M'→`C$,
where $M'⊂M$ is an ideal of~$M$ such that $N⊂M'$,
the map $μ$ vanishes on~$N$,
and the ideal~$M'$ satisfies the following saturation condition:
for any countable family $\{m_i\}_{i∈I}$ of disjoint elements of~$M'$
such that the sum $∑_{i∈I} μ(m'_i)$ converges absolutely for any family $\{m'_i\}_{i∈I}$ satisfying $m'_i⊂m_i$ and $m'_i∈M'$,
we have $$⋃_{i∈I} m_i∈M',\qquad
μ\Biggl(⋃_{i∈I} m_i\Biggr)=∑_{i∈I} μ(m_i).$$
We say that $μ$ is a
\li ^={faithful measure[|s]} if $μ|_m=0$ implies $m∈N$ for any $m∈M'$;
\li ^={finite measure[|s]} if $M'=M$;
^^={faithful finite measure[|s]}
\li ^={semifinite measure[|s]} if for any $m∈M∖N$ there is $m'∈M'∖N$ such that $m'⊂m$;
^^={faithful semifinite measure[|s]}
\li ^={real measure} (alias ^={signed measure} or ^={charge}) if $μ$ factors through $`R→`C$;
\li ^={positive measure} if $μ$ factors through $[0,∞)→`C$.
\endlist
Here $μ|_m$ for $m∈M$ denotes the restriction of~$μ$ to $M'_m=\{m'∈M'\mid m'⊂m\}$,
which is a ^{measure} on the ^{induced enhanced measurable space} $(m,M_m,N_m)$.
Also, $0$ denotes the zero measure, with $M'=M$ and $μ(m)=0$ for all $m∈M'$.
Two ^{measures} are equal if they have the same ideal $M'$ and take the same values on all elements of~$M'$.

\proclaim Remark.
The definition of a ^{finite measure} can be simplified as follows:
a ^{finite measure} on an ^{enhanced measurable space} $(X,M,N)$ is a countably additive map
$M→`C$ that vanishes on~$N$.
Countable additivity means that for any countable family $\{m_i\}_{i∈I}$ of disjoint elements of~$M$
the sum $∑_{i∈I} μ(m_i)$ converges absolutely and $μ\bigl(⋃_{i∈I} m_i\bigr)=∑_{i∈I} μ(m_i)$.

\proclaim Remark.
By the Radon–Nikodym theorem, complex ^{semifinite measures} on a ^{localizable enhanced measurable space} $(X,M,N)$
form a free module of rank~1 over the complex algebra of all morphisms (in the category~$@EMS$) of the form $(X,M,N)→`C_~Borel$.
This necessitates the consideration of ^{semifinite measures} and not just ^{finite measures}.
Some definitions of infinite signed measures found in the literature
require that either the positive or negative part in the Jordan decomposition is finite.
Such a convention would preclude the strong version of the Radon–Nikodym given above,
which motivates our version.

\proclaim Remark.
The pushforward of a ^{semifinite measure} along a morphism in $@EMS$ need not be semifinite,
which requires us to also consider nonsemifinite measures in the above definition.

\label\traditionalmeasurespace
\proclaim Remark.
If $μ$ is a ^{positive measure} on an ^{enhanced measurable space} $(X,M,N)$,
then setting $μ(m)=∞$ for every $m∈M∖M'$ produces a countably additive map $M→[0,∞]$.
Thus, $(X,M,μ)$ is a measure space in the traditional sense.
Conversely, if $(X,M,μ)$ is a complete measure space in the traditional sense,
then setting $N=\{m∈M\mid μ(m)=0\}$ and $M'=\{m∈M\mid μ(m)≠∞\}$
produces an ^{enhanced measurable space} $(X,M,N)$ together with a positive ^{faithful measure} $μ|_{M'}$.

\subsection Localizable enhanced measurable spaces

The following definition was introduced (in the nonenhanced case) by Irving Segal [\EqMS, Definition~2.6].
The choice of terminology is motivated by Theorem~5.1 in the cited paper,
which proves that being localizable is equivalent to the Lebesgue decomposition property:
any σ-ideal~$I⊂M$ is {\it localized\/} at some measurable set~$P$
in the sense that $Q∈I$ if and only if $Q∖P$ is negligible.

\proclaim Definition.
A ^={localizable enhanced measurable space[|s]}
is an ^{enhanced measurable space} $(X,M,N)$
such that the ^{Boolean algebra}~$M/N$ is ^{localizable} (^!{localizable}).
The full subcategory of localizable objects in~$@EMS$ is denoted by $@=LEMS$.

Localizable measure spaces as defined by Segal [\EqMS] can be now identified
with ^{localizable enhanced measurable spaces} equipped with a positive ^{faithful semifinite measure}.

\proclaim Remark.
It is useful to reformulate the existence of suprema in $M/N$ without referring to the quotient ^{Boolean algebra}.
Suppose $\{m_i\}_{i∈I}$ is a family of elements of~$M$
and $\hat m∈M$.
We have $[\hat m]=\sup_{i∈I}[m_i]$
if and only if $m_i∖\hat m∈N$ for all $i∈I$
and whenever for some $\hat m'∈M$ we have $m_i∖\hat m'∈N$ for all $i∈I$,
then $\hat m∖\hat m'∈N$.
The set $\hat m$ is also known as the ^={essential suprem[um|a]} of the family $\{m_i\}_{i∈I}$.

\proclaim Remark.
Dixmier [\HS, §7 starting from Lemma~8]
constructs a ^{complete Boolean algebra}~$A$ whose ^{Stone spectrum}
is a ^{Stonean space}~$S$ in which every ^{meager subset} is ^{rare}
and the support of every measure is ^{rare}.
The ^{enhanced measurable space} $(X,M,N)="TM(S)$
is such that the ^{Boolean algebra} $M/N≅A$ is ^{complete}
and every measure on $(X,M,N)$ vanishes.
Thus, the ^{Boolean algebra} $A$ is ^{complete} but not ^{localizable}.
Therefore, we must include the requirement of existence of sufficiently many ^{continuous valuations} on~$M/N$ in the definition of ^{localizability}.
Segal [\EqMS] imposes this requirement implicitly, by virtue of working with a prespecified semifinite measure $m:M→[0,∞]$
whose class of measure~0 sets coincides with our~$N$.

\subsection Essential measures on enhanced measurable spaces

The goal of this section is to explain the “correct” definition of a measure
on an ^{enhanced measurable space} that is not ^{σ-finite}.
Here “correct” means, in particular, that \vmeasurevaluation\ holds.

\proclaim Definition.
An ^{enhanced measurable space} $(X,M,N)$ is ^={σ-finite}
^^={σ-finite enhanced measurable space[|s]}
^^={σ-finiteness}
if it admits a ^{faithful finite measure}.
An element $m∈M$ is ^{σ-finite} if the ^{induced enhanced measurable space} on~$m$ is ^{σ-finite}.

\proclaim Remark.
If $μ$ is a σ-finite measure on a measurable space $(X,M)$ (in the usual sense,
meaning $X$ is a countable union of elements of~$M$ with a finite $μ$-measure),
then the underlying ^{enhanced measurable space} in the sense of \vtraditionalmeasurespace\ is a ^{σ-finite enhanced measurable space}.
Indeed, in this case $X$ can be presented as a countable disjoint union $⋃_{i≥0}m_i$ of elements of~$M$ with a finite nonzero $μ$-measure,
and setting $ν(m')=∑_{i≥0}μ(m'∩m_i)/(2^i μ(m_i))$ produces a finite measure on $(X,M)$ with the same class of negligible sets.
Conversely, if $μ$ is a ^{faithful measure} on the ^{enhanced measurable space} $(X,M,N)$,
then the measure space constructed in \vtraditionalmeasurespace\ is σ-finite in the usual sense.

\proclaim Definition.
A (complex infinite) ^{measure} $μ:M'→`C$ on an ^{enhanced measurable space} $(X,M,N)$
is ^={essential}
^^={essential measure[|s]}
^^={essentiality}
if for any $m∈M$ such that $μ|_m≠0$ we can find a ^{σ-finite} $m'∈M$ such that
$m'⊂m$ and $μ|_{m'}≠0$.

\label\essentialcompletelyadditive
\proclaim Proposition.
Any ^{essential} (complex infinite) ^{measure} $μ:M'→`C$ on an ^{enhanced measurable space} $(X,M,N)$
is completely additive:
^^={completely additive measure[|s]}
for any family $\{m_i\}_{i∈I}$ (countable or not) of elements of~$M'$
disjoint up to a ^{negligible set} (meaning $m_i∩m_j∈N$ whenever $i≠j$)
that admits an ^{essential supremum} $S∈M$ and
such that the sum $∑_{i∈I} μ(m'_i)$ converges absolutely for any family $\{m'_i\}_{i∈I}$ satisfying $m'_i⊂m_i$ and $m'_i∈M'$,
we have $$S∈M',\qquad
μ(S)=∑_{i∈I} μ(m_i).$$
Conversely, if $(X,M,N)$ is ^{localizable}, any completely additive measure is ^{essential}.

\proof Proof.
Suppose $μ$ is an ^{essential measure} and $\{m_i\}_{i∈I}$ is a family with indicated properties.
The absolute convergence of $∑_{i∈I} μ(m'_i)$ implies that the set $J=\{i∈I\mid μ|_{m_i}≠0\}$ is countable.
Set $S'=⋃_{j∈J}m_j∈M$.
By assumption, we have $S'∈M'$ and $μ(S')=∑_{i∈I} μ(m_i)$.
It remains to show that $S∖S'∈M'$ and $μ(S∖S')=0$.
Indeed, $S∖S'$ is the ^{essential supremum} of $\{m_i∩(S∖S')\}_{i∈I∖J}$
and the restriction of~$μ$ to any of these sets vanishes.
Hence $μ$ vanishes on all ^{σ-finite} elements of~$M$,
therefore $μ|_{S∖S'}=0$ because $μ$ is essential.
\ppar
Conversely, suppose $μ$ is completely additive.
If $μ$ is not ^{essential},
then there is $m∈M$ such that $μ|_{m'}=0$ for any ^{σ-finite} $m'∈M$ such that $m'⊂m$.
Since the ^{induced enhanced measurable space} $(m,M_m,N_m)$ is ^{localizable},
we can find a family $\{m_i\}_{i∈I}$
of ^{σ-finite} elements of~$M_m$
disjoint up to a ^{negligible set}
whose ^{essential supremum} equals~$m$.
By completely additivity of~$μ$ we must have $μ(m)=0$, a contradiction.

\label\essentialmeasuresuprema
\proclaim Proposition.
Any ^{essential} (complex infinite) ^{measure} $μ:M'→`C$ on an ^{enhanced measurable space} $(X,M,N)$
preserves ^{essential suprema}:
for any directed subset $T⊂M'$ that is closed under passage to measurable subsets
and admits an ^{essential supremum} $S∈M'$ (^!{essential supremum}),
the limit $\lim_{t∈T}μ(t)$ exists and $μ(S)=\lim_{t∈T}μ(t)$.
Conversely, if $(X,M,N)$ is ^{localizable}, then any ^{measure} that preserves ^{essential suprema} is ^{essential}.

\proof Proof.
Suppose $μ$ is an ^{essential measure}
and $T$ is a directed subset of~$M'$ that is closed under passage to measurable subsets and admits an ^{essential supremum} $S∈M'$ (^!{essential supremum}).
We want to show that the limit $\lim_{t∈T}μ(t)$ exists and $μ(S)=\lim_{t∈T}μ(t)$.
Replacing $(X,M,N)$ by $(S,M_S,N_S)$, we can assume $μ$ to be finite, i.e., $M'=M$, and $X=S$ is now the ^{essential supremum} of~$T$.
By linearity, we can assume $μ$ to be real.
By the classical Hahn–Jordan decomposition (Fremlin [\MTii, Theorem~231E]),
we can further assume $μ$ to be positive and finite.
In particular, $\lim_{t∈T}μ(t)=\sup_{t∈T}μ(t)$ exists and $\sup_{t∈T}μ(t)≤μ(S)$.
\ppar
For any $n>0$, set $A_n$ to an element of~$T$ such that $\sup_{t∈T}μ(t)-μ(A_n)<1/n$ and $A_n⊃A_{n-1}$ if $n>1$.
Denote by $F$ the union $⋃_{n>0}A_n$.
By construction, $F∈M$ and for any $t∈T$ we have $μ(t∖F)<1/n$ for all $n>0$, so $μ(t∖F)=0$.
Furthermore, $S∖F$ is the ^{essential supremum} of $t∖F$ for all $t∈T$.
Replacing $(X,M,N)$ by $(S∖F,M_{S∖F},N_{S∖F})$, we may now assume $μ(t)=0$ for all $t∈T$,
and we have to show that $μ(X)=0$.
Any ^{σ-finite} subset $B⊂X$ is the ^{essential supremum} of $t∩B$ for all $t∈T$.
Since all ^{measures} on ^{σ-finite enhanced measurable spaces} are ^{essential}, this implies that $μ(B)=0$.
Thus, $μ(X)=0$ because $μ$ is ^{essential}.
\ppar
Conversely, suppose $μ$ is a ^{measure} that preserves ^{essential suprema}.
To show that $μ$ is ^{essential},
suppose $m∈M$ is such that $μ|_{m'}=0$ for all ^{σ-finite} $m'∈M$ such that $m'⊂m$.
Since $(X,M,N)$ is ^{localizable}, $m$ is the ^{essential supremum} of all such~$m'$.
This ^{essential supremum} is preserved by~$μ$, therefore, $μ(m)=0$.

\label\essentialradonnikodym
\proclaim Remark.
If $(X,M,N)$ is an ^{enhanced measurable space}
with a ^{faithful semifinite measure}~$μ$,
then a map $ν:M→`R$ is a (finite real) ^{measure} on $(X,M,N)$ if and only if $ν$ is countably additive and absolutely continuous with respect to~$μ$.
Likewise, a map $ν:M→`R$ is an ^{essential measure} on $(X,M,N)$ if and only if $ν$ is countably additive and truly continuous (Fremlin [\MTii, 232A(b)])
with respect to~$μ$.
Thus, ^{essential measures} are precisely those measures that admit a Radon–Nikodym derivative with respect to~$μ$ (Fremlin [\MTii, 232E)).

\label\essentialcounterexample
\proclaim Remark.
Countably additive measures on ^{localizable enhanced measurable spaces} that are not ^{essential}
exist if and only if real-valued-measurable cardinals exist (Fremlin [\MTiii, Theorem~363S]).
Indeed, if $X$ is such a cardinal with a probability measure $2^X→[0,1]$ that vanishes on all countable subsets of~$X$,
then $2^X$ is a ^{localizable Boolean algebra}
such that $μ$ induces a probability measure~$ν$ on $"Spec(2^X)$ that vanishes on all ^{rare} subsets
but does not preserve the essential supremum~$X$ of the family of singleton subsets $\{\{x\}∈2^X\mid x∈X\}$,
since $ν(\{x\})=0$, but $ν(X)=1$.
In particular, $ν$ is not a ^{normal measure} because it is not a Radon measure
due to its violation of the τ-smoothness property.

\label\measurevaluation
\proclaim Proposition.
Given a ^{localizable enhanced measurable space} $(X,M,N)$,
composition with the quotient map $M→M/N$
establishes a canonical bijective correspondence
from ^{continuous valuations} on~$M/N$
to (finite) ^{essential measures} on $(X,M,N)$.

\proof Proof.
The quotient map $M→M/N$ sends countable unions to countable joins.
A ^{valuation} on $M/N$ restricts to a finitely additive measure on~$M$ that vanishes on~$N$,
and this correspondence is bijective.
A ^{valuation} is ^{continuous} if it sends suprema of directed subsets of $M/N$ to limits in~$`C$,
which translates into mapping ^{essential suprema} of directed subsets of~$M$ to limits in~$`C$.
By \vessentialmeasuresuprema, the latter measures are precisely the ^{essential measures}.

\label\sigmafiniteessential
\proclaim Remark.
By definition, any ^{measure}~$μ$ on a ^{σ-finite enhanced measurable space} $(X,M,N)$ is ^{essential},
so induces a ^{continuous valuation} on~$M/N$ by \vmeasurevaluation.

\subsection Strictly localizable enhanced measurable spaces

The problems described in \vstrictloc\ justify our consideration of ^{strictly localizable enhanced measurable spaces} below.

\proclaim Definition.
An ^{enhanced measurable space} $(X,M,N)$ is ^={strictly localizable}
(alias ^={decomposable})
^^={strictly localizable enhanced measurable space[|s]}
^^={strict localizability}
if it is isomorphic in the category $@PreEMS$ to the coproduct of a small family of ^{σ-finite enhanced measurable spaces}.
A specific choice of such a family will be referred to as a ^={strictly localizing partition[|s]}.
The ^={category of strictly localizable enhanced measurable spaces} $@=SLEMS$
is the full subcategory of $@EMS$ consisting of strictly localizable ^{enhanced measurable spaces}.

\label\strictlocnotiso
\proclaim Example.
(See Fremlin [\MTii, 216D].)
^{Strictly localizable enhanced measurable spaces} are not closed under isomorphisms in the category $@EMS$.
Pick an uncountable set~$I$ and consider the ^{enhanced measurable space} $(X,M,N)$,
where $X=\{0,1\}⨯I$, $N$ is the set of countable subsets of $\{1\}⨯I$,
and $M$ consists of subsets of the form $(\{0,1\}⨯J)⊕n$, where $J⊂I$ and $n∈N$.
The projection map $p:\{0,1\}⨯I→I$
is an isomorphism of ^{enhanced measurable spaces} $(X,M,N)→(I,2^I,\{∅\})$ in the category $@EMS$.
The inverse map injects $I≅\{0\}⨯I$ into $\{0,1\}⨯I$.
The codomain $(I,2^I,\{∅\})$ is ^{strictly localizable},
but the domain $(X,M,N)$ is not.
Indeed, elements of any ^{strictly localizing partition} of $(X,M,N)$ must be countable subsets of~$X$.
But then the set $\{0\}⨯I$ has a measurable intersection with any element of such a partition,
and yet $\{0\}⨯I∉M$ because $I$ is uncountable.

The following definition will eliminate \vstrictlocnotiso.

\proclaim Definition.
The category $@=LDEMS$ of ^={locally determined enhanced measurable space[s|]}
^^={locally determined}
is the full subcategory of $@EMS$ consisting of ^{enhanced measurable spaces} $(X,M,N)$
such for any $A⊂X$,
if $A∩F∈M$ (respectively $A∩F∈N$) for all ^{σ-finite} $F∈M$,
then also $A∈M$ (respectively $A∈N$).

\label\strictlocdet
\proclaim Remark.
We have an inclusion of full subcategories $@SLEMS⊂@LDEMS$ by Fremlin [\MTii, Theorem~211L(d)].

\proclaim Remark.
Consider an ^{enhanced measurable space} $(X,M,N)$.
Adding to $M$ respectively $N$ the subsets of~$X$ that ^!{locally determined} says must be elements of~$M$ respectively~$N$
produces a ^{locally determined enhanced measurable space} $(X,M',N')$ with $M'⊃M$ and $N'⊃N$ (Fremlin [\MTii, Proposition~213D]).
Thus, the map $(X,M',N')→(X,M,N)$ induced by the identity map on~$X$
is a ^{morphism of enhanced measurable spaces}.
If $(X,M,N)$ is ^{localizable}, then the induced ^{homomorphism of Boolean algebras} $M/N→M'/N'$ is an isomorphism (Fremlin [\MTii, Proposition~213H(b)]).
However, the identity map on~$X$ does not induce a ^{morphism of enhanced measurable spaces} $(X,M,N)→(X,M',N')$ since $N'$ can be bigger than~$N$.
Arguably, we want the identity map to be an isomorphism nevertheless.
This can be achieved by modifying the definition of a ^{premap of enhanced measurable spaces} (^!{premap of enhanced measurable spaces}),
replacing the condition $f^*n∈N_X$ for all $n∈N_Y$
with the condition that $f^*n∩m∈N_X$ for all $n∈N_Y$ and ^{σ-finite} $m∈M_X$.
We did not adopt this modification for the same reason that we required ^{completeness} in the definition of ^{enhanced measurable spaces}:
the resulting definition of a ^{morphism of enhanced measurable spaces} is the simplest possible.
(Compare ^!{enhanced measurable space} with ^!{completeness} for the noncomplete case.)
Since we quickly restrict to ^{strictly localizable enhanced measurable spaces},
which are ^{locally determined} by \vstrictlocdet,
the difference is irrelevant for us anyway.
The category $@LDEMS$ is only used in purely expository \vstrictlocdetiso, \vcslemsdetiso, and \vstrictloc,
which are not used in the main results.

\label\strictlocdetiso
\proclaim Lemma.
^{Strictly localizable enhanced measurable spaces} are closed under isomorphisms in the category $@LDEMS$.

\proof Proof.
Suppose $f:(X,M_X,N_X)→(Y,M_Y,N_Y)$ is an isomorphism in $@LDEMS$ with a ^{strictly localizable} codomain.
Fix a ^{strictly localizing partition} $\{B_i\}_{i∈I}$ for $(Y,M_Y,N_Y)$
together with a finite measure $μ_i$ (for each $i∈I$) on $(Y,M_Y,N_Y)$ such that for any $m∈M_Y$ we have $μ_i(m)=0$ if and only if $m∩B_i∈N_Y$.
The sets $A_i=f^*B_i$ and (if nonempty) $X∖?pdom f$ form a partition of~$X$ into nonempty measurable subsets of~$X$.
For any $i∈I$ the pushforward measure $ν_i≔(f^{-1})_*(μ_i)$ satisfies the following condition: for any $m∈M_X$ we have $ν_i(m)=0$ if and only if $m∩A_i∈N_X$.
We claim that the partition $\{A_i\}_{i∈I}∪\{X∖?pdom f\}$ is a ^{strictly localizing partition} for $(X,M_X,N_X)$.
This follows from Fremlin [\MTii, Proposition~213O(a)].

\label\strictlocloc
\proclaim Proposition.
^{Strictly localizable enhanced measurable spaces} are ^{localizable}.

\proof Proof.
The Boolean algebra $M/N$ associated to a coproduct of ^{enhanced measurable spaces} in the category $@PreEMS$
is the product of the Boolean algebras associated to individual members of the family.
Since the inclusion $@LBAlg→@BAlg$ preserves products,
it suffices to show that ^{σ-finite enhanced measurable spaces}
are ^{localizable}.
By definition of a ^{σ-finite enhanced measurable space},
it admits a ^{faithful finite measure}~$μ$,
which descends to a faithful finite ^{continuous valuation}~$ν$ on~$M/N$ by \vsigmafiniteessential.
The ^{Boolean algebra} $M/N$ admits suprema of countable subsets
because the ^{Boolean algebra}~$M$ does
and the ideal $N⊂M$ is closed under suprema of countable subsets.
By \vbalgcountable, the ^{Boolean algebra} $M/N$ is complete.
Thus, $M/N$ belongs to $@LBAlg$.

\label\balgcountable
\proclaim Lemma.
Suppose a ^{Boolean algebra}~$A$ admits a faithful ^{continuous valuation} $ν:A→[0,∞)$.
If $A$ admits suprema of countable subsets, then $A$ is complete.
In addition, if $A'$ is a ^{complete Boolean algebra}
and $f:A→A'$ is a ^{homomorphism of Boolean algebras} that preserves countable suprema,
then $f$ is complete.

\proof Proof.
Suppose $S⊂A$ is an arbitrary subset of~$A$, which we may assume to be closed under finite joins, which does not alter suprema.
We want to show that $\sup S$ exists.
It suffices to construct a countable subset $T⊂S$ such that for any $s∈S$ we have $s≤\sup T$,
which implies that $\sup S=\sup T$ exists.
\ppar
Set $r=\sup_{s∈S}ν(s)$, which exists because $ν(s)≤ν(1)$.
Choose a countable subset $T⊂S$ (which we again may assume to be closed under finite joins) such that $r=\sup_{t∈T}ν(t)$ and set $τ=\sup T$, so that $ν(τ)≥r$.
\ppar
For any $s∈S$ we have $$ν(s∨τ)=ν(s∨\sup T)=ν(\sup_{t∈T}(s∨t))≤ν(\sup S)≤r$$ because $s∨t∈S$ and $ν(\sup S)=\sup_{s∈S}ν(s)=r$ by continuity of~$ν$.
Thus, $$ν(s∖τ)=ν((s∨τ)∖τ)=ν(s∨τ)-ν(τ)≤r-ν(τ)≤0,$$
so $ν(s∖τ)=0$ and $s∖τ=0$ because $ν$ is faithful.
Hence, $s≤τ=\sup T$ as desired.
\ppar
To show that the map $f$ preserves arbitrary suprema, consider an arbitrary subset $S⊂A$,
which we may assume to be directed.
Using the above construction, choose a countable directed subset $T⊂S$ such that $\sup S=\sup T$.
Now $f(\sup S)=f(\sup T)=\sup f_!T≤\sup f_!S$.
Also, $\sup f_!S≤f(\sup S)$ follows from $f(s)≤f(\sup S)$ for any $s∈S$.
Hence $f(\sup S)=\sup f_!S$, i.e., $f$ preserves arbitrary suprema.

\subsection Compact enhanced measurable spaces

The following definition was introduced (with an additional countability condition) by Edward Marczewski [\CMeas, §4].
An extensive treatment of the underlying theory was given by Pfanzagl–Pierlo [\CSS].

\proclaim Definition.
A collection $K⊂2^X$ of subsets of a set~$X$ is a ^={compact class} if for any $K'⊂K$ the following ^={finite intersection property} holds:
if for any finite $K''⊂K'$ we have $⋂K''≠∅$, then also $⋂K'≠∅$.
A ^={compact enhanced measurable space[|s]}
^^={compactness}
is an ^{enhanced measurable space} $(X,M,N)$
for which there is a compact class $K⊂M$ such that for any $m∈M∖N$ there is $k∈K∖N$ such that $k⊂m$.

\proclaim Remark.
The choice of terminology is motivated by the fact (Fremlin [\MTiii, Lemma~342D(a)]) that $K⊂2^X$
is a ^{compact class} if and only if there is a compact topology on the set~$X$
such that every element of~$K$ is a closed (and hence compact) subset of~$X$.
Thus, an ^{enhanced measurable space} is ^{compact} if and only if it admits
a compact topology such that every nonnegligible measurable set contains a nonnegligible subset that is closed (Fremlin [\MTiii, Corollary~342F]).

\proclaim Definition.
The ^={category of compact strictly localizable enhanced measurable spaces} $@=CSLEMS$
^^={compact strictly localizable enhanced measurable space[|s]}
is the full subcategory of $@EMS$ consisting of ^{enhanced measurable spaces}
that are ^{compact[| enhanced measurable space]} (^!{compact enhanced measurable space})
and ^{strictly localizable} (^!{strictly localizable enhanced measurable space}).

\proclaim Example.
The class of ^{compact strictly localizable enhanced measurable spaces} contains ^={Radon enhanced measurable space[s|]}, i.e., ^{enhanced measurable spaces}
equipped with a structure of a Hausdorff topological space such that
the ^{σ-algebra} of measurable sets contains open sets and there is a ^{faithful measure}~$μ$
that is locally finite (every point has a neighborhood of finite measure)
and inner regular with respect to compact subsets (the measure of any measurable subset is the supremum of measures of its compact subsets).
(See Fremlin [\MTiv, Definition~411H(b)].)
Thus, $@CSLEMS$ includes the vast majority of ^{enhanced measurable spaces} used in analysis.
\ppar
Indeed, any Radon ^{enhanced measurable space} is ^{strictly localizable} by Fremlin [\MTiv, Theorem~414J].
Also, compact subsets of a Hausdorff topological space form a ^{compact class} because compact subsets of Hausdorff spaces are closed.
Furthermore, the inner regularity property implies that for any $m∈M∖N$ we have $μ(m)=\sup_{k\subset m}μ(k)>0$,
so there is $k\subset m$ such that $k$ is compact and $k∈K∖N$, as required by the definition of a ^{compact enhanced measurable space}.
Slightly more generally, compactness continues to hold for non-Hausdorff Radon ^{enhanced measurable spaces},
which are defined using closed compact subsets instead of compact subsets.

\label\compactsubspace
\proclaim Lemma.
If $(X,M,N)$ is a ^{compact enhanced measurable space} and $m∈M$,
then the ^{induced enhanced measurable space} $(m,M_m,N_m)$ is also ^{compact}.

\proof Proof.
If $K⊂M$ is a ^{compact class} that exhibits the ^{compactness} of $(X,M,N)$,
then $K∩M_m$ is a ^{compact class} that exhibits the ^{compactness} of $(m,M_m,N_m)$.

\label\vcslemsnotiso
\proclaim Example.
^{Compact strictly localizable enhanced measurable spaces} are not closed under isomorphisms in the category~$@EMS$.
Indeed, \vstrictlocnotiso\ constructs an isomorphism~$p$ of ^{enhanced measurable spaces}
whose domain is not ^{strictly localizable} and codomain is ^{strictly localizable}.
Both domain and codomain are ^{compact}, as witnessed by the compact class $\{\{0\}⨯\{i\}\}_{i∈I}$ for the domain
and $\{\{i\}\}_{i∈I}$ for the codomain.

The proof of the following lemma resembles
Part~($γ$) of the proof of Theorem in~\S3 of Fremlin [\CMS]
or Part~(g) of the proof of Theorem~343B in Fremlin [\MTiii].

\label\cslemsdetiso
\proclaim Lemma.
^{Compact strictly localizable enhanced measurable spaces} are closed under isomorphisms in the category $@LDEMS$ (^!{locally determined}).

\proof Proof.
Suppose $f:(X,M_X,N_X)→(Y,M_Y,N_Y)$ is an isomorphism in $@EMS$
whose codomain is a ^{compact strictly localizable enhanced measurable space}.
By \vstrictlocdetiso, the ^{enhanced measurable spaces} $(X,M_X,N_X)$ is ^{strictly localizable} because it is ^{locally determined}.
If each element of a ^{strictly localizing partition} of $(X,M,N)$ is compact, then the union of compact classes of each part is a compact class for $(X,M_X,N_X)$.
Thus we may assume $(Y,M_Y,N_Y)$ to be ^{σ-finite}.
\ppar
Denote by $g=f^{-1}:(Y,M_Y,N_Y)→(X,M_X,N_X)$ the inverse of~$f$ in the category $@EMS$.
Observe that $M_X=N_X$ if and only if $M_Y=N_Y$, in which case $(X,M_X,N_X)$ is compact because $M_X∖N_X=∅$.
Otherwise, we can always choose point-set representatives for $f$ and $g$
so that their ^{point-set domains} have empty complements, i.e., both $f$ and $g$ are everywhere defined.
Observe that $f^*(M_Y∖N_Y)⊂M_X∖N_X$ and $g^*(M_X∖N_X)⊂M_Y∖N_Y$.
\ppar
Suppose $K_Y⊂M_Y$ is a ^{compact class} that exhibits the compactness of~$Y$.
If $K'_Y⊃K_Y$ is a bigger ^{compact class}, then $K'_Y$ also exhibits the compactness of~$Y$.
In particular, we may assume $K_Y$ to be closed under finite unions and countable intersections.
Take $K_X$ to be the set of $k_X∈M_X$ such that there is $k_Y∈K_Y$ for which $g_!k_Y⊂k_X⊂f^*k_Y$.
We fix such $k_Y$ for each $k_X∈K_X$ and denote it by $L_{k_X}$.
We claim that the set $K_X$ is a ^{compact class} that exhibits the compactness of~$(X,M_X,N_X)$.
\ppar
To show that $K_X$ is a ^{compact class}, suppose that $K'⊂K_X$ has the ^{finite intersection property},
i.e., for any finite $K''⊂K'$ we have $⋂K''≠∅$.
We have $f^*(⋂_{k∈K''}L_k)=⋂_{k∈K''}f^*(L_k)⊃⋂_{k∈K''}k≠∅$, so $⋂_{k∈K''}L_k≠∅$.
Thus, the family $\{L_k\}_{k∈K'}$ has the ^{finite intersection property} and $W=⋂_{k∈K'}L_k≠∅$.
Now $g_!W⊂g_!(L_k)⊂k$ for all $k∈K'$,
so $g_!W⊂⋂_{k∈K'}k$, and hence $⋂_{k∈K'}k≠∅$, i.e., $K_X$ is a ^{compact class}.
\ppar
We now show that $K_X$ exhibits the compactness of~$X$.
Fix a ^{faithful finite measure}~$μ$ on $(X,M_X,N_X)$ and its pushforward~$ν=f_*μ$ on $(Y,M_Y,N_Y)$.
Given $m∈M_X∖N_X$, we must demonstrate that there is $k_X∈K_X∖N_X$ such that $k_X⊂m$.
To this end, we also construct $k_Y∈K_Y∖N_Y$ such that $g_!k_Y⊂k_X⊂f^*k_Y$.
Fix some $γ>0$ such that $γ<μ(m)$.
\ppar
Construct sequences $\{V_n\}_{n≥0}$, where $V_n∈M_Y∖N_Y$, $ν(V_n)>γ$,
and $\{F_n\}_{n≥0}$, where $F_n∈M_X∖N_X$, $μ(F_n)>γ$
by induction as follows.
Set $F_0=m$.
Suppose we already constructed $F_k$ for $k≤n$ and $V_k$ for $k<n$.
Observe that $g^*F_n∈M_Y∖N_Y$ because $ν(g^*F_n)=μ(F_n)>γ$.
Consider the set~$Q$ of all $R∈K_Y∖N_Y$ such that $R⊂g^*F_n$ and $ν(R)>0$,
which is nonempty by compactness of $(Y,M_Y,N_Y)$ exhibited by~$K_Y$
and is closed under finite unions by definition of~$K_Y$.
Denote by $P$ the ^{essential supremum} of~$Q$.
We have $ν(P)=ν(g^*F_n)$, since otherwise $g^*F_n∖P∈M_Y∖B_Y$
would contain some $R∈K_Y∖N_Y$ such that $R⊂g^*F_n∖P$ and $ν(R)>0$,
and hence $R∈Q$, contradicting the definition of~$P$.
Now set $V_n$ to an element $R∈Q$ such that $ν(R)>γ$.
Such an element exists because $ν(P)=ν(g^*F_n)=μ(F_n)>γ$ and $ν$ is an ^{essential measure}.
\ppar
Set $F_{n+1}$ to $f^*V_n∖(f^*g^*F_n∖F_n)$.
Observe that $f^*g^*F_n∖F_n∈N_X$, so $μ(F_{n+1})=μ(f^*V_n)=ν(V_n)>γ$.
Also $F_{n+1}⊂F_n$ by construction.
\ppar
After the induction, set $k_X=⋂_{n≥0}F_n$ and $k_Y=⋂_{n≥0}V_n$.
We have $k_Y∈K_Y$ because $V_n∈K_Y$ by construction and $K_Y$ is closed under countable intersections.
We have $$k_Y=⋂_{n≥0}V_n⊂⋂_{n≥0}g^*F_n=g^*⋂_{n≥0}F_n=g^*k_X,$$ so $g_!k_Y⊂k_X$.
Also $$k_X=⋂_{n≥0}F_n⊂⋂_{n≥0}F_{n+1}⊂⋂_{n≥0}f^*V_n=f^*⋂_{n≥0}V_n=f^*k_Y.$$
Altogether, $$g_!k_Y⊂k_X⊂f^*k_Y,$$
which shows that $k_X∈K_X$.
Furthermore, $μ(k_X)=\inf_{n≥0}μ(F_n)≥γ>0$, so $k_X∈K_X∖N_X$.
By construction, $k_X⊂F_0=m$, which completes the proof.

\proclaim Remark.
The argument in the last paragraph of the above proof
heavily uses the existence of a ^{faithful finite measure} obtained by exploiting ^{strict localizability}.
Although it is not necessary below,
it would be interesting to know whether ^{compact enhanced measurable spaces}
are closed under isomorphisms in $@LDEMS$,
and if not, whether one can impose a condition weaker than ^{strict localizability}
that would guarantee such closedness.

\proclaim Proposition.
Suppose $(X,M,N)$ is a ^{compact enhanced measurable space} that admits a ^{faithful semifinite measure},
$(X',M',N')$ is a ^{strictly localizable enhanced measurable space},
and $f:(X,M,N)→(X',M',N')$ is a ^{morphism of enhanced measurable spaces}.
Then $f$ has a well-defined ^={measurable image[|s]} $m∈M'$ in the following sense:
$f^*(X∖m)∈N$ and
for any $m'∈M'$ such that $m'⊂m$ and $f^*m'∈N$ we have $m'∈N'$.
The class of~$m$ in $M'/N'$ is unique.

\proof Proof.
Uniqueness is easy to show: if $m_0$ and $m_1$ are two elements with such properties,
then $f^*(m_0∖m_1)⊂f^*(X∖m_1)∈N$, so $m_0∖m_1∈N'$.
Likewise, $m_1∖m_0∈N'$ and hence $m_0⊕m_1∈N'$, so $m_0$ and $m_1$ map to the same element in $M'/N'$.
To show existence, denote by $p∈M'$ the ^{essential supremum}
of all $m'∈M'$ such that $f^*m'∈N$.
We claim that the element $m=X∖p∈M'$ has the desired properties.
First, if $m'⊂m$ and $f^*m'∈N$, then $m'∖p∈N'$.
Since $m'∖p=m'$ because $m'⊂m$ and $m∩p=∅$, we have $m'∈N'$.
\ppar
Secondly, to show that $f^*p∈N$, choose a ^{strictly localizing partition} $\{q_i\}_{i∈I}$ for $(X',M',N')$.
Set $r_i=q_i∩p$.
We claim that $f^*r_i∈N$.
Indeed, $r_i$ is the ^{essential supremum} of the collection $R_i=\{m∈M'\mid m⊂q_i, f^*m∈N\}$.
The ^{induced enhanced measurable space} of $q_i$ is ^{σ-finite},
so \vbalgcountable\ supplies a countable directed subset $R_i'⊂R_i$ with the same ^{essential supremum}.
Since $R_i'$ is countable, its ^{essential supremum} can be computed as the union.
Since $f^*$ preserves unions and $N$ is a ^{σ-ideal}, we deduce that $f^*r_i∈N$.
\ppar
Thus, the ^{induced enhanced measurable space} of $f^*p$ is a ^{compact enhanced measurable space}
that admits a ^{faithful semifinite measure}~$μ$
and is partitioned into a disjoint family of negligible subsets $\{f^*r_i\}_{i∈I}$.
Furthermore, by definition of a ^{strictly localizing partition},
for any subset $J⊂I$ we have $⋃_{i∈J}f^*r_i∈M$.
By Fremlin [\MTiv, Lemma~451Q], we have $⋃_{i∈I}f^*r_i=f^*p∈N$, which completes the proof.

\label\eqmeasspaceloc
\section Equivalence between measurable spaces and measurable locales

\subsection From hyperstonean spaces to enhanced measurable spaces

Recall that a continuous map of topological spaces is {\it open\/}
if images of open subsets are open.
The category of topological spaces and open maps is denoted by $@=TopOpen$.

The construction of an ^{enhanced measurable space} out of a topological space
given below was already known to
Loomis [\RepL] and Sikorski [\RepS], who defined the ^{enhanced measurable space} denoted below by $"Spec(A)$
for an arbitrary σ-complete Boolean algebra~$A$ by applying \vTMdef\ to the ^{Stone space} of~$A$,
and proved that its Boolean algebra of equivalences classes of measurable sets modulo negligible sets
is isomorphic to~$A$.

\label\TMdef
\proclaim Definition.
We define a functor $$"=TM:@TopOpen→@EMS$$
(TM for topology-to-measure) by sending a topological space $(X,U)$
to the ^{enhanced measurable space} $"TM(X,U)=(X,M,N)$ defined as follows.
Define $N$ to be the collection of all ^{meager subsets} of~$X$.
Define $M$ to be the collection of all subsets of~$X$ with the property of Baire, i.e., symmetric differences of elements of~$U$ and~$N$.
An open continuous map $f:(X,U)→(X',U')$ of topological spaces
is sent to the morphism $f:(X,M,N)→(X',M',N')$ with the same underlying map of sets.

\label\TMdefcorrect
\proclaim Lemma.
This definition is correct.

\proof Proof.
\vhyperstoneanmeasurable\ proves that $M$ is a ^{σ-algebra}
and $N⊂M$ is a ^{σ-ideal} of subsets of~$X$.
Thus, $(X,M,N)$ is an ^{enhanced measurable space}.
\ppar
Now suppose $f:S→S'$ is an open map
and $"TM(f):(X,M,N)→(X',M',N')$ is the induced ^{map of enhanced measurable spaces}.
We have to show that $"TM(f)^*$ sends elements of $M'$ to~$M$ and $N'$ to~$N$.
Since $f$ is continuous and $f^*$ preserves symmetric differences and countable unions,
it suffices to show that $f^*$ preserves closed subsets with empty interiors.
Indeed, if $V⊂X'$ is a closed subset with an empty interior
and $U⊂f^*V$ is an open subset,
then $f(U)⊂V$ is also an open subset because $f$ is an open map,
and hence $f(U)=∅$ and $U=∅$, so $f^*V⊂X$ is a closed subset with an empty interior.

\label\speccompactstrictloc
\proclaim Proposition.
The functor $$"TM:@TopOpen→@EMS$$
restricts to a functor $$"TM:@HStonean→@CSLEMS.$$

\proof Proof.
Recall that $@HStonean$ is a full subcategory of $@TopOpen$,
so it suffices to show that $"TM$ sends a ^{hyperstonean space}~$S$ to a ^{compact strictly localizable enhanced measurable space}.
To show strict localizability,
use Zorn's lemma to choose a maximal disjoint family $\{a_i\}_{i∈I}$ of nonzero elements
of the ^{complete Boolean algebra} $A="COpen("Ω(S))$
such that the ^{complete Boolean algebra} $a_iA$ admits a faithful finite ^{continuous valuation} for any $i∈I$.
In particular, every $a_i⊂S$ can be chosen to be a ^{clopen subset} of~$S$,
which admits a finite ^{measure} supported on it.
The ^{clopen subsets} $a_i⊂S$ are disjoint and their union is an open subset $b⊂S$
such that $¬¬b=S$, so $S∖b$ is ^{meager}.
The partition of~$S$ into disjoint ^{σ-finite} measurable subsets $\{a_i\}_{i∈I}$ and $S∖b$
exhibits $"TM(S)$ as a ^{strictly localizable enhanced measurable space}.
Indeed, if $\{n_i\}_{i∈I}$ is a collection of negligible (i.e., ^{rare}) subsets of~$a_i$,
then $⋃_{i∈I}n_i$ is ^{rare}, hence negligible.
Likewise, if $\{u_i⊕n_i\}_{i∈I}$ is a collection of measurable subsets of~$a_i$,
then $⋃_{i∈I}(u_i⊕n_i)=⋃_{i∈I}u_i⊕⋃_{i∈I}n_i$ is measurable because $⋃_{i∈I}u_i$ is open and $⋃_{i∈I}n_i$ is ^{rare}.
\ppar
To show ^{compactness},
recall that ^{hyperstonean spaces} are ^{compact}, so their ^{clopen subsets} are also ^{compact}.
Thus the class of all ^{clopen subsets} of~$S$ is a ^{compact class}.
It remains to show that any nonnegligible measurable subset $a⊂"TM(S)$ contains a nonempty ^{clopen subset}.
Since any open subset of $"TM(S)$ is a union of ^{clopen subsets}, it suffices to show that $a$ contains a nonempty open subset.
By definition of~$"TM$, the subset $a$ is the symmetric difference $u⊕n$ of an open subset~$u$ and a ^{meager subset}~$n$.
The subset~$u$ is nonempty because $a$ is nonnegligible.
It suffices to show that the subset~$n$ is ^{rare}, since then $u∖\bar n$ is a nonempty open subset of~$a$.
By definition of a ^{rare} set, it suffices to show that $\bar n$ is ^{rare}.
Denote by $\hat n$ the interior of $\bar n$.
The difference $\bar n∖\hat n$ is ^{rare}, so for any ^{normal measure}~$μ$ on~$S$
we have $0=μ(n)=μ(\bar n)=μ(\hat n)$, which by definition of a ^{hyperstonean space} implies $\hat n=∅$,
i.e., $n$ is ^{rare}.
This shows that $"TM(S)$ is a ^{compact enhanced measurable space}.

\label\defSpec
\proclaim Definition.
We define a functor $$"=Spec:@MLoc→@CSLEMS$$
as the composition $$@MLoc\lto7{"Ideal}@HStoneanLoc\lto7{"Sp}@HStonean\lto7{"TM}@CSLEMS.$$

\proclaim Proposition.
For any ^{hyperstonean space}~$S$
there is a canonical bijective correspondence
between ^{normal measures} on~$S$
and ^{essential measures} on $"TM(S)=(S,M,N)$.

\proof Proof.
Combine \vmeasurevaluation, \vnormalmeasurevaluation, and \vhnormal.
Alternatively, observe directly that any ^{essential measure} on $"TM(S)$
restricts to a ^{normal measure} on~$S$
(with the essentiality property implying τ-smoothness, which is equivalent to being a Radon measure)
and conversely, any ^{normal measure} on~$S$ extends to a unique
^{essential measure} on $"TM(S)$ because it automatically vanishes on elements of~$N$.

\subsection From enhanced measurable spaces to measurable locales

\proclaim Definition.
The functor $"=PreML:@PreEMS→@BAlg^@op$ sends an ^{enhanced measurable space} $(X,M,N)$ to the quotient ^{Boolean algebra}~$M/N$
and a ^{map of enhanced measurable spaces} $$f:(X,M,N)→(X',M',N')$$
to the morphism $$"PreML(f):M/N→M'/N'$$ in $@BAlg^@op$
given by the ^{map of Boolean algebras} $M'/N'→M/N$
induced by the map $f^*:M'→M$.
Furthermore, the functor~$"PreML$ descends to functors
$$"=StrictML:@StrictEMS→@BAlg^@op$$
and
$$"=WeakML:@EMS→@BAlg^@op.$$

\proclaim Lemma.
The above definition is correct.

\proof Proof.
By definition of a ^{σ-algebra} and ^{σ-ideal}, $N$ is an ideal of the ^{Boolean algebra}~$M$,
so the quotient $M/N$ is a ^{Boolean algebra}.
The map $f^*:M'→M$ preserves all Boolean operations
and it sends the ^{σ-ideal}~$N'$ to the ^{σ-ideal}~$N$ by definition of an ^{enhanced measurable space}.
Hence, $"PreML(f)^*:M'/N'→M/N$ is a ^{homomorphism of Boolean algebras}.
Composition is preserved because $(g∘f)^*=f^*∘g^*$.
Likewise, identities are preserved, which proves that we indeed have a functor~$"PreML$.
\ppar
^{Premaps of enhanced measurable spaces} that are ^{weakly equal almost everywhere}
by definition induce identical morphisms $M'/N'→M/N$.
This shows that the functor $$"PreML:@PreEMS→@BAlg^@op$$
descends to a functor $$"WeakML:@EMS→@BAlg^@op,$$
and hence also to a functor $$"StrictML:@StrictEMS→@BAlg^@op,$$
since $@EMS$ is a quotient of $@StrictEMS$.

\proclaim Definition.
The functor
$$"=ML:@CSLEMS→@MLoc.$$
(ML for measure-to-locale)
is the restriction of the functor $"WeakML$ to the corresponding subcategories.

\label\MLfactors
\proclaim Lemma.
The above definition is correct.

\proof Proof.
By \vstrictlocloc, the functor $"WeakML$ sends ^{strictly localizable enhanced measurable spaces} to ^{localizable Boolean algebras}.
It also sends a morphism $f:(X,M,N)→(X',M',N')$
of ^{strictly localizable enhanced measurable spaces}
to the homomorphism $"WeakML(f)^"op:M'/N'→M/N$ of ^{Boolean algebras} induced by the homomorphism $f^*:M'→M$.
The map $f^*$ preserves countable suprema, and hence so does $"WeakML(f)^"op$.
\ppar
To show that the map $"WeakML(f)^"op$ also preserves arbitrary suprema,
we reduce to the case of ^{σ-finite enhanced measurable spaces}.
By ^!{strictly localizable enhanced measurable spaces}, we can represent $(X,M,N)$
as a coproduct $∐_{i∈I}(X_i,M_i,N_i)$ of ^{σ-finite enhanced measurable spaces}
so that $M/N=∏_{i∈I}M_i/N_i$.
These spaces are ^{compact} by \vcompactsubspace.
The map $$"WeakML(f)^"op:M'/N'→M/N=∏_{i∈I}M_i/N_i$$
preserves arbitrary suprema if and only if its individual components
$M'/N'→M_i/N_i$
preserve suprema for all $i∈I$.
Thus, we may assume that $(X,M,N)$ is a ^{compact} ^{σ-finite enhanced measurable space}.
\ppar
By ^!{measurable image}, we have an element $X'_1∈M'$
such that $f^*(X'∖X'_1)∈N$ and for any $m'∈M'$ such that $m'⊂X'_1$ and $f^*m'∈N$ we have $m'∈N'$.
The ^{Boolean algebra} $M'/N'$ splits as the product of Boolean algebras $M'_1/N'_1$ and $M'_2/N'_2$,
where $(X'_1,M'_1,N'_1)$ and $(X'_2,M'_2,N'_2)$ are the ^{induced enhanced measurable spaces} of $X'_1$ and $X'_2=X'∖X'_1$.
By construction, the ^{homomorphism of Boolean algebras} $"WeakML(f)^"op$
is injective on $M'_1/N'_1$
and vanishes on $M'_2/N'_2$.
Thus, we can replace $(X',M',N')$ with $(X'_1,M'_1,N'_1)$ and assume $"WeakML(f)^"op$ to be injective.
\ppar
If $μ$ is a ^{faithful finite measure} on $(X,M,N)$, then $f_*μ=μ∘f^*$ is a ^{faithful finite measure} on $(X',M',N')$.
Thus, the map $"WeakML(f)^"op:M'/N'→M/N$ is an injective ^{homomorphism of Boolean algebras}
that preserves countable suprema and its domain and codomain admit faithful (finite) ^{continuous valuations}.
By \vbalgcountable\ the map $"WeakML(f)^"op$ preserves arbitrary suprema, which completes the proof.

\label\MLnotcomplete
\proclaim Remark.
Consider an uncountable set~$X$ with a ^{measure} $μ:2^X→[0,1]$ that vanishes on all singletons.
Such measures can be constructed in presence of sufficiently large cardinals, see Fremlin [\LC].
Set $N$ to the ^{σ-ideal} of subsets of~$X$ on which $μ$ vanishes.
The identity map $X→X$ yields a morphism $(X,2^X,N)→(X,2^X,\{∅\})$ in the category $@EMS$.
The domain is a ^{σ-finite enhanced measurable space} and the codomain is a ^{strictly localizable enhanced measurable space}.
The induced ^{homomorphism of Boolean algebras} $2^X→2^X/N$ does not preserve suprema
because all singletons in~$X$ map to~0, but their supremum in~$2^X$ is~$X$, which maps to the nonzero class of $X$ in $2^X/N$.
Thus, we really need ^{compactness} in the statement of \vMLfactors.

\subsection Equivalence of compact strictly localizable enhanced measurable spaces and measurable locales

\proclaim Definition.
The natural isomorphism~$ε$ is defined as a natural transformation
from the composition $$@MLoc\lto7{"Spec}@CSLEMS\lto7{"ML}@MLoc$$
to the identity functor on~$@MLoc$
that sends a ^{measurable locale}~$L$
to the isomorphism of ^{locales} $$ε_L:"ML("Spec(L))→L$$
whose associated ^{map of frames} $ε_L^*:L→"ML("Spec(L))$
sends $m∈L$ to the equivalence class of the ^{clopen subset} of $"Spec(L)$
corresponding to the open element of $"Ideal(L)$ given by the principal ideal of~$m$.

\proclaim Lemma.
This definition is correct.

\proof Proof.
\vhyperstoneanmeasurable\ shows that any measurable subset of $"Spec(L)$
has a unique presentation as the symmetric difference of a ^{clopen} subset and a ^{rare} subset of $"Spec(L)$.
Since negligible subsets of $"Spec(L)$ coincide with ^{rare} subsets,
the Boolean algebra $"ML("Spec(L))$ is canonically isomorphic to the ^{Boolean algebra} of ^{clopen} subsets of $"Spec(L)$,
equivalently, ^{clopen subsets} of $"Sp("Ideal(L))$,
i.e., $"COpen("Ideal(L))≅L$.
By construction, $ε_L^*$ is a ^{homomorphism of Boolean algebras}.
Thus, $ε_L$ is an isomorphism of ^{measurable locales}.
\ppar
To show the naturality of~$ε$, consider the following square for an arbitrary ^{map of measurable locales} $f:L→L'$:
$$\arrowsize15 \sqcd{
"ML("Spec(L))&\mapright{"ML("Spec(f))}&"ML("Spec(L'))\cr
\lmapdown{ε_L}&&\mapdown{ε_{L'}}\cr
L&\mapright{f}&L'.\cr
}$$
It suffices to show the commutativity of the induced diagram of inverse image maps:
$$\arrowsize15 \sqcd{
"ML("Spec(L))&\mapleft{"ML("Spec(f))^*}&"ML("Spec(L'))\cr
\lmapup{ε_L^*}&&\mapup{ε_{L'}^*}\cr
L&\mapleft{f^*}&L'.\cr
}$$
Pick an arbitrary element $m'∈L'$.
The element $ε_L^*(f^*(m'))$ is the equivalence class of the ^{clopen subset} of $"Spec(L)$ corresponding to~$f^*m'∈L$.
The element $"ML("Spec(f))^*(ε_{L'}^*(m'))$ is the equivalence class of the preimage under $"Spec(f)$ of the ^{clopen subset} of $"Spec(L')$ corresponding to~$m'$.
Using the equality $"Spec="TM∘"Sp∘"Ideal$ and the adjoint equivalence $"Ω⊣"Sp$ of \vHStonean,
we can equivalently describe this element as
the equivalence class of the ^{clopen subset} of $"Spec(L)$ corresponding to the ideal $"Idl(f^*)(m')$,
i.e., the ideal of~$L$ generated by the image of the principal ideal of~$m'$ under~$f^*$.
The latter ideal is precisely the principal ideal of~$f^*m'$.
Thus, the diagram commutes.

The following \vlifting\ is the technical heart of the paper.
It makes essential use of both ^{strict localizability} and ^{compactness}.

The map~$η$ in \vlifting\ is supplied by the von Neumann–Maharam lifting theorem.
See Maharam [\VN] for the ^{σ-finite} case and Fremlin [\MTiii, Corollary~341Q] for the ^{strictly localizable} case.
The ^{strict localizability} of $(X,M,N)$ is crucial at this point,
since Fremlin [\MTiii, Proposition~341M] shows that the existence of the map~$η$ implies ^{strict localizability}.

The construction of the inverse map~$ψ$ in \vlifting\ is due to
von Neumann [\MA], C.~Ionescu Tulcea [\SCE], Vesterstrøm–Wils [\PReal], Edgar [\MWS, Proposition~3.4], and Graf [\EPM, Theorem~1].
For an exposition, see Fremlin [\MTiii, Theorem~343B(iv)].
The ^{compactness} of $(X,M,N)$ is crucial for these results to be applicable,
since Fremlin [\CMS, §3, Theorem] and Rinkewitz [\CMMWS, Theorem~2.4] show that the existence of the map~$ψ$ implies ^{compactness}.

\label\lifting
\proclaim Proposition.
For any $(X,M,N)∈@CSLEMS$ there is an isomorphism
$$η_{(X,M,N)}:(X,M,N)→"Spec("ML(X,M,N)),$$
with the inverse map denoted by $ψ_{(X,M,N)}$.

\proof Proof.
Given $E∈M$, denote by $E^⋆$ the ^{clopen} subset of $"Spec("ML(X,M,N))$ corresponding to the element $[E]∈"ML(X,M,N)=M/N$.
Since $(X,M,N)$ is ^{strictly localizable}, by Fremlin [\MTiii, Corollary~341Q],
if $X∉N$, there is a map of sets
$$η:X→"Spec("ML(X,M,N))$$
such that $E⊕η^*(E^⋆)∈N$ for all $E∈M$ and $η^*F∈N$ for all ^{negligible subsets} $F⊂"Spec("ML(X,M,N))$.
(Choose an arbitrary positive faithful semifinite measure~$μ$ on $(X,M,N)$ to satisfy the conditions of the cited result.)
We recall the construction of~$η$: if $X∉N$, then by the von Neumann–Maharam lifting theorem (Fremlin [\MTiii, Theorem~341K])
the quotient homomorphism $M→M/N$ has a section $θ:M/N→M$, which is a ^{homomorphism of Boolean algebras}.
Then the point $η(x)∈"Spec("ML(X,M,N))$ is defined as the ^{homomorphism of Boolean algebras} $M/N→`Z/2$ that sends $m∈M/N$ to~$1$ if $x∈θ(m)$ and~$0$ otherwise.
(Recall that points of $"Spec(A)$ are precisely homomorphisms $A→`Z/2$.)
By construction, $θ([E])=η^*E^⋆$ for all $E∈M$, in particular, $η^*E^⋆∈M$ for all $E∈M$.
\ppar
Any measurable subset of $"Spec("ML(X,M,N))$ is by definition
the symmetric difference of $E^⋆$ for some $E∈M$ and a ^{meager subset},
and preimages of ^{meager subsets} under~$η$ belong to~$N$ (Fremlin [\MTiii, Proposition~341P(b)]),
so this implies that $η:(X,M,N)→"Spec("ML(X,M,N))$ is a ^{map of enhanced measurable spaces}.
\ppar
If $X∈N$, then $M/N=\{0\}$,
and hence $"Spec("ML(X,M,N))=∅$, so we take the (unique) map of sets $η:∅→"Spec("ML(X,M,N))$,
which yields a ^{map of enhanced measurable spaces}
$η:(X,M,N)→"Spec("ML(X,M,N))$ with $?pdom η=∅⊂X$.
It is precisely at this point that it is crucial to allow
the underlying maps of sets of ^{morphisms of enhanced measurable spaces} to have a ^{point-set domain} different from~$X$.
\ppar
Since $(X,M,N)$ is ^{compact}, by Fremlin [\MTiii, Theorem~343B(iv)],
we have a map of sets $$ψ_{(X,M,N)}:"Spec("ML(X,M,N))→X$$
such that $ψ^*E⊕E^⋆$ is negligible in the ^{enhanced measurable space} $"Spec("ML(X,M,N))$ for all $E∈M$.
The map~$ψ$ is constructed by setting $ψ(z)∈X$ to any element in the (nonempty) intersection of all sets $K∈K_X$ such that $z∈K^⋆$.
Here $K_X$ is the ^{compact class} that witnesses the ^{compactness} of $(X,M,N)$.
If there are no such~$K$, then $ψ(z)$ can be set to any point in~$X$.
The properties of~$ψ$ cited above imply that $$ψ:"Spec("ML(X,M,N))→(X,M,N)$$ is a ^{map of enhanced measurable spaces}
because $ψ^*n⊕n^⋆=ψ^*n⊕∅=ψ^*n∈N_{"Spec("ML(X,M,N))}$ for any $n∈N$
and $ψ^*m⊕m∈N_{"Spec("ML(X,M,N))}$ for any $m∈M$,
so $ψ^*m$ is the symmetric difference of $m$ and some negligible subset of $"Spec("ML(X,M,N))$ and therefore is measurable.
\ppar
For any $m∈M$ we have
$$η^*ψ^*m⊕m=η^*(m^⋆⊕n)⊕m=η^*m^⋆⊕m⊕η^*n=n'⊕η^*n∈N,$$
for some ^{negligible subset} $n⊂"Spec("ML(X,M,N))$ and some $n'∈N$.
Thus, $ψ∘η≈?id_{(X,M,N)}$.
Likewise, for any ^{measurable subset} of $"Spec("ML(X,M,N))$,
which can be assumed to be a ^{clopen} subset of the form~$E^⋆$,
the subset $ψ^*η^*E^⋆=ψ^*(E⊕n')=ψ^*E⊕ψ^*n'=E^⋆⊕n⊕ψ^*n'$ is a ^{negligible subset} of $"Spec("ML(X,M,N))$,
for some $n'∈N$ and some ^{negligible subset} $n⊂"Spec("ML(X,M,N))$.
Thus, $η∘ψ≈?id_{"Spec("ML(X,M,N))}$.
Hence, $η$ and $ψ$ are mutually inverse and both are isomorphisms.

\proclaim Proposition.
There is a natural isomorphism
$$η:?id_@CSLEMS→"Spec∘"ML$$
of functors $@CSLEMS→@CSLEMS$.

\proof Proof.
It remains to show that the maps~$η$ constructed in \vlifting\ are natural.
The square
$$\arrowsize19 \sqcd{
(X,M,N)&\mapright{h}&(X',M',N')\cr
\lmapdown{η_{(X,M,N)}}&&\mapdown{η_{(X',M',N')}}\cr
"Spec("ML(X,M,N))&\mapright{"Spec("ML(h))}&"Spec("ML(X',M',N'))\cr
}$$
commutes if and only if its image under~$"ML$ (depicted by the upper square below) does:
$$\arrowsize19 \sqcd{
"ML(X,M,N)&\mapright{"ML(h)}&"ML(X',M',N')\cr
\lmapdown{"ML(η_{(X,M,N)})}&&\mapdown{"ML(η_{(X',M',N')})}\cr
"ML("Spec("ML(X,M,N)))&\mapright{"ML("Spec("ML(h)))}&"ML("Spec("ML(X',M',N')))\cr
\lmapdown{ε_{"ML(X,M,N)}}&&\mapdown{ε_{"ML(X',M',N')}}\cr
"ML(X,M,N)&\mapright{"ML(h)}&"ML(X',M',N')\cr
}$$
By construction, we have $"ML(η_{(X,M,N)})=ε_{"ML(X,M,N)}^{-1}$ and likewise for the bottom right map.
Thus, the vertical compositions are identities, which completes the proof since the top and bottom maps must be the same.

\label\strictloc
\proclaim Remark.
The functor $"Spec$ lands inside $@CSLEMS$ and by \vcslemsdetiso\ the full subcategory $@CSLEMS⊂@LDEMS$ is closed under isomorphisms.
Thus, $(X,M,N)$ must be ^{compact} and ^{strictly localizable} for \vlifting\ to hold, assuming $(X,M,N)$ is ^{locally determined} (^!{locally determined}).
By Fremlin [\MTii, §216E], there are ^{localizable enhanced measurable spaces} that are not ^{strictly localizable},
so the category $@LEMS$ is not equivalent to $@MLoc$.
By \vMLnotcomplete, there is a morphism in $@LEMS$ whose image under $"WeakML$ does not preserve suprema.
If we discard such morphisms and pass to the nonfull subcategory $@=LEMSc$ with the same objects
and morphisms whose image under $"WeakML$ preserves suprema,
the functor $"WeakML$ restricts to a full functor $"ML:@LEMSc→@MLoc$.
If we now further perform a Gabriel–Zisman localization of $@LEMSc$ with respect to the morphisms~$f$ such that $"ML(f)$ is an isomorphism,
we obtain a new category $@LEMSc'$ together with the induced localization functor $@LEMSc'→@MLoc$.
The category $@CSLEMS$ is a full subcategory of $@LEMSc'$.
The functor $@LEMSc'→@MLoc$ is essentially surjective and full because its restriction to $@CSLEMS$ is essentially surjective and full.
It is faithful by definition of ^{weak equality almost everywhere}.
Thus, $@LEMSc'$ and $@MLoc$ are equivalent categories.
However, morphisms in $@CSLEMS$ (or in $@EMS$) admit a point-set description as equivalence classes of certain maps of sets,
whereas morphisms in $@LEMSc'$ do not (they are formal fractions of morphisms),
which is why we work with $@CSLEMS$.

\label\CSLEMSMLoc
\proclaim Theorem.
The functors
$$"ML:@CSLEMS→@MLoc$$
and
$$"Spec:@MLoc→@CSLEMS$$
together with natural isomorphisms
$$η:?id→"Spec∘"ML$$
and
$$ε:"ML∘"Spec→?id$$
form an adjoint equivalence of categories.

\proof Proof.
It remains to show that the exhibited equivalence is an adjoint equivalence.
The two triangle identities imply one another, so it suffices to establish just one of them.
To show that the composition
$$"ML(X,M,N)\lto{13}{"ML(η_{(X,M,N)})}"ML("Spec("ML(X,M,N)))\lto{13}{ε_{"ML(X,M,N)}}"ML(X,M,N)$$
equals identity it suffices to observe that
the map $ε_{"ML(X,M,N)}^*$ sends any $m∈"ML(X,M,N)=M/N$ to the equivalence class of the ^{clopen subset}~$m^⋆$ of $"Spec("ML(X,M,N))$.
Then the map $"ML(η_{(X,M,N)})$ sends this class to the equivalence class of $η_{(X,M,N)}^*m^⋆$,
which differs from~$m$ by an element of~$N$, and hence yields the same element of~$M/N$.

\section References

All bibliographic entries are equipped with hyperlinked back references.
An asterisk (*) after a back reference indicates that the corresponding entry is referenced after
the given numbered statement.

\medskip

\refs

\def\y#1{#1}
\yearkeytrue

\bib\MA
John von Neumann.
Einige Sätze über messbare Abbildungen.
Annals of Mathematics 33:3 (\y{1932}), 574--586.
\doi:10.2307/1968536.

\bib\RBA
Marshall~H.~Stone.
The theory of representations for Boolean algebras.
Transactions of the American Mathematical Society 40:1 (\y{1936}), 37--111.
\doi:10.1090/s0002-9947-1936-1501865-8.

\bib\BRGT
Marshall~H.~Stone.
Applications of the theory of Boolean rings to general topology.
Transactions of the American Mathematical Society 41:3 (\y{1937}), 375--481.
\doi:10.1090/s0002-9947-1937-1501905-7.

\bib\ACSBR
Marshall~H.~Stone.
Algebraic characterizations of special Boolean rings.
Fundamenta Mathematicae 29:1 (\y{1937}), 223--303.
\eudml:212943.
%\doi:10.4064/fm-29-1-223-303.

\bib\TRDL
Marshall~H.~Stone.
Topological representation of distributive lattices and Brouwerian logics.
Časopis pro pěstování matematiky a fysiky 67:1 (\y{1938}), 1--25.
\https://hdl.handle.net/10338.dmlcz/124080.

\bib\NR
Israel Gelfand.
Normierte Ringe.
Recueil Math\'ematique 9(51):1 (\y{1941}), 3--24.
\http://mi.mathnet.ru/msb6046.

\bib\OMCM
Paul~R.~Halmos, John von Neumann.
Operator methods in classical mechanics, II.
Annals of Mathematics 43:2 (\y{1942}), 332--350.
\doi:10.2307/1968872.

\bib\RepL
Lynn~H.~Loomis.
On the representation of σ-complete Boolean algebras.
Bulletin of the American Mathematical Society 53:8 (\y{1947}), 757--760.
\doi:10.1090/s0002-9904-1947-08866-2.

\bib\RepS
Roman Sikorski.
On the representation of Boolean algebras as fields of sets.
Fundamenta Mathematicae 35 (\y{1948}), 247--258.
\doi:10.4064/fm-35-1-247-258.

\bib\HS
Jacques Dixmier.
Sur certains espaces considérés par M.~H.~Stone.
Summa Brasiliensis Mathematicae 2 (\y{1951}), 151–182.
\https://dmitripavlov.org/scans/dixmier.pdf.

\bib\EqMS
Irving~E.~Segal.
Equivalences of measure spaces.
American Journal of Mathematics 73:2 (\y{1951}), 275--313.
\doi:10.2307/2372178.

\bib\DOA
Irving~E.~Segal.
Decompositions of operator algebras.  II: Multiplicity theory.
Memoirs of the American Mathematical Society 9 (\y{1951}).
\doi:10.1090/memo/0009.

\bib\CMeas
Edward Marczewski.
On compact measures.
Fundamenta Mathematicae 40:1 (\y{1953}), 113--124.
\doi:10.4064/fm-40-1-113-124.

\bib\CharW
Shôichirô Sakai.
A characterization of $W^*$-algebras.
Pacific Journal of Mathematics 6:4 (\y{1956}), 763--773.
\doi:10.2140/pjm.1956.6.763.

\bib\VN
Dorothy Maharam.
On a theorem of von Neumann.
Proceedings of the American Mathematical Society 9:6 (\y{1958}), 987--987.
\doi:10.1090/s0002-9939-1958-0105479-6.

\bib\SCE
Cassius Ionescu Tulcea.
Sur certains endomorphismes de ${\rm L}^∞_{\rm C}(Z,μ)$.
Comptes Rendus Hebdomadaires des Séances de l'Académie des Sciences 261 (\y{1965}), 4961--4963.
\https://gallica.bnf.fr/ark:/12148/bpt6k4027h/.

\bib\CSS
Johann Pfanzagl, W.~Pierlo.
Compact Systems of Sets.
Lecture Notes in Mathematics 16 (\y{1966}).
\doi:10.1007/BFb0078990.

\bib\PReal
Jørgen Vesterstrøm, Wilbert Wils.
On point realizations of $L^∞$-endomorphisms.
Mathematica Scandinavica 25 (\y{1969}), 178--180.
\doi:10.7146/math.scand.a-10954.

\bib\DT
Joan~W.~Negrepontis.
Duality in analysis from the point of view of triples.
Journal of Algebra 19:2 (\y{1971}), 228--253.
\doi:10.1016/0021-8693(71)90105-0.

\bib\CWstar
Shôichirô Sakai.
$C^*$-Algebras and $W^*$-Algebras.
Ergebnisse der Mathematik und ihrer Grenzgebiete 60 (\y{1971}).
\doi:10.1007/978-3-642-61993-9.

\bib\MWS
Gerald~A.~Edgar.
Measurable weak sections.
Illinois Journal of Mathematics 20:4 (\y{1976}), 630--646.
\doi:10.1215/ijm/1256049654.

\bib\EPM
Siegfried Graf.
Induced σ-homomorphisms and a parametrization of measurable sections via extremal preimage measures.
Mathematische Annalen 247:1 (\y{1980}), 67--80.
\doi:10.1007/bf01359867.

\bib\SS
Peter~T.~Johnstone.
Stone spaces.
Cambridge Studies in Advanced Mathematics 3 (\y{1982}).
\https://b-ok.cc/md5/CFB5295DDC01BCE9CAE1256930A62B5A.

\bib\PPtop
Peter~T.~Johnstone.
The point of pointless topology.
Bulletin of the American Mathematical Society 8:1 (\y{1983}), 41--53.
\doi:10.1090/s0273-0979-1983-15080-2.

\bib\ArtP
Peter~T.~Johnstone.
The art of pointless thinking: a student's guide to the category of locales.
Category theory at work.
Research and Exposition in Mathematics 18 (\y{1991}), 85--107.
\http://www.heldermann.de/R&E/RAE18/ctw06.pdf.

\bib\MIFHS
Michael Albert Wendt.
On measurably indexed families of Hilbert spaces.
Ph.D.~dissertation, Dalhousie University, \y{1992}.
\https://hdl.handle.net/10222/55327.

\bib\LC
David~H.~Fremlin.
Real-valued-measurable cardinals.
Set theory of the reals.
Israel Mathematical Conference Proceedings 6 (\y{1993}), 151--304.
\https://www1.essex.ac.uk/maths/people/fremlin/rvmc.pdf.

\bib\HCA
Francis Borceux.
Handbook of Categorical Algebra 3.
Categories of Sheaves.
Encyclopedia of Mathematics and its Applications 52 (\y{1994}).
\doi:10.1017/CBO9780511525872.

\bib\CDis
Michael~A.~Wendt.
The category of disintegration.
Cahiers de Topologie et Géométrie Différentielle Catégoriques 35:4 (\y{1994}), 291--308.
\eudml:91550.

\bib\MHS
Michael~A.~Wendt.
Measurable Hilbert sheaves.
Journal of the Australian Mathematical Society (Series A) 61:2 (\y{1996}), 189--215.
\doi:10.1017/s1446788700000197.

\bib\CBase
Michael~A.~Wendt.
Change of base for measure spaces.
Journal of Pure and Applied Algebra 128 (\y{1998}), 185--212.
\doi:10.1016/s0022-4049(97)00048-0.

\bib\CMS
David~H.~Fremlin.
Compact measure spaces.
Mathematika 46:2 (\y{1999}), 331--336.
\doi:10.1112/S0025579300007798.

\bib\MTii
David~H.~Fremlin.
Measure theory.  Volume~2.
Torres Fremlin, Colchester, \y{2001}.
\https://www1.essex.ac.uk/maths/people/fremlin/mt.htm.

\bib\CMMWS
Werner Rinkewitz.
Compact measures and measurable weak sections.
Mathematica Scandinavica 91:1 (\y{2002}), 150--160.
\doi:10.7146/math.scand.a-14383.

\bib\FLat
David~A.~Edwards.
On function lattices and commutative von Neumann algebras.
The Quarterly Journal of Mathematics 53:1 (\y{2002}), 19--29.
\doi:10.1093/qjmath/53.1.19.

\bib\ExtVal
Mauricio Alvarez-Manilla.
Extension of valuations on locally compact sober spaces.
Topology and its Applications 124 (\y{2002}), 397--433.
\doi:10.1016/s0166-8641(01)00249-8.

\bib\MTiii
David~H.~Fremlin.
Measure theory.  Volume~3.
Torres Fremlin, Colchester, \y{2002}.
\https://www1.essex.ac.uk/maths/people/fremlin/mt.htm.

\bib\MTiv
David~H.~Fremlin.
Measure theory.  Volume~4.
Torres Fremlin, Colchester, \y{2003}.
\https://www1.essex.ac.uk/maths/people/fremlin/mt.htm.

\bib\ObsSto[2005]
Hans~F.~de~Groote.
Observables I: Stone spectra.
\arXiv:math-ph/0509020v1.

\bib\ShMeas
Matthew Jackson.
A sheaf theoretic approach to measure theory.
Ph.D.~dissertation, University of Pittsburgh, \y{2006}.
ISBN: 978--0542--74794--6.
\https://d-scholarship.pitt.edu/7348/1/Matthew_Jackson_Thesis_2006.pdf.

\bib\GelfandDuality
Dmitri Pavlov.
Answer 3279 (revision~1) on MathOverflow.
October 29, \y{2009}.
\https://mathoverflow.net/revisions/3279/1.

\bib\PushPull
Dmitri Pavlov.
Answer 20820 (revision~1) on MathOverflow.
April 9, \y{2010}.
\https://mathoverflow.net/revisions/20820/1.

\bib\Duality
Dmitri Pavlov.
Question 23480 (revision~1) on MathOverflow.
May 4, \y{2010}.
\https://mathoverflow.net/revisions/23408/1.

\bib\MLoc
Dmitri Pavlov.
Answer 49542 (revision~1) on MathOverflow.
December 15, \y{2010}.
\https://mathoverflow.net/revisions/49542/1.

\bib\CQO
Hans~F.~de~Groote.
Classical and Quantum Observables.
Deep Beauty, Cambridge University Press (\y{2011}), 239--269.
\doi:10.1017/cbo9780511976971.

\bib\FLoc
Jorge Picado, Aleš Pultr.
Frames and Locales.  Topology without points.
Frontiers in Mathematics (\y{2012}), Birkhäuser.
\doi:10.1007/978-3-0348-0154-6.

\bib\LMSLGD[2014]
Simon Henry.
Localic metric spaces and the localic Gelfand duality.
Advances in Mathematics 294 (2016), 634--688.
\arXiv:1411.0898v1.
\doi:10.1016/j.aim.2016.03.006.